\documentclass[reqno,twoside,11pt,english]{amsart}
\usepackage{amsmath,amsfonts,amssymb,amsthm,epsfig}
\newcommand{\RNum}[1]{\uppercase\expandafter{\romannumeral #1\relax}}
\usepackage{comment}

\usepackage{color}
\usepackage[final,colorlinks,linkcolor=blue,anchorcolor=red,citecolor=blue]{hyperref}
\usepackage{graphicx}
\usepackage{titletoc,epsf}
\usepackage{amsmath,amsfonts,latexsym,amsthm,amsxtra,amssymb,bbm}
\allowdisplaybreaks
\usepackage{graphicx,float}
\usepackage{tikz}
\usepackage{ifthen}

\allowdisplaybreaks

\newtheorem{thm}{Theorem}[section]
\newtheorem{lem}[thm]{Lemma}
\newtheorem{op}[thm]{Open Problem}
\newtheorem{cor}[thm]{Corollary}
\newtheorem{prop}[thm]{Proposition}

\newtheorem{step}{Step}[section]

\newtheorem{cl}{Claim}[section]
\newtheorem{conc}[cl]{Conclusion}
\newtheorem{ca}{Case}
\newtheorem{sca}[section]{Subcase}
\newtheorem{scl}[section]{Subclaim}
\newtheorem{conj}[equation]{Conjecture}

\theoremstyle{definition}
\newtheorem{defn}[thm]{Definition}

\newtheorem{ques}[equation]{Question}
\newtheorem{rem}[thm]{Remark}
\newtheorem{exam}[thm]{Example}

\newcounter {own}
\def\theown {\thesection       .\arabic{own}}
\numberwithin{equation}{section}

\newenvironment{pf}[1][]{%
	\vskip 3mm
	\noindent
	\ifthenelse{\equal{#1}{}}%
	{{\slshape Proof. }}%
	{{\slshape #1.} }%
}%
{\qed\bigskip}

\newtheorem{Thm}{Theorem}

\newtheorem{Prop}[Thm]{Property}


\newcommand{\IR}{{\mathbb R}}

\newcommand{\IB}{{\mathbb B}}

\newcommand{\diam}{{\operatorname{diam}}}

\DeclareMathOperator{\CQH}{CQH}

\def\be{\begin{equation}}
	\def\ee{\end{equation}}

\newcommand{\bee}{\begin{enumerate}}
	\newcommand{\eee}{\end{enumerate}}

\newcommand{\blem}{\begin{lem}}
	\newcommand{\elem}{\end{lem}}
\newcommand{\bthm}{\begin{thm}}
	\newcommand{\ethm}{\end{thm}}
\newcommand{\bcor}{\begin{cor}}
	\newcommand{\ecor}{\end{cor}}
\newcommand{\beg}{\begin{exam}}
	\newcommand{\eeg}{\end{exam}}
\newcommand{\begs}{\begin{examples}}
	\newcommand{\eegs}{\end{examples}}
\newcommand{\bdefe}{\begin{defn}}
	\newcommand{\edefe}{\end{defn}}
\newcommand{\bprob}{\begin{prob}}
	\newcommand{\eprob}{\end{prob}}
\newcommand{\bques}{\begin{ques}}
	\newcommand{\eques}{\end{ques}}
\newcommand{\bei}{\begin{itemize}}
	\newcommand{\eei}{\end{itemize}}
\newcommand{\bcon}{\begin{conj}}
	\newcommand{\econ}{\end{conj}}
\newcommand{\bop}{\begin{op}}
	\newcommand{\eop}{\end{op}}

\newcommand{\bstep}{\begin{step}}
	\newcommand{\estep}{\end{step}}

\newcommand{\bca}{\begin{ca}}
	\newcommand{\eca}{\end{ca}}
\newcommand{\bsca}{\begin{sca}}
	\newcommand{\esca}{\end{sca}}

\newcommand{\bcl}{\begin{cl}}
	\newcommand{\ecl}{\end{cl}}

\newcommand{\bscl}{\begin{scl}}
	\newcommand{\escl}{\end{scl}}

\newcommand{\bcons}{\begin{conjs}}
	\newcommand{\econs}{\end{conjs}}
\newcommand{\bprop}{\begin{prop}}
	\newcommand{\eprop}{\end{prop}}
\newcommand{\br}{\begin{rem}}
	\newcommand{\er}{\end{rem}}
\newcommand{\brs}{\begin{rems}}
	\newcommand{\ers}{\end{rems}}
\newcommand{\bo}{\begin{obser}}
	\newcommand{\eo}{\end{obser}}
\newcommand{\bos}{\begin{obsers}}
	\newcommand{\eos}{\end{obsers}}
\newcommand{\bpf}{\begin{pf}}
	\newcommand{\epf}{\end{pf}}
\newcommand{\ba}{\begin{array}}
	\newcommand{\ea}{\end{array}}
\newcommand{\beq}{\begin{eqnarray}}
	\newcommand{\beqq}{\begin{eqnarray*}}
		\newcommand{\eeq}{\end{eqnarray}}
	\newcommand{\eeqq}{\end{eqnarray*}}

\newcommand{\ra}{\rightarrow}

\newcommand{\ds}{\displaystyle}

\newcounter{minutes}
\divide\time by 60
\newcounter{hours}
\multiply\time by 60 \addtocounter{minutes}{-\time}

\begin{document}

	\title[The dimension-free Gehring-Hayman inequality]{The dimension-free Gehring-Hayman inequality for quasigeodesics}
	
	\author[C.-Y. Guo, M. Huang and X. Wang]{Chang-Yu Guo, Manzi Huang and Xiantao Wang}

	\address[C.-Y. Guo]{Research Center for Mathematics and Interdisciplinary Sciences, Shandong University, 266237, Qingdao, P. R. China, and Department of Physics and Mathematics, University of Eastern Finland, 80101, Joensuu, Finland}
	\email{{\tt changyu.guo@sdu.edu.cn}}
	
	\address[M. Huang]{MOE-LCSM, School of Mathematics and Statistics, Hunan Normal University, Changsha, Hunan 410081, P. R. China} \email{mzhuang@hunnu.edu.cn}
	
	\address[X. Wang]{MOE-LCSM, School of Mathematics and Statistics, Hunan Normal University, Changsha, Hunan 410081, P. R. China} \email{xtwang@hunnu.edu.cn}
	
	\date{}

\begin{abstract}
	
	A well-known theorem of J. Heinonen and S. Rohde in 1993 states that if $D\subset \mathbb{R}^n$ is quasiconformally equivalently to an uniform domain, then the Gehring-Hayman inequality holds in $D$: quasihyperbolic geodesics in $D$ minimizes the Euclidean length among all curves in $D$ with the same end points, up to a universal dimension-dependent multiplicative constant.
	
	In this paper, we develop a new approach to strengthen the above result in the following three aspects: 1) obtain a dimension-free multiplicative constant in the Gehring-Hayman inequality; 2) relax the class of quasihyperbolic geodesics to more general quasigeodesics; 3) relax the quasiconformal equivalence to more general coarsely quasihyperbolic equivalence.
	
	As a byproduct of our general approach, we are able to prove that the above improved Gehring-Hayman inequality indeed holds in Banach spaces. This answers affirmatively an open problem raised by J. Heinonen and S. Rohde in 1993 and reformulated by J. V\"{a}is\"{a}l\"{a} in 2005. 

\end{abstract}


\maketitle

\tableofcontents

\section{Introduction}\label{sec-1}
\subsection{Background}
The unit disk $\mathbb{D}\subset \mathbb{C}$ plays a significant role in planar conformal analysis, as seen from the celebrated Riemannian mapping theorem, it is the distinguished conformal representative of all proper simply connected domains. In particular, many analysis on proper simply connected domains can be done via the analysis of conformal maps on unit disks \cite{Pommeranke-1992}. Replacing conformal mappings with the less rigid class of quasiconformal mappings, one would naturally expect that the class of quasiconformal images of unit ball inherits similar effects as planar disks in higher dimensions. Partially based on this point of view and on some ideas of F. John \cite{Jo}, O. Martio and J. Sarvas \cite{MS} introduced the important class of \emph{uniform domains}, which roughly speaking requires that the domain is quasiconvex and that each pair of points be joinable by a twisted double cone which is not too crooked nor too thin. Since then, uniform domains have received considerably large interests in the past decades and they have formed the ``correct" class of nice higher dimensional domains in analysis appearing in surprising many contexts \cite{Banuelos-1987,Gehring-1987,Geo,GP,GGGKN-2024,Guo-2015,Hajlasz-Koskela-2000,Herron-2021,Huang-Li-P-Wang-2015,Jo81,Martin-1985}.

Another more geometric view of the unit disk $\mathbb{D}\subset \mathbb{C}$ is to regard it as a complete Riemannian 2-manifold of negative curvature -1, when changing the underlying Euclidean metric to standard hyperbolic metric. Due to the Riemann mapping theorem, in terms of hyperbolic geometry, the unit disk becomes the distinguished hyperbolic isometric representative of all proper simply connected domains. In the celebrated work \cite{Gr-1}, M. Gromov introduced the notion of \emph{Gromov hyperbolic} metric spaces, which generalizes the metric properties of classical hyperbolic geometry and of trees. Gromov hyperbolicity is a large-scale property and turns out to be very useful in geometric group theory and metric geometry \cite{Bonk-Kleiner-2005,Bridson-H-1999}.

In an important and fundamental work \cite{BHK}, M. Bonk, J. Heinonen and P. Koskela have shown that the above two point of view, relating uniform spaces and Gromov hyperbolic spaces, are indeed conformally equivalent in the sense that there is a one-to-one conformal correspondence between the quasiisometric classes of certain Gromov hyperbolic spaces and the quasisimilarity classes of bounded locally compact uniform spaces. A basic starting point here is to explain uniformity through negative curvature with the aid of quasihyperbolic metric and the quasihyperbolization of a uniform space is Gromov hyperoblic. Based on this, they have developed an interesting uniformization theory, which provides conformal deformations of proper geodesic Gromov hyperbolic spaces into bounded uniform spaces. 
%

As was pointed out in \cite[Page 1]{BHK}, the main results of Bonk-Heinonen-Koskela relies heavily on the idea around \emph{Gehring-Hayman inequality}, which in particular depends on finite dimensional techniques such as Lebesgue measure and integration. Thus to extend the work of Bonk-Heinonen-Koskela to infinite dimensional spaces, one needs to develop new dimension-independent ideas and techniques.
As the Gehring-Hayman inequality plays a fundamental role in the theory, in this paper, we shall establish a dimension-free Gehring-Hayman inequality. It is thus useful to recall the definition here.

A metric space $(X,d)$ is said to satisfy the Gehring-Hayman inequality if there exists a constant $b>0$ such that for all $x,y\in X$ and for each quasihyperbolic geodesic $\vartheta$ with end points $x$ and $y$ and each curve $\gamma$ with the same end points, it holds
\begin{equation}\label{eq:def for Gehring-Hayman}
	\ell_d(\vartheta)\leq b\ell_d(\gamma).
\end{equation}
The Gehring-Hayman inequality was initially proved in \cite{Geh}, which states that the hyperbolic geodesic in a simply connected planar domain essentially minimizes the Euclidean length among all curves in the domain with the same end points, up to a universal multiplicative constant.
This inequality has later been generalized to quasihyperbolic geodesics in domains $\Omega\subset \IR^n$ that are quasiconformally equivalent to uniform domains by J. Heinonen and S. Rohde in \cite{HR}. This inequality, together with its conformally deformed variant, has found successful applications in the study of lengths of radii or boundary behavior of conformal and quasiconformal mappings \cite{BB-0,BB,BK,BKR,HeiN,Herron-2006,Koskela-Lammi-2013,Li}. Important extensions of this inequality to general (locally compact) metric spaces were obtained in \cite{BHK,Klm2014}. For our later reference, we state the main theorem of J. Heinonen and S. Rohde as follows.

\begin{Thm}[{\cite[Theorem 1.1]{HR})}]\label{thm:A}
	Let $D\subset \mathbb{R}^n$ be a domain that is $K$-quasiconformally equivalent to a $c$-uniform domain. If $\vartheta$ is a quasihyperbolic geodesic in $D$ and if $\gamma$ is any other arc in $D$ with the same end points as $\vartheta$, then
	\begin{equation*}
		\ell(\vartheta)\leq b\ell(\gamma),
	\end{equation*}
	where $b>0$ depends only on the dimension $n$, the uniformility constant $c$ and the quasiconformality coefficient $K$.
\end{Thm}

Recall that 
\begin{defn}\label{def:uniform domain}
	A domain $D$ in a metric space $X=(X,d)$ is called  {\it $c$-uniform}, $c\geq 1$, if each pair of points $z_{1},z_{2}$ in $D$ can be joined by a rectifiable curve $\gamma$ in $D$ satisfying
	\begin{enumerate}
		\item\label{con1} $\ds\min_{j=1,2}\{\ell (\gamma [z_j, z])\}\leq c\, d_D(z)$ for all $z\in \gamma$, and
		\item\label{con2} $\ell(\gamma)\leq c\,d(z_{1},z_{2})$,
	\end{enumerate}
	where $\ell(\gamma)$ denotes the arc-length of $\gamma$ with respect to the metric $d$,
	$\gamma[z_{j},z]$ the part of $\gamma$ between $z_{j}$ and $z$, and $d_D(z):=d(z,\partial D)$. In a $c$-uniform domain $D$, any curve $\gamma\subset D$, which satisfies conditions (1) and (2) above, is called a {\it $c$-uniform curve}.
\end{defn}

As straightforward applications of Theorem \ref{thm:A}, J. Heinonen and S. Rohde also obtained the following Gehring-Hayman inequality for image of quasihyperbolic geodesics.
\begin{Thm}[{\cite[Theorem 1.5]{HR})}]\label{thm:B}
	Suppose that $D\subsetneq \mathbb{R}^n$ is a $c$-uniform domain, and that $f:
	D\to D'$ is $M$-bi-Lipschitz in the quasihyperbolic metric. If
	$\vartheta$ is a quasihyperbolic geodesic in $D$, and if $\gamma$ is any other arc
	in $D$ with the same end points as $\vartheta$, then
	$$\ell(f(\vartheta))\leq b_1\;\ell(f(\gamma)),
	$$
	where $b_1$  depends only on the dimension $n$, the uniformility constant $c$ and the quasihyperbolicity coefficient $M$.
\end{Thm}

As far as we know, in all the previously mentioned works, the validity of Gehring-Hayman inequality depends (either explicitly or implicitly) on the finite dimensionality of the underlying space.  Our aim of this paper is to derive dimension-free estimates in both Theorems \ref{thm:A} and \ref{thm:B}. This is partially motivated by extending (part of) the theory of quasiconformal maps established in the seminal work of Heinonen and Koskela \cite{HeKo} to a larger class of non-Loewner metric spaces (such as infinite dimensional Banach spaces) and to solve an open problem of Bonk-Heinonen-Koskela \cite[Last paragraph on page 5]{BHK} regarding the relationship between (inner) uniformality and Gromov hyperbolicity in infinite-dimensional spaces,  where new methods giving dimension-free estimates are extremely valuable. Notice that the proof of Theorem \ref{thm:A} (and thus also Theorem \ref{thm:B}) relies crucially on \emph{the Whitney decomposition of domains} and \emph{the modulus of curve families}. Both these two techniques requires finite dimensionality of the underlying space $\mathbb{R}^n$ (in particular, the $n$-dimensional Lebesgue measure $\mathcal{L}^n$ on $\mathbb{R}^n$) and thus produces an extra dependence on the dimension $n$. Therefore, in order to obtain dimension-free estimate, a new approach is necessary.

\subsection{Main results}


From now on, we assume that $E$ and $E'$ are Banach spaces with dimensions no less than $2$ and $D$ is a proper subdomain in $E$.

	A homeomorphism $f:D\to D'$, where $D'\subset E'$, is called
	{\it $C$-coarsely $M$-quasihyperbolic}, abbreviated $(M,C)$-CQH,
	if for all $x$, $y\in D$, it holds
	$$\frac{1}{M}(k_D(x,y)-C)\leq k_{D'}(f(x),f(y))\leq M\;k_D(x,y)+C,$$
	where $k_D$ and  $k_{D'}$ denote the quasihyperbolic metrics in $D$ and $D'$, respectively (see Section \ref{sec-2} for the definition). In particular, an $(M,0)$-CQH mapping is said to be
	{\it $M$-bi-Lipschitz} in the quasihyperbolic metric.

\begin{defn}\label{def:quasigeodesic}
	A curve $\gamma\subset D$ is said to be a {\it
		$c$-quasigeodesic}, $c\geq 1$, if for all $x$, $y$ in $\gamma$, it holds
	$$\ell_{k}(\gamma[x,y])\leq c\, k_D(x,y),
	$$ where $\ell_{k}(\gamma[x,y])$ denotes the quasihyperbolic length of $\gamma[x,y]$ (see Section \ref{sec-2} for the definition).
\end{defn}

It is clear that, a $c$-quasigeodesic is a quasihyperbolic geodesic if and
only if $c=1$.

Our main result of this paper reads as follows.
\begin{thm}\label{thm:Heinonen-Rohde}
	Suppose that $D\subsetneq \mathbb{R}^n$ is homeomorphic to a $c$-uniform domain via an $(M,C)$-$\CQH$ mapping.  If $\vartheta$ is a $c_0$-quasigeodesic in $D$ and if $\gamma$ is any other arc in $D$ with the same end points as $\vartheta$, then
	\begin{equation*}
		\ell(\vartheta)\leq b\ell(\gamma),
	\end{equation*}
	where $b=b(c_0,c,M,C)$ depends only on the given parameters $c_0$, $c$, $M$ and $C$.
\end{thm}

Theorem \ref{thm:Heinonen-Rohde} improves on Theorem \ref{thm:A} in the following three aspects:
\begin{enumerate}
	
	\item The constant $b$ is independent of the dimension $n$.
	
	\item $\vartheta$ in Theorem \ref{thm:Heinonen-Rohde} lies in the larger class of quasigeodesics than quasihyperbolic geodesics.
	
	\item The class of CQH mappings is \emph{strictly larger} than the class of quasiconformal mappings. Indeed, by \cite[Theorem 6.4]{Vai3} and \cite[Theorem 4.14]{Vai4}, quasiconformal mappings in Euclidean spaces are quantitatively coarsely quasihyperbolic, but the converse is not true; see Example \ref{ex1} below for a counter-example.
\end{enumerate}

As an application of Theorem \ref{thm:Heinonen-Rohde}, we obtain the following improved Gehring-Hayman inequality for image of quasigeodesics under bi-Lipschitz quasihyperbolic homeomorphisms.
\begin{thm}\label{thm:Heinonen-Rohde-2}
	Suppose that $D\subsetneq \mathbb{R}^n$ is a $c$-uniform domain, and that $f:
	D\to D'$ is $M$-bi-Lipschitz in the quasihyperbolic metric. If
	$\vartheta$ is a $c_0$-quasigeodesic in $D$, and if $\gamma$ is any other arc
	in $D$ with the same end points as $\vartheta$, then
	$$\ell(f(\vartheta))\leq b_1\;\ell(f(\gamma)),
	$$
	where $b_1=b_1(c, c_0, M)$ depends only on the given parameters $c$, $c_0$ and $M$.
\end{thm}

Theorem \ref{thm:Heinonen-Rohde-2} improves on Theorem \ref{thm:B} in the following two aspects:
\begin{enumerate}	
	\item The constant $b_1$ is independent of the dimension $n$.
	
	\item $\vartheta$ is a quasigeodesic, which is weaker than quasihyperbolic geodesics unless $c_0=1$.
\end{enumerate}

Since our new approach is rather general, we are able to prove the dimension-free Gehring-Hayman inequality in the setting of Banach spaces, from which Theorem \ref{thm:Heinonen-Rohde} follows as a special case.
\begin{thm}\label{thm1.1}
	Let $E$ be a Banach space and $D\subsetneq E$ be homeomorphic to a $c$-uniform domain via an $(M,C)$-$\CQH$ mapping. If $\vartheta$ is a $c_0$-quasigeodesic in $D$, and if $L$
	is any other arc in $D$ with the same end points as $\vartheta$, then
	$$\ell(\vartheta)\leq b\;\ell(L),
	$$
	where $b=b(c, c_0, C, M)$ depends only on the given parameters $c$, $c_0$, $C$ and $M$.
\end{thm}
Theorem \ref{thm1.1} provides an affirmative answer to the following question: \emph{Does Theorem \ref{thm:A} hold for domains that are freely quasiconformally equivalent to an uniform domain in general Banach spaces?} This question was first asked by J. Heinonen and S. Rohde in \cite[Remark 2.25]{HR} and then by J. V\"ais\"al\"a in \cite[Open problem 5 in Section 6]{Vai9} with a more precise formulation using the class of coarsely quasihyperbolic homeomorphism (CQH). Theorem \ref{thm1.1} actually answers this question in a slightly stronger form, where quasihyperbolic geodesics are relaxed to the larger class of quasigeodesics.

As an application of Theorem \ref{thm1.1}, we obtain the following improved Gehring-Hayman inequality for image of quasigeodesics in the setting of Banach spaces, from which Theorem \ref{thm:Heinonen-Rohde-2} follows as a special case.
\begin{thm}\label{thm1.3}
	Let $E$ be a Banach space. Suppose that $D\subsetneq E$ is a $c$-uniform domain, and that $f:D\to D'$ is $M$-bi-Lipschitz in the quasihyperbolic metric. If
	$\vartheta$ is a $c_0$-quasigeodesic in $D$, and if $L$ is any other arc
	in $D$ with the same end points as $\vartheta$, then
	$$\ell(f(\vartheta))\leq b_1\;\ell(f(L)),
	$$
	where $b_1=b_1(c, c_0, M)$ depends only on the given parameters $c$, $c_0$ and $M$.
\end{thm}

In both Theorems \ref{thm1.1} and \ref{thm1.3}, one can replace the class of $c_0$-quasigeodesic to the larger class of $(c_0,\mu)$-quasigeodesic. Here a curve $\vartheta\subset D$ with end points $x,y$ is called a  $(c_0,\mu)$-quasigeodesic, $c_0\geq 1$ and $\mu\geq 0$, if for all $u,v\in \gamma$, it holds
$$\ell_{k}(\vartheta[u,v])\leq c_0\, k_D(u,v)+\mu\min\{|x-y|,1\}.$$
In this case, the constant $b$ or $b_1$ will depend additionally on the parameter $\mu$. To simplify our exposition, we shall not consider $(c_0,\mu)$-quasigeodesics in this paper, but leave the details to interested readers.  

Part of the techniques developed in this paper has been successfully applied in our following-up works \cite{Guo-Huang-Wang-2025-II,Guo-Huang-Wang-2025-III}. But in these works, in particular \cite{Guo-Huang-Wang-2025-III}, we address completely different questions and new innovations are necessary to overcome the difficulties arising there.

\subsection{Outline of the proof}
In this section, we briefly outline our approach and main innovations on the proof of Theorem \ref{thm1.1}. Let $\vartheta$ be a quasigeodesic in $D$ with end points $w$ and $\varpi$ and let $\gamma$ be any other curve in $D$ with the same end points as $\vartheta$. Our aim is to prove
\begin{equation}\label{eq:aim inequality}
	\ell(\vartheta)\leq b\;\ell(L)\tag{*}
\end{equation}
and we may assume $k_D(w,\varpi)>1$ as the other case follows from simple computations.

Our strategy for the proof of \eqref{eq:aim inequality} is to adapt the \emph{``contradiction-compactness argument"} to (locally noncompact) Banach spaces. This type of argument has appeared in many different contexts in geometric analysis or metric geometry, but usually in the finite dimensional setting. To run the contradiction argument, suppose for some large $b$, $\ell(\vartheta)\geq b\;\ell(L)$. Then we want to construct some ``good subcurves" $\vartheta_i$ of $\vartheta$ and $L_i$ of $L$, with the same end points, such that $\frac{\ell(\vartheta_i)}{\ell(L_i)}\gtrsim \mu^i\to \infty$ as $i\to \infty$. On the other hand, we can use the quasihyperbolic length of $L$ to give uniform upper bound for $\frac{\ell(\vartheta_i)}{\ell(L_i)}\leq M<\infty$. This contradiction would give the desired inequality.

To realize the above idea, we introduce the new class of \emph{$\lambda$-curves} (see Definition \ref{def:lambda curve}) and \emph{$\lambda$-pairs} (see Definition \ref{def:lambda pair and c condition}). $\lambda$-curves are essentially quasigeodesics with good control on the quasigeodesic coefficient. Comparing with quasihyperbolic geodesics, $\lambda$-curves are less sensitive with respect to small variations, which make the coarse analysis much more effective. A \textit{$\lambda$-pair} $(\gamma,L)$ consists of two curves in $D$ with the same end points, where $\gamma$ is a $\lambda$-curve, and we say that a $\lambda$-pair $(\gamma,L)$ satisfies \emph{$\mathfrak{C}$-condition} if $\ell(\gamma)\geq \mathfrak{C}\ell(L)$.

The most important step, to realize the previous strategy, is to prove the following ``compactness theorem" (see Theorem \ref{9-15-1}) for $\lambda$-pairs: \textit{if a $\lambda$-pair $(\gamma,L)$ satisfies \emph{$\mathfrak{C}$-condition} with $\mathfrak{C}$ sufficiently large, then there exists a $\lambda$-subpair $(\gamma_1,L_1)$ of $(\gamma,L)$, which satisfies $\mu\mathfrak{C}$-condition with a large $\mu\gg 1$, quantitatively.} More precisely, it says that if the length of a $\lambda$-curve $\gamma$ is too large compared with a fixed curve $L$ with the same end points, then we may find two points $w_1,w_2$ on $L$ such that for any $\lambda$-curve $\beta$ connecting these two points, we have $\ell(\beta)\geq \mu\mathfrak{C}\ell(L[w_1,w_2])$. Geometric insight suggests us that in certain small portion $L_1$ on $L$, the curve $\gamma$ has to oscillate a lot and so we may construct a $\lambda$-curve $\gamma_1$ with end points on $L$ (based on $\gamma$ and $L$) such that $\ell(\gamma_1)/\ell(L_1)$ becomes definitely larger than $\ell(\gamma_1)/\ell(L_1)$. However, to write down a rigorous proof  turns out to be not easy at all.

In our approach, the proof of Theorem \ref{9-15-1} is again based on a contraction argument. One key auxiliary result for the proof is Theorem \ref{22-11-16}, which says that the length of a class of $\lambda$-curve is uniformly bounded from above by the inner distance between its end points. This auxiliary result allows us to construct a sequence of points with prescribed metric controls (see Lemma \ref{9-14-1}), using basic properties of quasihyperbolic metrics, $\lambda$-curves, Gromov hyperbolicity and uniformality. A careful analysis on these points would lead to a contradiction.

The proof of Theorem \ref{22-11-16} is extremely technical and is again by a contradiction argument. The main innovation in the proof is to introduce the notion of \emph{six-tuples} (see Definition \ref{def:six-tuple}). A six-tuple $\Omega=[w_1,w_2,\gamma,\varpi,w_3,w_4]$ consists of five points $w_1,w_2,\varpi,w_3,w_4$ of special positions with controlled metric properties and a $\lambda$-curve $\gamma$ connecting $w_1$ and $w_2$. We shall analyze the basic properties of $\lambda$-curves via studying \emph{six-tuples which satisfies Property \ref{Property-A}.} The main ingredient is to construct iteratively new six-tuples based on older ones with prescribed metric properties. This step requires very delicate and careful estimates of quasihyperbolic distances using basic properties of $\lambda$-curves, Gromov hyperbolicity and $\CQH$-mappings. One feature of our analysis is that we reserve a large amount of constants (see Section \ref{subsec:notation} for the list) in the proof to compensate the lack of compactness in locally noncompact spaces and then derive a contradiction based on a careful selection of these constants.

\smallskip

\textbf{Structure.} The arrangement of this paper is as follows. In Section \ref{sec-2},
we introduce the list of reserved constants used in this paper, recall some known results, and develop a series of
elementary properties concerning the quasihyperbolic metric, quasigeodesics, $\lambda$-curves, solid arcs etc.
In Section \ref{sec-3}, we prove the compactness theorem for $\lambda$-pairs satisfying $\mathfrak{C}$-Condition, namely Theorem \ref{9-15-1} mentioned above, assuming the technical Theorem \ref{22-11-16}. Section \ref{sec-4} is devoted to the proofs of main results in this paper. In Section \ref{sec-5}, we give a complete proof of Theorem \ref{22-11-16}.
\smallskip

\textbf{Notations.} Throughout this paper, the norm of a vector $z\in E$ is written as $|z|$ and the distance between two points $z_1,z_2\in E$ is denoted by $|z_1-z_2|$ and $\sigma_D$ is the {\it inner distance} defined by
$$\sigma_D(z_1,z_2)=\inf \{\ell(\alpha):\; \alpha\subset D\;
\mbox{is a rectifiable curve joining}\; z_1\; \mbox{and}\; z_2 \}.$$ 
We write $[z_1, z_2]$ for the line segment with end points $z_1$ and $z_2$.

Let $D$ and $D'$ be proper domains in $E$ and $E'$, respectively, and let $f:D\to D'$ be a homeomorphism. For convenience, in what follows, the symbols $x, y, z,
\ldots$ represents points in $D$, while the primes $x',
y', z', \ldots$ denote their images in $D'$ under $f$. Similarly, we use $\alpha, \beta, \gamma, \ldots$ to denote arcs in $D$, while the primes $\alpha', \beta', \gamma', \ldots$ refer to their images in $D'$ under $f$. For a set $A$ in $D$, $A'$ often denotes its image in $D'$ under $f$. Throughout this paper, we call $\gamma$ an arc or a curve in $D$ if it is a continuous map from an interval to $D$ and we use $\ell_d(\gamma)$ and $\diam_d(\gamma)$ to denote the $d$-length and diameter of $\gamma$, respectively. When the content is clear, we omit the subscript $d$.

\section{Preliminaries}\label{sec-2}

\subsection{List of reserved constants}\label{subsec:notation}


In the following, we list the constants used in this paper.

\begin{enumerate}
	\item\label{ap-1}  $C$ and $M$: $\min\{C,M\}\geq 1$; 
	\smallskip
	
	\item\label{ap-2} $c_0$: $c_0> 1$;
	\smallskip
	
	\item\label{ap-3} $c_1$ and $h_1$:    $c_1=\frac{9}{4}M(M+1)$ and $h_1=C(\frac{9}{4}M+1)$ (see Lemma \ref{10-1});
	\smallskip
	
	\item\label{ap-4}  $c_2$:  $c_2=c_2(c, c_1, h_1)\geq e^{1+c_1+h_1}$ (see Theorem \ref{Theo-B} and Remark \ref{re-2.1});
	\smallskip

	
	\item\label{ap-6} $c_3$:  $c_3=c_3(c)\geq 8$ (see Theorem \ref{Theo} and Remark \ref{re-2.3});
	\smallskip
	
	\item\label{ap-7} $\upsilon_0$:  $\upsilon_0=\upsilon_0(c_0,\lambda,\delta)\geq 16$ (see Lemma \ref{Thm-23(1)} and Remark \ref{re-2.4});
	\smallskip
	
	\item\label{ap-8} $\mu_0$ and $\rho_0$:  $\mu_0=\mu_0(c,C,M,\delta)\geq 16$ and $\rho_0=\rho_0(c,C,M,\delta)>1$ (see Lemma \ref{Thm-7-2} and Remark \ref{re-2.5});
	\smallskip
	
	\item\label{ap-9} $\mu_1$ and $\mu_2$: $\mu_1=3\mu_0$ (see Lemma \ref{lem-8-12}) and $\mu_2=\frac{1}{4}+\mu_0+2\upsilon_0$ (see Lemma \ref{lem-2.4-2});
	\smallskip
	
	\item\label{ap-10} $\mu_3$, $\mu_4$ and $\mu_5$: $\mu_3=\frac{1}{8}c_0^2c_2c_3^3\mu_2e^{4+C+h_1+M}$, $\mu_4= 12c_2c_3(1+\mu_3)^2$ and $\mu_5= 2\mu_4e^{2\mu_4e^{8\mu_4}}$;
	\smallskip
	
	
	
	\item\label{ap-13}
	$\mathfrak{c}$: $\mathfrak{c}\geq 64\mu_1$ and $\mu_6=\mu_5^{\frac{1}{6}}$.
	\smallskip
	
	\item\label{ap-14}
	$\tau_1$, $\tau_2$ and $\tau_3$: $2^{11}\mu_1\leq 4\tau_1\leq \tau_2\leq \frac{1}{4}\tau_3$.
\end{enumerate}

\subsection{Coarse quasihyperbolicity and quasiconformality}
In the section, we shall construct an example which shows that coarsely quasihyperbolic homeomorphisms are in general not necessarily quasiconformal, even in the finite dimensional case.


\begin{exam}\label{ex1}
	For $s>0,$ let $\mathbb{B}(s)=\{x\in \mathbb{C}:\; |x|<s\}$, where $\mathbb{C}$ denotes the complex plane, and $O$ is the origin of $\mathbb{C}$. Set
	\begin{displaymath}
		f(x) = \left\{ \begin{array}{ll}
			g(|x|)\frac{x}{|x|}, & \textrm{$x\in \mathbb{B}(2)\backslash\{O\}$},\\
			O, & \textrm{$x=O$},
		\end{array} \right.
	\end{displaymath}
	where  $g: [0,2)\ra [0,2)$ is a homeomorphism defined as follows:
	\begin{displaymath}
		g(t) = \left\{ \begin{array}{ll}
			g_0(t), & \textrm{$t\in [0,\frac{1}{2})$},\\
			g_n(t), & \textrm{$t\in [1-\frac{1}{2^{n!}},1-\frac{1}{2^{(n+1)!}})$}\;\;\mbox{for}\;\; n\in \{1,2,...\},\\
			t, & \textrm{$t\in [1,2)$},
		\end{array} \right.
	\end{displaymath}
	$g_0(t)=\frac{3}{2}t,$
	and for each $n\in \{1,2,...\}$,
	$$g_n(t)=a_nt+b_n$$ with $$a_n=\frac{2^{(n+1)!}(2^{(n+2)!-(n+1)!}-1)}{2^{(n+2)!}(2^{(n+1)!-n!}-1)}\;\;\mbox{and}\;\;b_n=1-\frac{1}{2^{(n+1)!}}-a_n(1-\frac{1}{2^{n!}}).$$
\end{exam}
\noindent Then the following statements hold.
\begin{enumerate}
	\item\label{ex1-1}
	$f\colon \IB(2)\to \IB(2)$ is $(1,4)$-CQH.
	\item\label{ex1-2}
	$f\colon \IB(1)\to \IB(1)$ is not $\varphi$-FQC for any self-homeomorphism $\varphi\colon [0,+\infty)\to [0,+\infty)$. In particular, $f\colon \IB(2)\to \IB(2)$ is not quasiconformal by \cite[Remark 3.5]{Vai3}.
\end{enumerate}

Recall that a homeomorphism $f\colon \IB(1)\to \IB(1)$  is said to be $\varphi$-FQC, if for any $x,y\in \IB(2)$, it holds
\beq\label{4}
\varphi^{-1}(k_{\IB{(1)}}(x,y))\leq k_{\IB(1)}(f(x),f(y))\leq \varphi(k_{\IB{(1)}}(x,y)).
\eeq

\bpf
Firstly, we show that $f$ is $(1,4)$-CQH. For this, let $x$, $y\in \IB(2)$.
\begin{itemize}
	\item If $x$, $y\in \IB(1)$, then $f(x)$, $f(y)\in \IB(1)$
	and \cite[Lemma 2.2(2)]{Vai3} give
	\beqq\label{1}
	k_{\IB(2)}(f(x),f(y))\leq 2\leq k_{\IB(2)}(x,y)+2.
	\eeqq
	
	\item If $x$, $y\in \IB(2)\backslash\IB(1)$, then $f(x)=x$, $f(y)=y$, and so,
	\beqq\label{2}
	k_{\IB(2)}(f(x),f(y))= k_{\IB(2)}(x,y).
	\eeqq
	
	\item If $x\in \IB(1)$, $y\in \IB(2)\backslash\IB(1)$, then $f(x)\in \IB(1)$, $f(y)=y$. It follows that
	\beqq\label{3}
	k_{\IB(2)}(f(x),f(y))&\leq& k_{\IB(2)}(f(x),O)+k_{\IB(2)}(O,f(y)) \leq 4+ k_{\IB(2)}(x,y).
	\eeqq
\end{itemize}




Therefore, we conclude from the above estimates that for any $x$, $y\in \IB(2)$,
$$k_{\IB(2)}(f(x),f(y))\leq k_{\IB(2)}(x,y)+4.$$
Similarly, we may prove that for any $x$, $y\in \IB(2)$,
$$k_{\IB(2)}(x,y)\leq k_{\IB(2)}(f(x),f(y))+4.$$
This shows that $f$ is $(1,4)$-CQH.
\medskip

Suppose on the contrary that $f\colon \IB(1)\to \IB(1)$ is $\varphi$-FQC for a self-homeomorphism $\varphi:[0,+\infty)\to [0,+\infty)$. To get a contradiction, for each $n\in \mathbb{N}$, select $x_n, y_n\in \mathbb{S}(1-\frac{1}{2^{n!}})$ with
$$|x_n-y_n|=\frac{1}{2^{n!}}.$$
Then it follows from \cite[Lemma 2.2(2)]{Vai3} that
$$k_{\IB(1)}(x_n,y_n)\leq 1.$$

On the other hand, the definition of $f$ ensures that both $f(x_n)$ and $f(y_n)$ are on $\mathbb{S}(1-\frac{1}{2^{(n+1)!}})$, and
$$|f(x_n)-f(y_n)|\geq |x_n-y_n|.$$
This, together with \cite[Lemma 2.2(1)]{Vai3} and \eqref{4}, implies that
\beqq
\varphi(1)\geq k_{\IB(1)}(f(x_n),f(y_n))\geq  n\log2,
\eeqq
which is impossible if $n$ is sufficiently large. This completes the proof.
\epf

\subsection{Quasigeodesics, $\lambda$-curves and uniform domains}
We start this section by recalling the definition of quasihyperbolic metric, which was initially introduced by Gehring and Palka \cite{GP} for a domain in $\IR^n$ and then has been extensively studied in \cite{Geo}. Since we are concerning it in a Banach space, we shall follow the presentation from \cite{Vai3}.

The {\it quasihyperbolic length} of a rectifiable arc
$\gamma$ 
in a proper domain $D\subsetneq E$ is defined as
$$
\ell_{k}(\gamma):=\int_{\gamma}\frac{|dz|}{d_D(z)}.
$$
For any $z_1$, $z_2$ in $D$, the {\it quasihyperbolic distance}
$k_D(z_1,z_2)$ between $z_1$ and $z_2$ is set to be
$$k_D(z_1,z_2)=\inf_{\gamma}\{\ell_k(\gamma)\},
$$
where the infimum is taken over all rectifiable arcs $\gamma$
joining $z_1$ and $z_2$ in $D$.

An arc $\gamma$ from $z_1$ to $z_2$ is called a {\it quasihyperbolic geodesic} if
$\ell_k(\gamma)=k_D(z_1,z_2)$. Clearly, each subarc of a quasihyperbolic
geodesic is a quasihyperbolic geodesic. It is a well-known fact that a
quasihyperbolic geodesic between any two points in $E$ exists if the
dimension of $E$ is finite; see \cite[Lemma 1]{Geo}. This is not always 
true in infinite dimensional spaces; see \cite[Example 2.9]{Vai4}. In order to remedy this shortage, J. V\"ais\"al\"a \cite{Vai4}
introduced the class of \textit{quasigeodesics} in $E$; see Definition \ref{def:quasigeodesic}. Moreover, he proved in \cite[Theorem 3.3]{Vai4} that
for each pair of points $z_1$ and $z_2$ in $D$ and each $c>1$, there exists a
$c$-quasigeodesic in $D$ joining $z_1$ and $z_2$.

For any $z_1$, $z_2$ in $D$, let $\gamma\subset D$ be an arc with end points $z_1$ and $z_2$. Then we have the following elementary estimates (see for instance \cite[Section 2]{Vai3})

\beq\label{(2.1)}
\ell_{k}(\gamma)\geq
\log\Big(1+\frac{\ell(\gamma)}{\min\{d_D(z_1), d_D(z_2)\}}\Big)
\eeq
and
\beq\label{(2.2)}
\begin{aligned}
	k_{D}(z_1, z_2) &\geq  \log\Big(1+\frac{\sigma_D(z_1,z_2)}{\min\{d_D(z_1), d_D(z_2)\}}\Big)
	\\
	&\geq
	\log\Big(1+\frac{|z_1-z_2|}{\min\{d_D(z_1), d_D(z_2)\}}\Big)
	\geq
	\Big|\log \frac{d_D(z_2)}{d_D(z_1)}\Big|.
\end{aligned}
\eeq

\begin{lem}\label{lem-ne-1}
	Suppose that $x$ and $y$ are two points in $D$ with $k_D(x,y)\geq 1$.
	\begin{enumerate}
		\item\label{9-10-2}
		Then $\max\{d_D(x),d_D(y)\}\leq 2|x-y|$.
		\item\label{9-10-1}
		Let $\gamma$ be an arc in $D$ joining $x$ and $y$. Then for any $w\in \gamma$, $d_D(w)\leq 2\ell(\gamma).$
	\end{enumerate}
\end{lem}

\bpf (1). Assume that $$\max\{d_D(x),d_D(y)\}=d_D(x).$$
Suppose on the contrary that $d_D(x)> 2|x-y|.$ Then the segment $[x,y]\subset \mathbb{B}(x, \frac{1}{2}d_D(x))$ and thus for any $z\in [x,y]$, $$d_D(z)\geq d_D(x)-|x-z|>\frac{1}{2}d_D(x),$$
which implies $$k_D(x,y)\leq \int_{[x,y]}\frac{|dz|}{d_D(z)}<1.$$ This contradiction shows that the first statement holds.

(2). We may assume that $\gamma$ is rectifiable. Suppose on the contrary that there is some point $w_0\in \gamma$ such that $$\ell(\gamma)<\frac{1}{2}d_D(w_0).$$
Then for any $w\in \gamma$, we have
$$d_D(w)\geq d_D(w_0)-\ell(\gamma[w,w_0])>\frac{1}{2}d_D(w_0).$$
This gives $$k_D(x,y)\leq \int_{\gamma}\frac{|dw|}{d_D(w)}<1.$$
This contradiction proves the second statement, and hence, the proof is complete.
\epf

Next we discuss the coarse length of a curve introduced by J. V\"ais\"al\"a in \cite[Section 4]{Vai4}. Let $\gamma$ be an (open, closed or half open) arc in $D$. Write $\overline{x}=(x_0,\ldots,x_m)$ with
$m\geq 1$, whose coordinates form a finite sequence of $m$ successive points of $\gamma$.
For $h\geq 0$, we say that $\overline{x}$ is {\it $h$-coarse} if
$k_D(x_{j-1},x_j)\geq h$ for all $j\in\{ 1, \ldots, m\}$. Let $\Phi_k(\gamma,h)$
be the family of all $h$-coarse sequences of $\gamma$. Set
$$s_k(\overline{x}):=\sum^{m}_{j=1}k_D(x_{j-1},x_j)\;\; \mbox{and}\;\;
\ell_k(\gamma, h):=\sup \left\{s_k(\overline{x}):\; \overline{x}\in \Phi_k(\gamma,h)\right\}$$
with the convention that $\ell_k(\gamma, h)=0$ if
$\Phi_k(\gamma,h)=\emptyset$. Then the number $\ell_k(\gamma, h)$ is called the
{\it $h$-coarse quasihyperbolic length} of $\gamma$.

\bdefe[Solid arc] \label{def:h-solid} An arc
$\gamma$ in $D$ is called {\it $(c,h)$-solid} with $c\geq 1$ and
$h\geq 0$ if
$$\ell_k(\gamma[z_1,z_2], h)\leq c\;k_D(z_1,z_2)$$ for all pairs of points $z_1$ and $z_2$ in $\gamma$.
{We denote the class of all $(c,h)$-solid arcs in $D$ by $S^{c,h}$. The symbol $S^{c,h}_{xy}$ stands for those arcs in $S^{c,h}$ with end points $x,y$. }
\edefe

By \cite[Section 2.2]{Vai7} or \cite[Section 4]{Vai4}, we know that for any arc $\gamma \subset D$, $\ell_k(\gamma, 0)=\ell_k(\gamma)$.
Notice that $\gamma$ is a $(c,0)$-solid arc if and only if it is a $c$-quasigeodesic.

\bdefe[Short arc] \label{def:h-short}
A curve $\gamma$ in $D$ is said to be {\it $h$-short}
if for all $x,y\in\gamma$, $$\ell_{k}(\gamma[x,y])\leq  k_D(x,y)+h.$$
\edefe

It is clear from Definition \ref{def:h-short} that every subarc of an $h$-short arc is $h$-short, every $h$-short arc must be $h_1$-short for any $h_1\geq h$, and every quasihyperbolic geodesic is $h$-short with $h=0$.

\bdefe[$\lambda$-curve]\label{def:lambda curve}
For each pair of distinct points $x$, $y\in D$ and $\lambda>0$, we set $$c_{x,y}^\lambda:=\min\left\{1+\frac{\lambda}{2k_D(x,y)},\frac{9}{8}\right\}.$$ An arc $\gamma\subset D$ with end points $x$ and $y$ is called a $\lambda$-\textit{curve} if it is a $c_{x,y}^\lambda$-quasigeodesic in $D$. We use $\Lambda_{xy}^{\lambda}$ to stand for the set of all $\lambda$-curves in $D$ with end points $x$ and $y$ and use $\gamma_{xy}^\lambda$ to denote a general $\lambda$-curve in $\Lambda_{xy}^{\lambda}$. Throughout this paper, we always assume that $\lambda\in (0,\frac{1}{2}]$.
\edefe
By \cite[Theorem 3.3]{Vai4}, we know that for $x,y\in D$ and $\lambda>0$, if $x\not=y$, then $\Lambda_{xy}^{\lambda}\not=\emptyset$.
For $c_0\geq 1$,
	we use $Q^{c_0}$ to denote all $c_0$-quasigeodesics in $D$ and $Q^{c_0}_{xy}$ for $c_0$-quasigeodesics in $D$ with end points $x$ and $y$.
	
	The next result shows that all $\lambda$-curves possess the following shortness.
	
	\blem\label{lem-2.2}
	For any two distinct points $x\not=y\in D$, every $\lambda$-curve $\gamma\in \Lambda_{xy}^{\lambda}$ is $\lambda$-short.
	\elem
	\bpf Take any two points $u\not=v\in \gamma$. Then
	$$\ell_k(\gamma)\leq \big(1+\frac{\lambda}{2k_D(x,y)}\big)k_D(u,v)=k_D(u,v)+\frac{\lambda k_D(u,v)}{2k_D(x,y)}.$$
	Observe that $\gamma$ is a $\frac{9}{8}$-quasigeodesic and thus
	$$k_D(u,v)\leq \ell_k(\gamma)\leq \frac{9}{8}k_D(x,y).$$
	Combining the above two estimates gives the desired inequality:
	$$\ell_k(\gamma[u,v])\leq k_D(u,v)+\lambda.$$
	\epf
	
	The following lemma shows that every subcurve of a $\lambda$-curve is a $2\lambda$-curve.
	
	\blem\label{lem-2.3}
	Let $x\neq y$ be two distinct points in $D$ and fix a $\lambda$-curve $\gamma\in \Lambda_{xy}^\lambda$. Then for any $z_1\not=z_2\in\gamma$, the subarc $\gamma[z_1,z_2]$ of $\gamma$ belongs to $ \Lambda_{z_1z_2}^{2\lambda}$.
	\elem
	\bpf
	Let $z_1\not=z_2\in\gamma$ and observe that for any $u\not=v\in \gamma$,
	$$
	\ell_k(\gamma[u,v])\leq \min\left\{1+\frac{\lambda}{2k_D(x,y)},\frac{9}{8}\right\}k_D(u,v).
	$$
	To prove the lemma, it suffices to show that
	\[
	\min\left\{1+\frac{\lambda}{2k_D(x,y)},\frac{9}{8}\right\}\leq \min\left\{1+\frac{\lambda}{k_D(z_1,z_2)},\frac{9}{8}\right\}.
	\]
	
	As $\gamma$ is a $\frac{9}{8}$-quasigeodesic, we have
	$$
	k_D(z_1,z_2)\leq \ell_k(\gamma)\leq \frac{9}{8}k_D(x,y),
	$$
	which clearly implies
	$$1+\frac{\lambda}{2k_D(x,y)}\leq 1+\frac{\lambda}{k_D(z_1,z_2)}.$$ The proof is thus complete.
	\epf
	
	\blem\label{lem-5-14}
	Let $x\neq y$ be two distinct points in $D$. For any curve $\beta$ connecting $x$ and $y$ and any $\lambda$-curve $\gamma\in \Lambda_{xy}^\lambda$, if $\ell(\gamma)> \nu_1\ell(\beta)$ with $\nu_1\geq 18e^9$, then
	for any $w\in \beta$, $d_D(w)\leq 2\ell(\beta).$
	\elem
	\bpf By Lemma \ref{lem-ne-1}, it suffices to show that $k_D(x,y)> 1$. Suppose on the contrary that $k_D(x,y)\leq 1$. Since $\gamma$ is a $\frac{9}{8}$-quasigeodesic, it follows from (the proof of)
	\cite[Lemma 2.31]{Vai7} with $c=\frac{9}{8}$ and $C=1$ that $\gamma$ is a $ 18e^9$-uniform curve. This implies that
	$$\ell(\gamma)\leq  18e^9 |x-y|\leq  18e^9\ell(\beta),$$
	which clearly contradicts with the assumption that $\ell(\gamma)> \nu_1\ell(\beta)\geq 18e^9\ell(\beta)$.

	\epf
	
	From now on, we make the following convention: For two conditions, we say that Condition $\Phi$ quantitatively implies Condition $\Psi$ if Condition $\Phi$ implies Condition $\Psi$ and the data of Condition $\Psi$ depends only on that of Condition $\Phi$.
	\blem\label{10-1}
	$(1)$ For any two distinct points $x\not=y\in D$, each element in $\Lambda_{xy}^{\lambda}$ is a $(c_0,0)$-solid arc and belongs to $Q_{xy}^{c_0}$ with $c_0=\frac{9}{8}$.

	$(2)$ Suppose that $f:$ $D\to D'$ is $(M,C)$-$\CQH$.
	Then the image of any $\lambda$-curve under $f$ is a $(c_1,h_1)$-solid arc, where $c_1=\frac{9}{4}M(M+1)$ and $h_1=C\big(\frac{9}{4}M+1\big)$.
	\elem
	\begin{proof}	
		By \cite[Theorem 4.15]{Vai4}, we see that under an $(M, C)$-CQH mapping, the solidness of an arc implies the solidness of its image, quantitatively. The explicity formulation of $c_1$ and $h_1$ in terms of $M$ and $C$ follows from direct computations in \cite[Theorems 4.11 and 4.14]{Vai4}.
	\end{proof}
	
	By \cite[Theorem 6.22]{Vai4} and Lemma \ref{10-1}$(1)$, we have the following result concerning solid arcs and $\lambda$-curves in uniform domains.
	
	\begin{thm}[{\cite[Theorem 6.22]{Vai4}}]\label{Theo-B}  Suppose that $D$ is a $c$-uniform domain. For $z_1, z_2\in D$, let $\gamma\in S^{c_1,h_1}_{z_1z_2}\cup \Lambda_{z_1z_2}^\lambda$.
		Then for some $c_2=c_2(c, c_1, h_1)$, we have
		\bee
		\item\label{THB-1} for any $z\in \gamma$,
		$$\ds\min\{\diam(\gamma [z_1, z]),\;\; \diam(\gamma [z, z_2])\}\leq c_2\, d_D(z);$$
		\item\label{THB-2}  $\diam(\gamma)\leq c_2\max\big\{|z_1-z_2|, 2(e^{h_1}-1)\min\{d_D(z_1),d_D(z_2)\}\big\}$.
		\eee
	\end{thm}
	
	\br\label{re-2.1}
	In the rest of this paper, we assume that $$c_2=c_2(c, c_1, h_1)\geq e^{1+c_1+h_1}.$$
	\er

	We close this section with the following metric characterization of uniform domains in $E$, which we shall apply in our later proofs.
	
	\begin{thm}[{\cite[Theorem 6.16]{Vai4}}]\label{Theo}  For a domain $D\varsubsetneq E$, the following statements are
		equivalent: \bee
		\item $D$ is a $c$-uniform domain;
		\item
		For any $x,y\in D$, $k_D(x,y)\leq c_3\;
		\log\Big(1+\frac{|x-y|}{\ds\min\{d_D(x),d_D(y)\}}\Big)$,
		\eee
		where the constants $c$ and $c_3$ depend only on each other.
	\end{thm}
	
	\br\label{re-2.3}
	In the rest of the paper, we assume that $$c_3=c_3(c)\geq 8.$$
	\er

	\subsection{Gromov hyperbolicity and Rips spaces}
	There are many equivalent definitions of Gromov hyperbolicity and we shall mainly follow the one used in \cite{Vai10}.
	\begin{defn}\label{def:Gromov hyperbolicity}
		A domain $D\varsubsetneq E$ is called {\it $\delta$-Gromov hyperbolic}, $\delta>0$, if $(D,k_D)$ is $\delta$-Gromov hyperbolic.
		This means that for all $x, y, z, p\in D$, $$(x|z)_p\geq \min\{(x|y)_p,(y|z)_p\}-\delta,$$ where $(x|y)_p$ is the Gromov product defined by $$2(x|y)_p=k_D(p,x)+k_D(p,y)-k_D(x,y).$$
	\end{defn}
	
	It is well-known that every inner uniform domain is Gromov hyperbolic; see \cite[Remark 2.16]{Vai9'} or \cite[Theorem 1.11]{BHK}.
	Since Gromov hyperbolicity is preserved by CQH mappings \cite[Theorem 3.18]{Vai10}, it follows that every quasiconformal image of an inner uniform domain in $\IR^n$ is Gromov hyperbolic. In particular, this implies by the Riemann mapping theorem that every simply connected proper subdomain of the complex plane is Gromov hyperbolic.

	\blem\label{Thm-23(1)}
	Suppose that $D$ is $\delta$-Gromov hyperbolic, and $x\not=y\in D$. Then there exists a constant $\upsilon_0=\upsilon_0(c_0,\delta)>0$ such that
	\begin{itemize}
		\item for any $\gamma\in \Lambda_{xy}^{\lambda}$, any $\beta\in Q_{xy}^{c_0}$ and any $u_1\in\gamma$, there is $u_2\in \beta$ such that $k_D(u_1, u_2)\leq \upsilon_0$;
		
		\item for any $\gamma\in \Lambda_{xy}^{\lambda}$, any $\beta\in Q_{xy}^{c_0}$ and any $v_1\in \beta$, there is $v_2\in \gamma$ such that $k_D(v_1, v_2)\leq \upsilon_0$.
	\end{itemize}
	\elem
	\begin{proof}
		By \cite[Theorem 3.11]{Vai10}, we know that in a $\delta$-Gromov hyperbolic domain, if an $h$-short arc $\gamma$ and a $c$-quasigeodesic $\beta$ have the same end points, then the Hausdorff distance (with respect to the quasihyperbolic metric) between
		these two arcs is bounded by a constant which depends only on $c$, $h$ and $\delta$. Then the conclusions follow from the above fact by noticing Lemma \ref{lem-2.2}.
	\end{proof}
	
	\br\label{re-2.4}
	In the rest of the paper, we assume that $$\upsilon_0=\upsilon_0(c_0,\delta)\geq 16.$$
	\er
	
	Let $x_1,$ $x_2$ and $x_3$ be a triple of points in $D$. For each $i\in \{1,2,3\}$, let $\alpha_i$ denote an arc joining $x_{i}$ and $x_{i+1}$, where $x_4=x_1$.
	We write $\Delta=(\alpha_1,\alpha_2, \alpha_3)$ to represent a {\it triangle} in $D$ with sides $\alpha_1$, $\alpha_2$, $\alpha_3$ and the points $x_1,x_2,x_3$ are called the {\it vertices} of $\Delta$. A triangle $\Delta$ in $D$ is called {\it $h$-short} if all its three sides are $h$-short.
	
	The following result is a direct consequence of Lemma \ref{lem-2.2}.
	
	\blem\label{lem-8-1-3}
	Let $x_1,$ $x_2$ and $x_3$ be a triple of points in $D$. For any $\lambda$-curves $\gamma_1\in \Lambda_{x_1x_{2}}^\lambda,$ $\gamma_2\in \Lambda_{x_2x_{3}}^\lambda$ and $\gamma_3\in \Lambda_{x_3x_{1}}^\lambda$, the triangle $\Delta=(\gamma_1, \gamma_2, \gamma_3)$ is $\lambda$-short.
	\elem
	
	Next, we recall the definition of Rips space introduced in \cite[Section 2.26]{Vai10}. In what follows, we use the notation $k_D(w,\alpha)$ to denote the quasihyperbolic distance between the point $w\in D$ and the curve $\alpha$, i.e., $k_D(w,\alpha)=\inf\{k_D(w,v):\; v\in \alpha\}$.
	\bdefe[Rips space] \label{def12-1}
	Let $h\geq 0$ be a constant. If there is a constant $C>0$ such that for every $h$-short triangle $\Delta=(\alpha_1,\alpha_2, \alpha_3)$ in $D$ and for all $w\in \alpha_i$, $i\in \{1,2,3\}$, with $\alpha_{4}=\alpha_{1}$ and $\alpha_{5}=\alpha_{2}$, it holds
	$$k_D(w, \alpha_{i+1}\cup \alpha_{i+2})\leq C,$$
	then we call $(D, k_D)$ a {\it $(C,h)$-Rips space}.
	\edefe

	The following lemma gives uniform control on distance between ``close" points on short triangle in a Rips space.
	\begin{lem}\label{lem-12-21-1} Suppose that $(D, k_D)$ is a $(C,h)$-Rips space. Then for every $h$-short triangle $\Delta=(\alpha_1,\alpha_2, \alpha_3)$ in $D$, and for each $i\in \{1,2,3\}$, there is $w_i\in \alpha_i$ such that for any $i\not=j\in \{1,2,3\}$,
		$$k_D(w_i,w_j)\leq 3C.$$
	\end{lem}
	\bpf Assume that the vertices of $\Delta=(\alpha_1,\alpha_2, \alpha_3)$ are $x_1$, $x_2$ and $x_3$ such that for each $i\in \{1,2,3\}$, $\alpha_i$ connects $x_i$ and $x_{i+1}$, where $x_4=x_1$. We consider two cases.\medskip
	
	\noindent {\bf Case $1.$}  For any $z\in \alpha_1$, $k_D(z,\alpha_2)\leq C$.\medskip
	
	Since $\Delta$ is a triangle, there are $w_1\in \alpha_1\setminus \{x_1,x_2\}$ and $w_3\in \alpha_3\setminus \{x_1,x_3\}$ such that
	$$k_D(w_1,w_3)<C.$$ The assumption in this case ensures that there is $w_2\in \alpha_2\setminus \{x_2,x_3\}$ such that
	$$k_D(w_1,w_2)<2C.$$
	These yield that
	$$k_D(w_2,w_3)\leq k_D(w_1,w_2)+k_D(w_1,w_3) \leq 3C.$$
	Clearly, $w_1$, $w_2$ and $w_3$ are the desired points as claimed in the lemma.
	\medskip
	
	\noindent {\bf Case $2.$}  There is $w_0\in \alpha_1$ such that $k_D(w_0,\alpha_2)> C$.\medskip
	
	Since $k_D(x_2,\alpha_2)=0$, the assumption in this case implies that there are $u_1$ and $u_2\in \alpha_1[x_2,w_0]$ such that
	for any $u\in \alpha_1[x_2,u_1]$, $k_D(u,\alpha_2)\leq C$ and for any $w\in \alpha_1[u_1,u_2]\setminus \{u_1\}$, $k_D(w,\alpha_2)> C$. Then there is $w_1\in \alpha_1[u_1,u_2]\setminus \{u_1,w_0\}$ such that
	$k_D(w_1,\alpha_2)\leq \frac{4}{3}C,$ which implies that there is $w_2\in \alpha_2\setminus \{x_2,x_3\}$ such that
	\[
	k_D(w_1,w_2)\leq \frac{5}{3}C.
	\]

	Since $k_D(w_1,\alpha_2\cup \alpha_3)\leq C$, we see from the fact $k_D(w_1,\alpha_2)> C$ that $k_D(w_1,\alpha_3)\leq C,$ which implies that there is $w_3\in \alpha_3\setminus\{x_1,x_3\}$ such that
	$$k_D(w_1,w_3)\leq \frac{4}{3}C.$$
	By triangle inequality, we have
	$$k_D(w_2,w_3)\leq k_D(w_2,w_1)+k_D(w_1,w_3)\leq 3C.$$
	These show that the points $w_1$, $w_2$ and $w_3$ satisfy our requirements and hence the proof is complete.
	\epf
	
	The next result would be useful in the proofs of our later results.
	\begin{lem}\label{Thm-7-2}
		The following statements hold.
		\begin{enumerate}
			\item\label{Thm-7-2-1}
			If $D$ is $\delta$-Gromov hyperbolic, then $(D,k_D)$ is a $(\mu_0,\rho_0)$-Rips space for each $\rho_0>0$ with $\mu_0=\mu_0(\delta,\rho_0)$;
			\item\label{Thm-7-2-2}
			If $f:D\to D'$ is $(M,C)$-$\CQH$ and $D'$ is a $c$-uniform domain, then $(D, k_D)$ is a $(\mu_0,\rho_0)$-Rips space with $\mu_0=\mu_0(c,C,M)$ and $\rho_0=\rho_0(c,C,M)$. In particular, if $D'$ is $c$-uniform, then it is a $(\mu_0,\rho_0)$-Rips space with $\mu_0=\mu_0(c)$ and $\rho_0=\rho_0(c)$.
		\end{enumerate}
	\end{lem}
	\begin{proof}
		(1). This follows from \cite[Theorem 2.35]{Vai10}.

		(2). For this assertion, note first that by \cite[Theorem 2.12]{Vai9'}, a uniform domain is quantitatively Gromov hyperbolic in the quasihyperbolic metric, and so $(D',k_{D'})$ is $\delta'$-Gromov hyperbolic with $\delta'=\delta'(c)$. Next, by \cite[Theorem 3.18]{Vai10}, Gromov hyperbolicity is quantitatively preserved by CQH mappings. Consequently, $(D,k_D)$ is $\delta$-Gromov hyperbolic with $\delta=\delta(c,C,M)$. The claim follows then from the first assertion.
		
	\end{proof}
	
	\br\label{re-2.5}
	In the rest of the paper, we assume that $$\mu_0=\mu_0(c,C,M,\delta)\geq 16\;\;\mbox{and}\;\;\rho_0=\rho_0(c,C,M,\delta)>1.$$
	\er

	\begin{figure}[htbp]
		\begin{center}
			\begin{tikzpicture}[scale=0.8]
				\draw (-4,0) to [out=-25,in=205]  node [above left, pos=0, yshift=-0.1cm]{$v_3$} node [below, pos=0.15, xshift=-0.1cm]{$x_3$} coordinate[pos=0.131] (x3) node [below, pos=0.5, xshift=0cm]{$w_2$} coordinate[pos=0.5] (w2) node [above right, pos=1, yshift=-0.1cm]{$u_3$} (4,0);
				\filldraw  (x3) circle (0.04);
				\filldraw  (w2) circle (0.04);
				
				\draw (-4.1,-1) to [out=20,in=-110] node [left, pos=0, yshift=-0.1cm]{$v_1$} node [left, pos=0.5, yshift=0.1cm]{$w_3$} coordinate[pos=0.5] (w3) node [right, pos=0.8, xshift=0cm]{$x_1$} coordinate[pos=0.803] (x1) node [above, pos=1, xshift=0.1cm]{$u_1$} (0.6,4.5);
				\filldraw  (x1) circle (0.04);
				\filldraw  (w3) circle (0.04);
				
				\draw (4.1,-1) to [out=160,in=-70] node [right, pos=0, yshift=-0.1cm]{$v_2$} node [right, pos=0.5, yshift=0.1cm]{$w_1$} coordinate[pos=0.5] (w1) node [below, pos=0.15, xshift=-0.1cm]{$x_2$} coordinate[pos=0.164] (x2) node [above, pos=1, xshift=-0.1cm]{$u_2$} (-0.6,4.5);
				\filldraw  (x2) circle (0.04);
				\filldraw  (w1) circle (0.04);
				
			\end{tikzpicture}
		\end{center}
		\caption{Illustration of Lemma \ref{lem-8-12}} \label{fig-4-2}
	\end{figure}
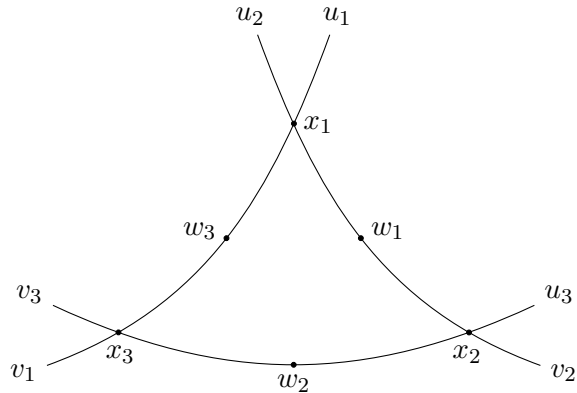
	
	Now, we prove two elementary estimates on distance between certain points in a Rips space.
	\blem\label{lem-8-12}
	Suppose that $(D,k_D)$ is a $(\mu_0,\rho_0)$-Rips space.
	For each $i\in \{1,2,3\}$, let $\{u_i,v_i\}$ denote a pair of points in $D$, $\gamma_i\in \Lambda_{u_iv_i}^\lambda$ and let $x_i\in \gamma_i\cap \gamma_{{i+1}}$, where $\gamma_4=\gamma_1$. Then for each $i\in \{1,2,3\}$, there is $w_i\in \gamma_{i+1}[x_i,x_{i+1}]$ such that for any $r\not=s\in \{1,2,3\}$,
	$$k_D(w_r,w_s)\leq \mu_1,$$ where $x_4=x_1$ and $\mu_1=3\mu_0$.
	\elem
	\begin{proof}
		This is a direct consequence of Lemmas \ref{lem-2.2}, \ref{lem-2.3}, \ref{lem-8-1-3} and \ref{lem-12-21-1}; see Figure \ref{fig-4-2}.
	\end{proof}
	
	\begin{lem}\label{lem-2.4-1}
		Suppose that $(D,k_D)$ is a $(\mu_0,\rho_0)$-Rips space.
		For each $i\in \{1,2\}$, let $\{u_i,v_i\}$ be a pair of points in $D$, $\gamma_i\in \Lambda_{u_iv_i}^\lambda$ and let $x\in \gamma_1\cap \gamma_2$, $y\in \gamma_1[x,v_1]$ and $z\in \gamma_2[x,v_2]$. If $k_D(y,z)\leq \tau$, then for any $w\in\gamma_1[x,y]$, there exists $u\in\gamma_2[x,z]$ such that
		\[
		k_D(w,u)\leq \mu_0+2\lambda+\tau.
		\]
	\end{lem}
	\bpf Let $w\in\gamma_1[x,y]$ and fix any $\gamma_3\in \Lambda_{yz}^\lambda$. Since $(D,k_D)$ is a $(\mu_0,\rho_0)$-Rips space, it follows from Lemmas \ref{lem-2.2} and  \ref{lem-2.3} that there exists
	$v\in\gamma_2[x,z]\cup \gamma_3$ satisfying
	\[
	k_D(w,v)\leq \lambda+\mu_0.
	\]
	If $v\in\gamma_2[x,z]$, then we simply take $u=v$. 
	If $v\in\gamma_3$, then we take $u=z$. It follows from Lemma \ref{lem-2.2} and the previous estimate that
	\begin{eqnarray*}k_D(w,u)&\leq&   k_D(w,v)+\ell_k(\gamma_3)\leq k_D(w,v)+k_D(y,z)+\lambda
		\\ \nonumber&\leq& \mu_0+2\lambda+\tau.
	\end{eqnarray*}
	In either cases, we see that the conclusion of lemma holds. The proof is thus complete.
	\epf
	
	The next result concerns a similar estimate in a Gromov hyperbolic space.
	\begin{lem}\label{lem-2.4-2}
		Suppose that $D$ is $\delta$-Gromov hyperbolic, and $x,y,z\in D$ are distinct points. Then there exists a constant $\mu_2=\mu_2(\delta,c_0)$ such that for each $\gamma_1\in Q_{xy}^{c_0}$, $\gamma_2\in Q_{yz}^{c_0}$, $\gamma_3\in Q_{xz}^{c_0}$ and each  $u\in \gamma_1$, there is $v\in \gamma_{2}\cup \gamma_{3}$ satisfying the inequality
		$$k_D(u,v)\leq \mu_2.$$
		Indeed, we may take $\mu_2=\frac{1}{4}+\mu_0+2\upsilon_0$ with $\mu_0$ and $\upsilon_0$ given by Lemma \ref{Thm-7-2} and Lemma \ref{Thm-23(1)}, respectively.
	\end{lem}
	\bpf Fix $u\in \gamma_{1}$, $\beta_1\in \Lambda_{xy}^\lambda$, $\beta_2\in \Lambda_{yz}^\lambda$ and $\beta_3\in \Lambda_{xz}^\lambda$ with $\lambda=\frac{1}{4}$. Then Lemma \ref{Thm-23(1)} shows that there is $u_1\in\beta_1$
	such that
	\be\label{23-10-11-0}k_D(u,u_1)\leq \upsilon_0,\ee
	where $v_0=v_0(c_0,\delta)$ is a uniform constant depending only on $c_0$ and $\delta$.
	
	By Lemma \ref{Thm-7-2}\eqref{Thm-7-2-1}, $(D,k_D)$ is a $(\mu_0,\rho_0)$-Rips space with $\mu_0=\mu_0(\delta)$ and $\rho_0=\rho_0(\delta)$. Then Lemma \ref{lem-2.2} implies that there is $u_2\in \beta_2\cup\beta_3$
	such that \be\label{23-10-11-1}k_D(u_1,u_2)\leq \frac{1}{4}+\mu_0.\ee
	Without loss of generality, we may assume $u_2\in\beta_3$. Applying Lemma \ref{Thm-23(1)} once again, we find $v\in \gamma_{3}$
	such that $$k_D(u_2,v)\leq \upsilon_0,$$
	which, together with \eqref{23-10-11-0} and \eqref{23-10-11-1}, yields
	$$k_D(u,v)\leq k_D(u,u_1)+k_D(u_1,u_2)+k_D(u_2,v)\leq \frac{1}{4}+\mu_0+2\upsilon_0.$$ This completes the proof.
	\epf
	
	We now establish the following technical lemma, which compares the length of quasigeodesics with one common end point in Gromov hyperbolic spaces.
	\begin{lem}\label{lem-2.4}
		Suppose that $D$ is $\delta$-Gromov hyperbolic, $x_1$, $x_2$ and $x_3$ are points in $D$, and that $\mu_2=\mu_2(\delta,c_0)$ is the constant given by Lemma \ref{lem-2.4-2}. Suppose further that there is a constant $\tau\geq 0$ such that $k_D(x_2,x_3)\leq\tau$. If $$\min\left\{k_D(x_1,x_2),k_D(x_1,x_3)\right\}> 1,$$ then for any $\vartheta\in Q_{x_1x_2}^{c_0}$ and $\gamma\in Q_{x_1x_3}^{c_0}$, we have
		$$\max\left\{\ell(\vartheta),\ell(\gamma)\right\}\leq 4e^{2c_0^2(2\tau+3\mu_2)+\frac{4}{3}\mu_2}\min\left\{\ell(\vartheta),\ell(\gamma)\right\}.$$
	\end{lem}
	
	\bpf
	Fix $\vartheta\in Q_{x_1x_2}^{c_0}$ and $\gamma\in Q_{x_1x_3}^{c_0}$.  If $\ell(\vartheta)=\ell(\gamma)$, then there is nothing to prove. Otherwise, we may assume without loss of generality that $$\ell(\vartheta)>\ell(\gamma).$$
	Under this assumption, to prove the lemma, we only need to verify the following inequality:
	\beq\label{7-8-10}
	\ell(\vartheta)\leq 4e^{2c_0^2(2\tau+3\mu_2)+\frac{4}{3}\mu_2}\ell(\gamma),
	\eeq
	where $\mu_2=\mu_2(\delta,c_0)$ is the constant given by Lemma \ref{lem-2.4-2}.
	
	In the following, we divide the discussions into two cases: $k_D(x_1,x_2)\leq 2c_0(2\tau+3\mu_2)$ and $k_D(x_1,x_2)> 2c_0(2\tau+3\mu_2)$. Since $k_D(x_1,x_3)> 1,$ it follows from Lemma \ref{lem-ne-1}\eqref{9-10-1} that
	\be\label{eq-l-ne-7}
	d_D(x_1)\leq 2\ell(\gamma).
	\ee
	
	In the first case, \eqref{(2.1)} gives
	$$2c_0^2(2\tau+3\mu_2)\geq  c_0k_D(x_1,x_2)\geq \ell_{k}(\vartheta)  \geq \log\Big(1+\frac{\ell(\vartheta)}{\min\{d_D(x_1),d_D(x_2)\}}\Big).$$
	This, together with \eqref{eq-l-ne-7}, implies
	\beq\label{7-8-11}\quad\quad\quad
	\ell(\vartheta)\leq e^{2c_0^2(2\tau+3\mu_2)}\min\{d_D(x_1),d_D(x_2)\}\leq 2e^{2c_0^2(2\tau+3\mu_2)}\ell(\gamma).
	\eeq
	
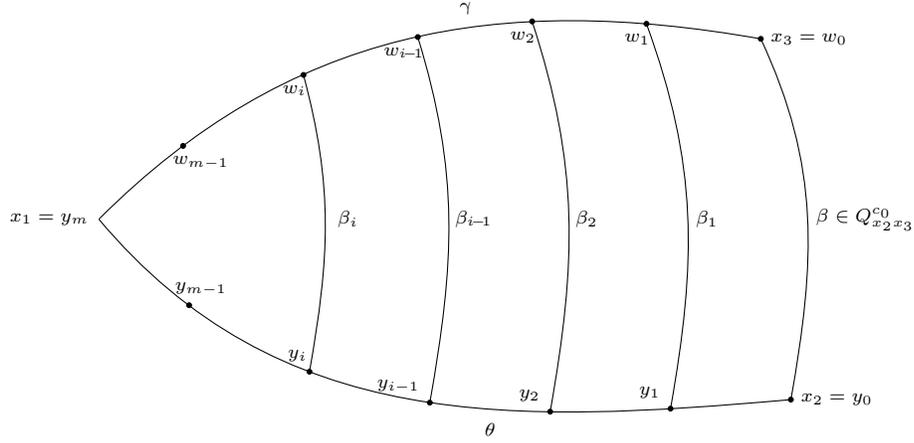
\begin{figure}[htbp]
		\begin{center}
			\begin{tikzpicture}[scale=0.8]
				\draw (-6,3)node[left] {\tiny $x_1=y_m$} to [out=-50,in=185] (5.5,0);
				
				\filldraw  (-4.5,1.57)node[above] {\tiny \;\;\;\;$y_{m-1}$} circle (0.04);
				\filldraw  (-2.5,0.465)node[above] {\tiny $y_{i}\;\;\;\;$} circle (0.04);
				\filldraw  (-0.5,-0.05)node[above left] {\tiny $y_{i-1}$} circle (0.04);
				\filldraw  (1.5,-0.2)node[above left] {\tiny $y_{2}$} circle (0.04);
				\filldraw  (3.5,-0.15)node[above left] {\tiny $y_{1}$} circle (0.04);
				\filldraw  (5.5,0)node[right] {\tiny $x_2=y_{0}$} circle (0.04);
				
				\draw (-6,3) to [out=45,in=170] (5,6);
				\filldraw  (-4.6,4.22)node[below] {\tiny \;\;\;\;\;\;$w_{m-1}$} circle (0.04);
				\filldraw  (-2.6,5.4)node[below] {\tiny $w_{i}\;\;\;$} circle (0.04);
				\filldraw  (-0.7,6.03)node[below] {\tiny $w_{i\!-\!1}\;\;\;\;\,$} circle (0.04);
				\filldraw  (1.2,6.29)node[below] {\tiny $w_{2}\;\;\;$} circle (0.04);
				\filldraw  (3.1,6.25)node[below] {\tiny $w_{1}\;\;\;$} circle (0.04);
				\filldraw  (5,6)node[right] {\tiny $x_3=w_{0}$} circle (0.04);
				
				\draw (-2.5,0.465) to [out=80,in=-75] (-2.6,5.4);
				\draw (-0.5,-0.05) to [out=80,in=-73] (-0.7,6.03);
				\draw (1.5,-0.2) to [out=80,in=-72] (1.2,6.29);
				\draw (3.5,-0.15) to [out=80,in=-70] (3.1,6.25);
				\draw (5.5,0) to [out=80,in=-65] (5,6);
				
				\node[rotate=10] at(0.1,6.5) {\tiny $\gamma$};
				\node[rotate=-5] at(0.5,-0.5) {\tiny $\theta$};
				
				\node[rotate=0,right] at(-2.2,3) {\tiny $\beta_{i}$};
				\node[rotate=0,right] at(-0.25,3) {\tiny $\beta_{{i\!-\!1}}$};
				\node[rotate=0,right] at(1.75,3) {\tiny $\beta_{2}$};
				\node[rotate=0,right] at(3.75,3) {\tiny $\beta_1$};
				\node[rotate=0,right] at(5.75,3) {\tiny $\beta\in Q_{x_2x_3}^{c_0}$};
				
			\end{tikzpicture}
		\end{center}
		\caption{Illustration for the proof of second case in Lemma \ref{lem-2.4}} \label{fig-12-3}
	\end{figure}

	In the second case, that is, $k_D(x_1,x_2)> 2c_0(2\tau+3\mu_2)$, we need some preparation. First, we determine a partition of the arc $\vartheta$. For this, we denote $x_2$ by $x_0$.
	Clearly, there exists a finite sequence of points $\{y_i\}_{i=0}^m$ in $\vartheta$ with $m\geq 3$ such that for each $i\in \{1,\cdots,m-1\}$, $y_i\in \vartheta[y_{i-1},x_1],$
	\beq\label{eq-l-ne-8}\quad\quad
	k_D(y_{i-1},y_i)=c_0(2\tau+3\mu_2)\;\;\mbox{and}\;\;0<k_D(y_{m-1},y_m)\leq c_0(2\tau+3\mu_2),
	\eeq
	where $y_m=x_1$; see Figure \ref{fig-12-3}.
	
	By \eqref{(2.1)}, for each $i\in \{1,\cdots,m\}$, we have
	$$\log\Big(1+\frac{\ell(\vartheta[y_{i-1},y_i])}{\min\{d_D(y_{i-1}),d_D(y_{i})\}}\Big)\leq \ell_{k}(\vartheta[y_{i-1},y_i])\leq c_0k_D(y_{i-1},y_i)\leq c_0^2(2\tau+3\mu_2),$$
	from which it follows
	\beq\label{eq-l-ne-12}
	\ell(\vartheta[y_{i-1},y_i]))\leq e^{c_0^2(2\tau+3\mu_2)}\min\{d_D(y_{i-1}),d_D(y_{i})\}.
	\eeq

	Next, we are going to find a partition of $\gamma$ based on $\{y_i\}_{i=1}^{m-1}$. For this, we fix $\beta\in Q_{x_2x_3}^{c_0}$ and select $p\in \beta$ such that
	\be\label{11-16-1}
	k_D(y_1,\beta)\leq k_D(y_1,p)<k_D(y_1,\beta)+\frac{7}{8}\tau,
	\ee
	where $k_D(y_1,\beta)$ is the quasihyperbolic distance from $y_1$ to $\beta$.
	
	From the assumption of $k_D(x_2,x_3)<\tau$, we note
	$$\ell_{k}(\beta)\leq c_0k_D(x_2,x_3)<c_0\tau,$$
	and thus
	$$k_D(y_1,p)\geq k_D(y_1,x_0)-k_D(x_0,p)\geq k_D(y_1,x_0)-c_0\tau.$$
	Then we infer from \eqref{eq-l-ne-8} and \eqref{11-16-1} that
	\beq\label{7-8-1}
	k_D(y_1,\beta)\geq k_D(y_1,p)-\frac{7}{8}\tau\geq c_0(3\mu_2+\tau)-\frac{7}{8}\tau>\mu_2.
	\eeq
	Since $D$ is $\delta$-Gromov hyperbolic, we know from Lemma \ref{lem-2.4-2} that there is $q_1\in \gamma\cup \beta$ such that
	$k_D(y_1,q_1)\leq \mu_2.$ By \eqref{7-8-1}, we see that $q_1\in \gamma$. This implies that
	\beq\label{7-8-2}
	k_D(y_1,\gamma)\leq k_D(y_1,q_1) \leq \mu_2.
	\eeq
	
	Now, we are ready to determine a point in $\gamma$ based on $y_1$.
	Denote $x_3$ by $w_0$, and let $w_1\not=w_0$ be a point in $\gamma$ (see Figure \ref{fig-12-3}) such that $$k_D(y_1,w_1)<\frac{4}{3}k_D(y_1,\gamma).$$
	Then
	it follows from \eqref{7-8-2} that
	\beq\label{7-8-3}
	k_D(y_1,w_1)\leq \frac{4}{3}\mu_2,
	\eeq which is the point we need.
	
	To find the next point in $\gamma$ based on $y_2$, we fix $\beta_1\in Q_{w_1y_1}^{c_0}$ and derive a lower bound for $k_D(y_2, \beta_1)$ as follows. Since
	$\ell_k(\beta_1)\leq c_0 k_D(y_1, w_1)$, we infer from \eqref{eq-l-ne-8} and  \eqref{7-8-3} that
	\beq\label{7-8-4}
	k_D(y_2, \beta_1)\geq k_D(y_1,y_2)-\ell_k(\beta_1)\geq c_0\big(\frac{5}{3}\mu_2+2\tau\big)>\mu_2.
	\eeq
	Once again, Lemma \ref{lem-2.4-2} implies that there is $q_2\in \gamma[x_1,w_1]\cup \beta_1$ such that
	\[
	k_D(y_2,q_2)\leq \mu_2.
	\]
	Then \eqref{7-8-4} ensures that $q_2\in \gamma[x_1,w_1]$, and thus,
	\beq\label{7-8-5}
	k_D(y_2,\gamma[x_1,w_1])\leq k_D(y_2,q_2)\leq \mu_2.
	\eeq
	
	Let $w_2$ be a point in $\gamma[x_1,w_1]\setminus\{w_1\}$ (see Figure \ref{fig-12-3}) such that $$k_D(y_2,w_2)<\frac{4}{3}k_D(y_2,\gamma[x_1,w_1]).$$ Then
	it follows from \eqref{7-8-5} that
	\beq\label{7-8-6}
	k_D(y_2,w_2)\leq \frac{4}{3}\mu_2,
	\eeq 
	which is the point we want.
	
	By repeating this procedure, we may find a finite sequence of points $\{w_i\}_{i=1}^{m-1}$ in $\gamma$ (see Figure \ref{fig-12-3}) such that for each $i\in\{1,\cdots,m-1\}$,
	\beq\label{eq-l-ne-9}
	w_i\in \gamma[x_1,w_{i-1}]\setminus\{w_{i-1}\}\;\;\mbox{and}\;\;k_D(y_i,w_i)\leq \frac{4}{3}\mu_2.
	\eeq
	
	In the following, we prove that for each $i\in\{1,\cdots,m-1\}$,
	\beq\label{eq-l-ne-11}
	\ell(\vartheta[y_{i-1},y_i])\leq 2e^{c_0^2(2\tau+3\mu_2)+\frac{4}{3}\mu_2}\ell(\gamma[w_{i-1},w_i]).
	\eeq
	
	First, \eqref{eq-l-ne-9}, together with \eqref{(2.2)}, ensures that for each $i\in\{1,\cdots,m-1\}$, we have
	$$\Big|\log \frac{d_D(y_i)}{d_D(w_i)}\Big|\leq k_D(y_i,w_i)\leq  \frac{4}{3}\mu_2,$$
	which gives
	\beq\label{eq-l-ne-10}
	d_D(y_i)\leq e^{\frac{4}{3}\mu_2} d_D(w_i).
	\eeq
	
	Since $$k_D(w_{i-1},w_i) \geq k_D(y_{i-1},y_i)-k_D(y_{i-1},w_{i-1})-k_D(y_i,w_i),$$
	we deduce from \eqref{eq-l-ne-8}, \eqref{eq-l-ne-9} and the assumption $k_D(x_0,w_0)=k_D(x_2,x_3)<\tau$ that
	\beqq\label{Ge}
	k_D(w_{i-1},w_i)
	\geq
	\left\{
	\begin{array}{ll}
		(2c_0-1)\tau+\frac{9c_0-4}{3}\mu_2, & {\rm when}\,\; i=1,
		\\
		\;\;\;\;2c_0\tau+\frac{9c_0-8}{3}\mu_2, & {\rm when}\,\; 2\leq i\leq m-1.
	\end{array}
	\right.
	\eeqq
	Since our choice of $c_0$ and $\mu_2$ in Section \ref{subsec:notation} implies that
	$$k_D(w_{i-1},w_i)>1,$$
	we know from Lemma \ref{lem-ne-1}\eqref{9-10-1} that
	$$d_D(w_i)\leq 2\ell(\gamma[w_{i-1},w_i]),$$
	and so, it follows from \eqref{eq-l-ne-12} and \eqref{eq-l-ne-10} that for each $i\in\{1,\cdots,m-1\}$,
	\beqq\label{7-9-1}
	\ell(\vartheta[y_{i-1},y_i]) \leq  2e^{c_0^2(2\tau+3\mu_2)+\frac{4}{3}\mu_2}\ell(\gamma[w_{i-1},w_i]),
	\eeqq which proves \eqref{eq-l-ne-11}. Moreover, by \eqref{eq-l-ne-7} and \eqref{eq-l-ne-12}, we have
	\beq\label{7-8-13}
	\ell(\vartheta[y_{m-1},y_m]) \leq  2e^{c_0^2(2\tau+3\mu_2)}\ell(\gamma).
	\eeq
	
	Now, we are ready to show \eqref{7-8-10}. Since
	$$\ell(\vartheta) =\sum_{i=1}^{m-1}\ell(\vartheta[y_{i-1},y_{i}])+\ell(\vartheta[y_{m-1},y_{m}]),
	$$
	it follows from \eqref{eq-l-ne-11} and \eqref{7-8-13} that
	\beq\label{7-8-14}
	\ell(\vartheta)\leq 4e^{c_0^2(2\tau+3\mu_2)+\frac{4}{3}\mu_2}\ell(\gamma).
	\eeq
	
	In either cases, we conclude from  \eqref{7-8-11} and \eqref{7-8-14} that the desired inequality \eqref{7-8-10} holds, and hence, the proof is complete.
	\epf
	
	Based on Lemma \ref{lem-2.4}, we now establish the following useful length comparison theorem for quasigeodesics with the same end points in Gromov hyperbolic spaces.
	\bthm\label{24-4-11-12}
	Suppose that $D$ is $\delta$-Gromov hyperbolic, $x_1$ and $x_2$ are points in $D$. Then for any $\vartheta$, $\gamma\in Q_{x_1x_2}^{c_0}$, it holds
	$$\max\{\ell(\vartheta),\ell(\gamma)\}\leq 4e^{8c_0^2\mu_2}\min\{\ell(\vartheta),\ell(\gamma)\}.$$
	\ethm
	\bpf Without loss of generality, we may assume that $\min\{d_D(x_1),d_D(x_2)\}=d_D(x_1)$ and $\max\{\ell(\vartheta),\ell(\gamma)\}=\ell(\vartheta)$. Note that if $k_D(x_1,x_2)\leq 1$, then we may get from \eqref{(2.1)} that
	$$\log\Big(1+\frac{\ell(\vartheta)}{d_D(x_1)}\Big)\leq \ell_k(\vartheta)\leq c_0k_D(x_1,x_2)\leq c_0,$$
	which gives
	\be\label{25H-01-1} \ell(\vartheta)\leq (e^{c_0}-1)d_D(x_1).\ee
	On the other hand, we know from \eqref{(2.2)} that
	$$\log \Big(1+\frac{|x_1-x_2|}{d_D(x_1)}\Big)\leq k_D(x_1,x_2)\leq 1,$$
	which implies that
	$$|x_1-x_2|\leq (e-1)d_D(x_1).$$
	Clearly, for each $x\in [x_1,x_2]$, it holds
	$$d_D(x)\geq \min\{d_D(x_1),d_D(x_2)\}-\frac{1}{2}|x_1-x_2|= d_D(x_1)-\frac{1}{2}|x_1-x_2|\geq \frac{3-e}{2}d_D(x_1),$$
	and so we have
	\beqq
	k_D(x_1,x_2)\leq \int_{x\in [x_1,x_2]}\frac{|dx|}{d_D(x)}\leq \frac{2|x_1-x_2|}{(3-e)d_D(x_1)}.
	\eeqq
	Combining this with \eqref{25H-01-1} gives
	$$\frac{\ell(\vartheta)}{e^{c_0}d_D(x_1)}\leq \ell_k(\vartheta)\leq c_0k_D(x_1,x_2)\leq \frac{2c_0|x_1-x_2|}{(3-e)d_D(x_1)},$$
	from which it follows
	$$\ell(\vartheta)\leq \frac{2c_0e^{c_0}}{3-e}|x_1-x_2|\leq \frac{2c_0e^{c_0}}{3-e}\ell(\gamma).$$
	
	Thus,  we only need to consider the case $k_D(x,y)> 1.$
	In this case, the result follows directly from Lemmas \ref{lem-2.4} (applied with $x_2=x_3$).
	\epf

	As an application of Theorem \ref{24-4-11-12}, we prove a technical lemma.
	\blem\label{24-4-22-3}
	Suppose that $D$ is $\delta$-Gromov hyperbolic. For $u, v\in D$, let $\beta\subset D$ be a curve connecting $u$ and $v$, and suppose that for all $w\in \beta$,
	\be\label{24-4-22-1}
	d_D(w)\leq 2\ell(\beta).
	\ee
	Let $u_1\in \beta$ and $v_1\in \beta[u_1,v]$. For any $\lambda$-curve $\gamma_1\in \Lambda_{u_1v_1}^{\lambda}$, if $\ell(\gamma_1)\geq \mu_6^{\frac{1}{2}}\ell(\beta)$, then there is $w_0\in \beta[u_1,v_1]$
	such that
	\begin{enumerate}
		\item
		for any $\lambda$-curve $\gamma_2\in \Lambda_{u_1w_0}^{\lambda}$,
		$$\mu_6^{\frac{1}{4}}\ell(\beta)\leq \ell(\gamma_2)\leq \mu_6^{\frac{1}{2}}\ell(\beta).$$
		\item
		for any $z\in \beta[w_0,v_1]\setminus\{w_0\}$, and for any $\lambda$-curve $\gamma_3\in \Lambda_{u_1z}^{\lambda}$ satisfying
		$$\ell(\gamma_3)> \mu_6^{\frac{1}{4}}\ell(\beta).$$
	\end{enumerate}
	\elem
	
	\bpf
	For each $w\in\beta[u_1,v_1]$,  let $\Lambda_w$ be a subclass of $\Lambda_{u_1w}^{\lambda}$ consisting of $\gamma_{w}\in \Lambda_{u_1w}^{\lambda}$ with
	$$\ell(\gamma_w)\leq 2\inf_{\gamma\in \Lambda_{u_1w}^{\lambda}}\{\ell(\gamma)\}.$$
	Then there exist $w_1\in\beta[u_1,v_1]$ and $\gamma_{w_1}\in \Lambda_{w_1}$ such that
	\begin{enumerate}
		\item[$(i)$]
		$\ell(\gamma_{w_1})\leq  2\mu_6^{\frac{1}{4}}\ell(\beta);$
		\item[$(ii)$]
		for any $w\in \beta[w_1,v_1]\setminus\{w_1\}$, there exists $\gamma_w\in \Lambda_w$ with
		$\ell(\gamma_{w})>2\mu_6^{\frac{1}{4}}\ell(\beta).$
	\end{enumerate}
	Indeed, if not, then by Theorem \ref{24-4-11-12}, $\ell(\gamma_1)\leq 8e^{11\mu_2}\mu_6^{\frac{1}{4}}\ell(\beta)<\mu_6^{\frac{1}{2}}\ell(\beta)$, contradicting with our assumption.

	Now select $w_0\in \beta[w_1,v_1]$ such that
	\be\label{24-4-27-1}
	0<k_D(w_0,w_1)<1
	\ee
	and let $\gamma_{w_0}\in \Lambda_{w_0}$ be given by $(ii)$.
	Then we deduce from Lemma \ref{10-1}$(1)$, \eqref{24-4-22-1} and  $(ii)$ that
	\[
	\begin{aligned}
		k_D(u_1,w_0) &\geq  \frac{8}{9}\ell_{k}(\gamma_{w_0})\geq  \frac{8}{9}\log\Big(1+\frac{\ell(\gamma_{w_0})}{\min\{d_D(u_1),d_D(w_0)\}}\Big)\\
		& \geq  \frac{8}{9}\log\Big(1+\frac{\ell(\gamma_{w_0})}{2\ell(\beta)}\Big)\geq \frac{8}{9}\log\Big(1+\mu_6^{\frac{1}{4}}\Big)>1.
	\end{aligned}
	\]
	Moreover, by \eqref{24-4-27-1}, we obtain that
	\beqq
	k_D(u_1,w_1)\geq k_D(u_1,w_0)-k_D(w_1,w_0)\geq \frac{8}{9}\log\Big(1+\mu_6^{\frac{1}{4}}\Big)-1>1.
	\eeqq
	The above inequalities yield \be\label{24-4-27-3}\min\{k_D(u_1,w_0), k_D(u_1,w_1)\}>1.\ee
	Since $D$ is $\delta$-Gromov hyperbolic, we infer from \eqref{24-4-27-1},
	\eqref{24-4-27-3}, $(i)$ and Lemma \ref{lem-2.4} (replacing $x_1$, $x_2$ and $x_3$ with $u_1$, $w_1$ and $w_0$, respectively) that
	$$\mu_6^{\frac{1}{4}}\ell(\beta)\leq \frac{1}{2}\ell(\gamma_{w_0})\leq \ell(\gamma_2)\leq e^{9\mu_2}\ell(\gamma_{w_1})< \mu_6^{\frac{1}{2}}\ell(\beta),$$
	which proves (1).
	
	For the second assertion, simply note that by  $(ii)$, it holds
	$$ \mu_6^{\frac{1}{4}}\ell(\beta)<\frac{1}{2}\ell(\gamma_z)\leq \ell(\gamma_3).$$
	The proof is thus complete.
	\epf

	\section{Compactness of $\lambda$-pairs with $\mathfrak{C}$-Condition}\label{sec-3}
	\subsection{Main result}\label{sec-3-1}
	Throughout this section, we make the following standing assumptions.
	\vspace{8pt}
	
	\noindent \textbf{Standing Assumptions.} {\it We assume that $f:$ $D\to D'$ is $(M,C)$-CQH and $D'$ is a $c$-uniform domain. Let $\omega\not=\varpi\in D$, and let $\gamma\in\Lambda_{\omega\varpi}^{\lambda}$ be a fixed element with $\lambda=\frac{1}{4}$. }
	\vspace{8pt}

	For the statement of our main result of this section, we need to recall the following definitions.
	\begin{defn}[$\lambda$-pair and $\mathfrak{C}$-Condition]\label{def:lambda pair and c condition}
		A $\lambda$-curve $\gamma\in \Lambda_{uv}^{\lambda}$ and a curve $\beta$ in $D$ is called a {\it $\lambda$-pair}, written as $(\gamma,\beta)$, if the end points of $\beta$ are $u$ and $v$. For $u_1\not=v_1\in \beta$, the pair $(\gamma[u_1v_1],\beta[u_1,v_1])$ is called a {\it $\lambda$-subpair} of $(\gamma,\beta)$. A $\lambda$-pair $(\gamma,\beta)$ is said to satisfy {\it $\mathfrak{C}$-Condition} if there is a constant $\mathfrak{C}>0$ such that $\ell(\gamma)\geq \mathfrak{C}\ell(\beta)$.
	\end{defn}

	
	The aim of this section is to prove the following two theorems, which play key roles for the discussions in the next section. In the statement below, $\mu_5$ is defined in Section \ref{subsec:notation}.
	
	\bthm\label{9-15-1}
	Suppose that $\mathfrak{l}$ is an arc in $D$ such that $(\gamma,\mathfrak{l})$ is a $\lambda$-pair and satisfies $\mathfrak{C}$-Condition $($see Figure \ref{fig-9-22}$)$. If $\mathfrak{C}\geq \frac{4}{5}\mu_5$, then there exist points $\omega_1\not=\varpi_1\in\mathfrak{l}$ such that for any $\lambda$-curve $\gamma[\omega_1\varpi_1]\in \Lambda_{\omega_1\varpi_1}^{\lambda}$, there is some point $\varpi_0\in\gamma[\omega_1\varpi_1]$ satisfying
	$$
	d_D(\varpi_0)\geq \mu_5^{\frac{1}{2}}\ell(\mathfrak{l}[\omega_1,\varpi_1]).$$
	\ethm
	
	\begin{figure}[htbp]
		\begin{center}
			\begin{tikzpicture}[scale=1]
				\draw (-6,-0.5) to [out=55,in=130] (0.5,0) to [out=-50,in=220] (6,0.5);
				\filldraw  (-6,-0.5)node[below] {\tiny $\eta_1=\omega=u_1$} circle (0.04);
				\filldraw  (6,0.5)node[below right] {\tiny $\varpi$} circle (0.04);
				\filldraw  (-5,0.48)node[below] {\tiny $\;\;\;\;u_2$} circle (0.04);
				\filldraw  (-4,1)node[below] {\tiny $\;\;u_3$} circle (0.04);
				\filldraw[color=red]  (2.7,-0.94)node[below] {\tiny $q$} circle (0.04);
				\node[rotate=0] at(4,-0.25) {$\mathfrak{l}$};
				
				\filldraw  (-2.8,1.265)node[below] {\tiny $u_{i-1}$} circle (0.04);
				\filldraw  (-1.25,1.12)node[below] {\tiny $u_i$} circle (0.04);
				
				\filldraw  (-0.2,0.633)node[below] {\tiny $u_m\;\;\;$} circle (0.04);
				\filldraw  (0.8,-0.318)node[below] {\tiny $u_{m+1}\;\;\;\;\;\;\;$} circle (0.04);
				
				\draw (-6,-0.5) to [out=69,in=195] (-4,1.8) to [out=15,in=110]  (-2.8,1.265);
				\filldraw  (-3.7,1.84)node[above] {\tiny $u_{0,i\!-\!1}\;\;\;$} circle (0.04);
				\node[rotate=38] at(-4.65,1.2) {\tiny $\gamma_{i-1}$};
				
				\draw  (-2.8,1.265) to [out=75,in=180] (-2.05,2.1) to [out=0,in=105] (-1.25,1.12);
				\filldraw  (-2.5,1.9) circle (0.04);
				\node at(-2.6,2.1) {\tiny$u_{2,i}\;\;\;\;$};
				\filldraw  (-2.05,2.1)node[above] {\tiny $u_{1,i}$} circle (0.04);
				\node[rotate=0] at(-1.85,1.7) {\tiny $\beta_i$};
				
				\draw  (-6,-0.5) to [out=85,in=180] (-3,4) to [out=0,in=90] (-1.25,1.12);
				\filldraw  (-3,4)node[above] {\tiny $u_{0,i}$} circle (0.04);
				\node[rotate=0] at(-1.6,3.7) {\tiny $\gamma_{i}$};
				
				\draw (-6,-0.5) to [out=95,in=180] (-0.4,8) to [out=0,in=95] (6,0.5);
				\filldraw[color=red]  (-6.022,1)node[left] {\tiny $z$} circle (0.04);
				\filldraw  (-5.915,2)node[left] {\tiny $\eta_2$} circle (0.04);
				\filldraw  (-5.715,3)node[left] {\tiny $\eta_3$} circle (0.04);
				\filldraw  (-5.32,4.2)node[left] {\tiny $\eta_{i-1}$} circle (0.04);
				\filldraw  (-4.825,5.2)node[left] {\tiny $\eta_i$} circle (0.04);
				\filldraw  (-4,6.35)node[above] {\tiny $\eta_m\;\;\;\;\;$} circle (0.04);
				\filldraw  (-3,7.225)node[above left] {\tiny $\eta_{m+1}=z_0$} circle (0.04);
				\filldraw  (3.5,6.39)node[above right] {\tiny $x_0$} circle (0.04);
				\node[rotate=0] at(0.8,8.1) {$\gamma$};
				
				\draw[color=red] (-6,-0.5) to [out=90,in=180] (-1.7,6) to [out=0,in=90]  (2.7,-0.94);
				\filldraw[color=red]  (0.5,5.25)node[above right] {\tiny $q_0$} circle (0.04);
				\node[color=red,rotate=0] at(2.5,2.8) {\tiny $\gamma_{0}$};
				
			\end{tikzpicture}
		\end{center}
		\caption{Illustration for the proof of Theorem \ref{9-15-1}}
		\label{fig-9-22}
	\end{figure}

	The proof of Theorem \ref{9-15-1} depends on Theorem \ref{22-11-16} below, which relates the length of a $\lambda$-curve and the inner distance between its end points.
	
	\bthm\label{22-11-16}
	Suppose that $f:$ $D\to D'$ is $(M,C)$-$\CQH$ and $D'$ is a $c$-uniform domain.
	For any distinct points $x\not=y\in D$ and any $\lambda$-curve $\gamma\in \Lambda_{xy}^{\lambda}$, if $k_D(x,y)>1$, then we have
	\be\label{24-5-6-1}d_D(x_0)< \mu_5^{\frac{1}{2}}\sigma_D(x,y)
	\ee
	for any $x_0\in\gamma$.
	\ethm
	
	As the proof of Theorem \ref{22-11-16} is rather technical and lengthy, we move it to Section \ref{sec-5}.

	\subsection{Auxiliary lemmas}\label{sec-3-4}
	In this section, under the \textbf{Standing Assumptions} in Section \ref{sec-3-1} and assumptions of Theorem \ref{9-15-1}, we shall establish a series of auxiliary results that are necessary for the proof of Theorem \ref{9-15-1}.

	For the convenience of later references, we collect some useful properties in Proposition \ref{11-19-1} below, where $(P_1)$ comes from Lemma \ref{10-1}$(1)$, $(P_2)$ from Lemma \ref{10-1}$(2)$, $(P_3)$ from \cite[Theorem 2.12]{Vai9'} and \cite[Theorem 3.18]{Vai10}, and finally $(P_4)$ from Lemma \ref{Thm-7-2}\eqref{Thm-7-2-2}.
	\bprop\label{11-19-1}
	$(P_1)$ For any $x\not=y\in D$, every $\lambda$-curve $\gamma\in \Lambda_{xy}^\lambda$ is a $(c_0,0)$-solid arc and belongs to $Q_{xy}^{c_0}$ with $c_0=\frac{9}{8}$;
	
	\noindent $(P_2)$ For any $\lambda$-curve $\gamma\in \Lambda_{uv}^{\lambda}$, its image $\gamma'$ under $f$ is $(c_1,h_1)$-solid, where $c_1=\frac{9}{4}M(M+1)$ and $h_1=C\big(\frac{9}{4}M+1\big)$;
	\medskip
	
	\noindent $(P_3)$
	Both $D$ and $D'$ are $\delta$-Gromov hyperbolic, where $\delta=\delta(c,C,M)$;
	\medskip
	
	\noindent $(P_4)$ Both $(D, k_D)$ and $(D', k_{D'})$ are $(\mu_0,\rho_0)$-Rips spaces.
	\eprop
	
	Throughout this section, we let $c_2=c_2(c,c_1,h_1)$ be the constant given by Theorem \ref{Theo-B} and $c_3=c_3(c)$ be a constant given by Theorem \ref{Theo}.
	
	The first lemma is an easy consequence of Lemma \ref{lem-5-14}.
	\blem\label{23-08-11-2}
	For any $w\in\mathfrak{l}$, $d_D(w)\leq 2\ell(\mathfrak{l}).$
	\elem
	\bpf Note that by our choice of constants in Section \ref{subsec:notation}, $\frac{4}{5}\mu_5> 18e^9$ and thus $\ell(\gamma) >\frac{4}{5}\mu_5\ell(\mathfrak{l})>18e^9\ell(\mathfrak{l})$. The lemma follows from Lemma \ref{lem-5-14}.
	\epf
	
	\blem\label{23-08-12-1}
	$\diam(\gamma')\leq 4c_2e^{h_1}|\omega'-\varpi'|$.
	\elem
	\bpf Lemma \ref{23-08-11-2} and the assumption $\ell(\gamma)>\mathfrak{C} \ell(\mathfrak{l})$ in Theorem \ref{9-15-1} yield $$\ell(\gamma)\geq \frac{\mathfrak{C}}{2} d_D(\omega).$$
	Then we get from Lemma \ref{lem-2.2} and \eqref{(2.1)} that
	$$k_D(\omega,\varpi)\geq \ell_k(\gamma)-\frac{1}{4}\geq \log\Big(1+\frac{\ell(\gamma)}{d_D(\omega)}\Big)-\frac{1}{4}\geq \log\Big(1+ \frac{\mathfrak{C}}{2}\Big)-\frac{1}{4}.$$
	Since $f$ is $(M,C)$-CQH and $\mathfrak{C}\geq \frac{4}{5}\mu_5$, it follows
	$$k_{D'}(\omega',\varpi')\geq \frac{1}{M}(k_D(\omega,\varpi)-C)\geq \frac{1}{M}\log\Big(1+ \frac{\mathfrak{C}}{2}\Big)-\frac{1+4C}{4M}>1.$$
	Thus an application of Lemma \ref{lem-ne-1}\eqref{9-10-2} gives
	\be\label{24-5-9-1}|\omega'-\varpi'|\geq \frac{1}{2}\max\{d_{D'}(\omega'),d_{D'}(\varpi')\}.\ee
	
	Since $\gamma'$ is $(c_1,h_1)$-solid by Proposition \ref{11-19-1}$(P_2)$ and $D'$ is $c$-uniform, Theorem \ref{Theo-B}(\ref{THB-2}) gives
	$$\diam(\gamma')\leq c_2\max\{|\omega'-\varpi'|, 2(e^{h_1}-1)\min\{d_{D'}(\omega'),d_{D'}(\varpi')\}\},
	$$ where $c_2=c_2(c,c_1,h_1)$ is from Theorem \ref{Theo-B}. Then \eqref{24-5-9-1} guarantees that
	$$\diam(\gamma')\leq 4c_2e^{h_1}|\omega'-\varpi'|.$$
	The proof is thus complete.
	\epf
	
	Let $z_0$ (resp. $x_0$) be the first point on $\gamma$ along the direction from $\omega$ to $\varpi$ (resp. from $\varpi$ to $\omega$) such that
	\be\label{9-11-2}
	|\omega'-z'_0|=\frac{1}{4}\diam(\gamma')\;\quad (\mbox{resp.}\; |\varpi^\prime-x_0^\prime|=\frac{1}{4}\diam(\gamma^\prime));
	\ee
	see Figure \ref{fig-9-22} for an illustration. Then, $z_0\in \gamma[\omega,x_0]$ and $x_0\in \gamma[\varpi,z_0]$.
	
	\blem\label{22-11-10-1}
	$(1)$ For any $z\in\gamma[\omega,z_0]$, $\diam(\gamma'[\omega',z'])\leq c_2 d_{D'}(z'),$ and for any $y\in\gamma[\varpi,x_0]$, $\diam(\gamma'[\varpi',y'])\leq c_2 d_{D'}(y').$
	
	$(2)$
	$\ell(\gamma[z_0, x_0])\leq (1+4c_2)^{c_3 M}e^{\frac{1}{2}+C}\min\{d_D(z_0), d_D(x_0)\}$.
	\elem
	\bpf (1). We only show the first claim in statement (1), since the proof of second claim is similar. For this, let $z\in\gamma[\omega,z_0]$. Then the choice of $z_0$ in \eqref{9-11-2} implies that
	$$\diam(\gamma'[\omega',z'])\leq 2|\omega'-z_0'|=\frac{1}{2}\diam(\gamma'),$$
	and thus,
	$$\diam(\gamma'[\varpi',z'])\geq\frac{1}{2}\diam(\gamma')\geq \diam(\gamma'[\omega',z']).$$
	Then we know from Theorem \ref{Theo-B}(\ref{THB-1}) that
	$$\diam(\gamma'[\omega',z'])\leq c_2 d_{D'}(z'),$$
	which verifies the desired claim.
	
	(2). Since $D'$ is $c$-uniform, Theorem \ref{Theo} gives
	\beqq
	k_{D^\prime}(z_0^\prime, x_0^\prime)
	&\leq& c_3\log\left(1+\frac{|x_0^\prime-z_0^\prime|}{\min\{d_{D^\prime}(x_0^\prime), d_{D^\prime}(z_0^\prime)\}}\right)\\
	&\leq& c_3\log\left(1+\frac{\diam(\gamma^\prime)}{\min\{d_{D^\prime}(x_0^\prime), d_{D^\prime}(z_0^\prime)\}}\right).
	\eeqq
	By the choice of $z_0$ and $x_0$ in \eqref{9-11-2} and statement (1) of Lemma \ref{22-11-10-1}, we have
	\beqq
	k_{D^\prime}(z_0^\prime, x_0^\prime)\leq c_3\log(1+4c_2).
	\eeqq
	Since $f$ is $(M,C)$-CQH, it follows from Lemmas \ref{lem-2.2} and \ref{lem-2.3} that
	\beqq
	\ell_k(\gamma[z_0, x_0])-\frac{1}{2}\leq k_D(z_0, x_0)\leq Mk_{D^\prime}(z_0^\prime, x_0^\prime)+C\leq c_3 M\log(1+4c_2)+C,
	\eeqq
	and thus, we infer from \eqref{(2.1)} that
	\beqq
	\log\left(1+\frac{\ell(\gamma[z_0, x_0])}{\min\{d_{D}(z_0), d_{D}(x_0)\}}\right)\leq c_3 M\log(1+4c_2)+C+\frac{1}{2}.
	\eeqq
	This proves the second statement, and hence, the proof is complete.
	\epf
	
	Without loss of generality, we assume that
	\beq\label{12-7-3}
	\ell(\gamma[\omega, z_0])\geq \ell(\gamma[\varpi, x_0]).
	\eeq
	Under this assumption, we have the following lemma.
	
	\blem\label{12-7-1}
	$\ell(\gamma) \leq 2\left(1+(1+4c_2)^{c_3 M}e^{\frac{1}{2}+C}\right)\ell(\gamma[\omega, z_0]).$
	\elem
	\bpf
	First, we know from \eqref{12-7-3} and Lemma \ref{22-11-10-1}$(2)$ that
	\beq\label{eq-1.4}
	\begin{aligned}
		\ell(\gamma) &= \ell(\gamma[\omega, z_0])+\ell(\gamma[\varpi, x_0])+\ell(\gamma[z_0, x_0])\\
		&\leq 2\ell(\gamma[\omega, z_0])+ (1+4c_2)^{c_3 M}e^{\frac{1}{2}+C}\min\{d_D(z_0), d_D(x_0)\}.
	\end{aligned}
	\eeq
	
	Next, we show the following estimate:
	\beq\label{eq-1.5}
	d_D(z_0)\leq 2\ell(\gamma[\omega, z_0]).
	\eeq
	Suppose on the contrary that
	\be\label{12-7-2} d_D(z_0)> 2\ell(\gamma[\omega, z_0]).\ee
	Then it holds
	\be\label{12-7-4}
	d_D(\omega)\geq d_D(z_0)-\ell(\gamma[\omega, z_0])\geq \frac{1}{2}d_D(z_0).
	\ee
	Consequently, \eqref{eq-1.4}, \eqref{12-7-2} and \eqref{12-7-4} give
	$$
	\ell(\gamma) \leq \left(1+(1+4c_2)^{c_3 M}e^{\frac{1}{2}+C}\right)d_D(z_0)\leq 2\left(1+(1+4c_2)^{c_3 M}e^{\frac{1}{2}+C}\right)d_D(\omega).
	$$

	On the other hand, since $\ell(\gamma)\geq \frac{4}{5}\mu_5\ell(\mathfrak{l})$, Lemma \ref{23-08-11-2} implies
	$$\ell(\gamma)\geq \frac{2}{5}\mu_5 d_D(\omega).$$
	This is impossible as our choice of constant in Section \ref{subsec:notation} implies
	$$\frac{2}{5}\mu_5>2\left(1+(1+4c_2)^{c_3 M}e^{\frac{1}{2}+C}\right).$$ Thus \eqref{eq-1.5} holds and the desired estimate in Lemma \ref{12-7-1} follows from \eqref{eq-1.4} and \eqref{eq-1.5}.
	\epf
	
	
	In the following lemmas, it is helpful to see Figure \ref{fig-9-22} for an intuition.
	\blem\label{22-11-12-1} Let $z\in\gamma[\omega,z_0]$ such that $k_D(\omega,z)\geq \mu_3$. Then
	there exists $q\in\mathfrak{l}$ such that
	$|\omega'-q'|=\mu_3^{-1}\diam(\gamma'[\omega',z'])$, and for any $u\in \mathfrak{l}[\omega,q]\setminus\{q\}$, $$|\omega'-u'|<\mu_3^{-1}\diam(\gamma'[\omega',z']).$$
	\elem
	\bpf
	Since $\mu_3> 4c_3e^{h_1}$, we see from Lemma \ref{23-08-12-1} that
	$$|\omega'-\varpi'|>\mu_3^{-1}\diam(\gamma'[\omega',z']).$$ This implies that there exists $q\in \mathfrak{l}$ such that $$|\omega'-q'|=\mu_3^{-1}\diam(\gamma'[\omega',z']),$$ and for any $u\in \mathfrak{l}[\omega,q]\setminus\{q\}$, $$|\omega'-u'|<\mu_3^{-1}\diam(\gamma'[\omega',z']).$$ This completes the proof of lemma.
	\epf
	
	
	\blem\label{9-19-1}
	There are an integer $m$ and a finite sequence of points $\{\eta_i\}_{i=1}^{m+1}$ on $\gamma[\omega,z_0]$ such that
	\begin{enumerate}
		\item\label{11-17-1}
		$m>e^{8\mu_4}$ and $m \mu_4\leq \ell_k(\gamma[\omega,z_0])< (m+1)\mu_4$;
		\item\label{22-11-25-2}
		for each $i\in\{2,\cdots,m\}$, $\eta_i\in\gamma[\eta_{i-1},z_0]$,
		$$\ell_k(\gamma[\eta_{i-1},\eta_i])=\mu_4\;\;\mbox{and}\;\;\mu_4\leq\ell_k(\gamma[\eta_m,\eta_{m+1}])<2\mu_4,$$ where $\eta_1=\omega$ and $\eta_{m+1}=z_0$;
		\item\label{22-11-25-1}
		for any $i\not= j\in \{1,2,\cdots,m+1\}$, $$k_{D'}(\eta'_i,\eta'_j)\geq  (2M)^{-1}\mu_4;$$
		\item\label{9-23-1} for each $i\in\{2,\cdots,m\}$,
		\begin{enumerate}
			\item\label{22-11-28-1}
			$\ell(\gamma[\eta_{i-1},\eta_i])\leq e^{2\mu_4}d_D(\eta_i);$
			\item\label{23-08-13-0}
			$\diam(\gamma'[\eta'_{i-1},\eta'_i])\leq  c_2d_{D'}(\eta'_i);$
			\item\label{23-08-13-7}
			$\diam(\gamma'[\eta'_{i-1},\eta'_i])\geq e^{(3c_3M)^{-1}\mu_4}d_{D'}(\eta'_{i-1});$
			\item\label{23-08-23-3}
			$\diam(\gamma'[\eta'_1,\eta'_i])\geq c_2^{-1}e^{(3c_3M)^{-1}\mu_4}\diam(\gamma'[\eta'_1,\eta'_{i-1}])$.
		\end{enumerate}
	\end{enumerate}
	\elem
	
	\bpf
	(1). By Lemmas  \ref{23-08-11-2}, \ref{12-7-1} and the assumption $\ell(\gamma)>\mathfrak{C}\ell(\mathfrak{l})$, we have
	$$\ell(\gamma[\omega,z_0])\geq \frac{1}{\wp_2}\ell(\gamma) \geq  \frac{\mathfrak{C}}{2\wp_2}d_D(w),$$
	where $\wp_2=2\left(1+(1+4c_2)^{c_3 M}e^{\frac{1}{2}+C}\right)$. Then it follows from \eqref{(2.1)} that
	$$\ell_k(\gamma[\omega,z_0])\geq \log\Big(1+\frac{\ell(\gamma[\omega,z_0])}{d_D(\omega)}\Big)\geq \log\big(1+\frac{\mathfrak{C}}{2\wp_2}\big)>\mu_4e^{8\mu_4},$$
	where in the last inequality we used the  the assumption $\mathfrak{C}>\frac{4}{5}\mu_5=\frac{8}{5}\mu_4e^{2\mu_4e^{8\mu_4}}$.
	
	Consequently, there exists an integer $m$ such that \eqref{11-17-1} holds.
	
	(2). By the last estimate in the proof of \eqref{11-17-1}, there exists a sequence $\{\eta_i\}_{i=1}^{m+1}$ in $\gamma[\omega,z_0]$ satisfying \eqref{22-11-25-2}.
	
	(3). Let $\{i,j\}\subset\{1,\cdots,m+1\}$ with $i<j$. Since  $f$ is $(M,C)$-CQH, we have
	$$ k_{D'}(\eta'_i,\eta'_j) \geq \frac{1}{M}\left(k_D(\eta_i,\eta_j)-C\right).$$
	It follows from Lemmas \ref{lem-2.2} and \ref{lem-2.3}, together with \eqref{22-11-25-2}, that
	\beqq 
	k_{D'}(\eta'_i,\eta'_j)\geq \frac{1}{M}(k_D(\eta_i,\eta_j)-C)\geq \frac{1}{M}\big(\ell_k(\gamma[\eta_i,\eta_j])-\frac{1}{2}-C\big)
	\geq  (2M)^{-1}\mu_4, 
	\eeqq
	as  $\mu_4>\frac{1}{2}+C$. This proves \eqref{22-11-25-1}.
	
	(4). For each $i\in\{2,\cdots,m+1\}$, \eqref{(2.1)} implies
	$$\log\Big(1+\frac{\ell(\gamma[\eta_{i-1},\eta_i])}{d_D(\eta_i)}\Big)\leq \ell_k(\gamma[\eta_{i-1},\eta_i]).$$
	Then \eqref{22-11-28-1}  follows from \eqref{22-11-25-2}.
	
	The statement \eqref{23-08-13-0}  follows immediately from Lemma \ref{22-11-10-1}:
	$$\diam(\gamma'[\eta'_i,\eta'_{i-1}])\leq \diam(\gamma'[\eta'_1,\eta'_{i-1}])\leq c_2d_{D'}(\eta'_{i-1}).$$
	
	Since $k_{D'}(\eta'_{i-1},\eta'_i)>1$ by \eqref{22-11-25-1}, it follows from Lemma \ref{lem-ne-1}\eqref{9-10-2} and  \eqref{23-08-13-0} that
	\be\label{11-18-2}
	d_{D'}(\eta'_{i-1})\leq 2|\eta'_{i-1}-\eta'_i|\leq 2c_2d_{D'}(\eta'_i),
	\ee
	Note that Theorem \ref{Theo} and \eqref{22-11-25-1} ensure that
	$$(2M)^{-1}\mu_4\leq  k_{D'}(\eta'_{i-1},\eta'_i)\leq c_3\log\Big(1+\frac{|\eta'_{i-1}-\eta'_i|}{\min\{d_{D'}(\eta'_{i-1}),d_{D'}(\eta'_i)\}}\Big).$$
	Thus we see from the previous estimate and \eqref{11-18-2} that
	\beqq
	\begin{aligned}
	\diam(\gamma'[\eta'_{i-1},\eta'_i])&\geq |\eta'_{i-1}-\eta'_i|\geq (e^{(2c_3M)^{-1}\mu_4}-1)\min\{d_{D'}(\eta'_{i-1}),d_{D'}(\eta'_i)\}
	\\ &\geq e^{(3c_3M)^{-1}\mu_4}d_{D'}(\eta'_{i-1}).
	\end{aligned}
	\eeqq
	This shows \eqref{23-08-13-7}.
	
	By Lemma \ref{22-11-10-1} and \eqref{23-08-13-7}, we have
	\beqq
	\begin{aligned}
	\diam(\gamma'[\eta'_1,\eta'_i])&\geq \diam(\gamma'[\eta'_{i-1},\eta'_i])\geq e^{(3c_3M)^{-1}\mu_4}d_{D'}(\eta'_{i-1})
	\\ &\geq c_2^{-1}e^{(3c_3M)^{-1}\mu_4}\diam(\gamma'[\eta'_1,\eta'_{i-1}]).
\end{aligned}
	\eeqq
	This proves \eqref{23-08-23-3} and thus completes the proof of lemma.
	\epf
	
	\blem\label{9-14-1} Let $u_1=\omega=\eta_1$. Then the following statements hold.
	\begin{enumerate}
		\item\label{24-5-9-2}
		For each $i\in\{2,\cdots,m+1\}$,
		there exists $u_i\in \mathfrak{l}[u_{i-1},\varpi]$ such that
		$|u'_1-u'_i|=\mu_3^{-1}\diam(\gamma'[\eta'_1,\eta'_i])$, and for any $u\in \mathfrak{l}[u_1,u_i]\setminus\{u_i\}$, $$|u'_1-u'|<\mu_3^{-1}\diam(\gamma'[\omega',\eta_i']).$$
		\item\label{24-5-9-3}
		For each $i\in\{2,\cdots,m+1\}$ and any $\lambda$-curve $\gamma_i\in \Lambda_{u_{i-1} u_i}^{\lambda}$, there exists $u_{0,i}\in\gamma_i$ such that
		\begin{enumerate}
			\item\label{9-9-2}
			$k_D(u_{0,i},\eta_i)\leq c_3M\log (1+16c_2^2 e^{h_1}\mu_3)+C.$
			\item\label{23-08-12-9}
			$d_D(\eta_i)\leq \mu_3^{2c_3M}d_D(u_{0,i}).$
		\end{enumerate}
	\end{enumerate}
	\elem
	
	\bpf
	(1). Let $i\in\{2,\cdots,m+1\}$. Since $f$ is $(M,C)$-CQH, we know from Lemma \ref{9-19-1}(\ref{22-11-25-1}) that
	$$k_D(\eta_1,\eta_i)\geq \frac{1}{M}(k_{D'}(\eta'_1,\eta'_i)-C)\geq (3M^2)^{-1}\mu_4>\mu_3.$$
	Applying  Lemma \ref{22-11-12-1}, with $z$ replaced by $\eta_i$, we obtain that there exists
	$u_i\in \mathfrak{l}$ such that \eqref{24-5-9-2} holds.
	
	(2). Before the proof of \eqref{24-5-9-3}, let $u'_{0,i}\in\gamma'_i$ be the first point along the direction from $u'_{i-1}$ to $u'_i$ such that
	\be\label{9-20-1}
	|u'_{i-1}-u'_{0,i}|=\frac{1}{4}\diam(\gamma'_i).
	\ee
	Then it follows from Theorem \ref{Theo-B} that
	\beqq\label{9-15-2}
	\frac{1}{4}|u'_{i-1}-u'_i|\leq \frac{1}{4}\diam(\gamma'_i)=|u'_{i-1}-u'_{0,i}|\leq c_2 d_{D'}(u'_{0,i}).
	\eeqq
	By Lemmas \ref{9-19-1}\eqref{23-08-23-3} and \ref{9-14-1}\eqref{24-5-9-2}, we get
	$$|u'_{i-1}-u'_i|\geq \mu_3^{-1}(\diam(\gamma'[\eta'_1,\eta'_i])-\diam(\gamma'[\eta'_1,\eta'_{i-1}]))\geq \mu_3^{-1}(1-c_2e^{-(3c_3M)^{-1}\mu_4})\diam(\gamma'[\eta'_1,\eta'_i]).$$
	Furthermore, (\ref{9-15-2}) yields that
	\be\label{H-26-2-21-1}
	\diam(\gamma'[\eta'_1,\eta'_i])<2\mu_3|u'_{i-1}-u'_i|\leq 8c_2 \mu_3d_{D'}(u'_{0,i}).
	\ee
	Since $$(1-c_2e^{-(3c_3M)^{-1}\mu_4})\diam(\gamma'[\eta'_1,\eta'_i])\leq\diam(\gamma'[\eta'_1,\eta'_i])-\diam(\gamma'[\eta'_1,\eta'_{i-1}])\leq \diam(\gamma'[\eta'_{i-1},\eta'_i]),$$
	we have by Lemmas \ref{9-19-1} \eqref{23-08-13-0} and \eqref{23-08-23-3} that
	\be\label{H-26-2-21-2}
	\diam(\gamma'[\eta'_1,\eta'_i])<2\diam(\gamma'[\eta'_{i-1},\eta'_i])\leq 2c_2 d_{D'}(\eta'_i).
	\ee
	Then, we obtain from Lemma \ref{9-14-1}\eqref{24-5-9-2} and \eqref{9-20-1} that
	\beqq
\begin{aligned}
	|u'_{0,i}-\eta'_i|&\leq |u'_1-\eta'_i|+|u'_1-u'_{i-1}|+|u'_{i-1}-u'_{0,i}|\\ &=
	\diam(\gamma'[\eta'_1,\eta'_i])+\mu_3^{-1}\diam(\gamma'[\eta'_1,\eta'_{i-1}])+\frac{1}{4}\diam(\gamma'_i)
	\\ &< (1+\mu_3^{-1})\diam(\gamma'[\eta'_1,\eta'_i])+\frac{1}{4}\diam(\gamma'_i),
\end{aligned}
	\eeqq
	and again by Theorem \ref{Theo-B} and Lemma \ref{9-14-1}\eqref{24-5-9-2}, we have
	\beqq
	\begin{aligned}
	|u'_{0,i}-\eta'_i|&< (1+\mu_3^{-1})\diam(\gamma'[\eta'_1,\eta'_i])+c_2e^{h_1}|u'_{i-1}-u'_i|
	\\ &\leq (1+\mu_3^{-1})\diam(\gamma'[\eta'_1,\eta'_i])+c_2e^{h_1} (|u'_1-u'_{i-1}|+|u'_1-u'_i|)
	\\ &< 2c_2e^{h_1} \diam(\gamma'[\eta'_1,\eta'_i]),
\end{aligned}
	\eeqq
	which, together with (\ref{H-26-2-21-1}) and (\ref{H-26-2-21-2}), shows that
	$$|u'_{0,i}-\eta'_i|\leq 16c_2^2 e^{h_1}\mu_3\min\{d_{D'}(u'_{0,i}),d_{D'}(\eta'_i)\}.$$
	Hence, it follows from Theorem \ref{Theo} that
	$$k_{D'}(u'_{0,i},\eta'_i)\leq c_3\log\big(1+\frac{|u'_{0,i}-\eta'_i|}{\min\{d_{D'}(u'_{0,i}),d_{D'}(\eta'_i)\}}\big)\leq c_3\log (1+16c_2^2 e^{h_1}\mu_3).$$
	Meanwhile, since $f$ is $(M,C)$-CQH, it follows from (\ref{(2.1)}) that
	$$\log \frac{d_D(\eta_i)}{d_D(u_{0,i})}\leq k_D(u_{0,i},\eta_i)\leq Mk_{D'}(u'_{0,i},\eta'_i)+C<c_3M\log (1+16c_2^2 e^{h_1}\mu_3)+C,$$
	and so $$d_D(\eta_i)< \mu_3^{2c_3M} d_D(u_{0,i}).$$ The above two inequalities give (\ref{9-9-2}) and (\ref{23-08-12-9})), and hence, the proof is complete.
	\epf

	\subsection{Proof of Theorem \ref{9-15-1}}\label{sec-3-5}
	In this section, we shall keep all the notations used in Section \ref{sec-3-4}.
	
	From the assumption in this theorem we have
	$$\ell(\gamma)\geq \frac{4}{5}\mu_5 \ell(\mathfrak{l}),$$
	and so Lemma \ref{12-7-1} implies that
	$$\ell(\gamma[\omega,z_0])>\mu_5^{\frac{3}{4}}\ell(\mathfrak{l}).$$
	Hence, there must exists $s\in\{2,\cdots,m+1\}$ such that
	$$
	\ell(\gamma[\eta_{s-1},\eta_s])> \mu_5^{\frac{3}{4}}\ell(\mathfrak{l}[u_{s-1},u_s]).
	$$
	Combining Lemmas \ref{9-19-1}(\ref{22-11-28-1}) and \ref{9-14-1}(\ref{9-9-2}), we have
	$$d_D(u_{0,s})\geq \mu_3^{-2c_3M}d_D(\eta_s)>e^{-3\mu_4} \ell(\gamma[\eta_{s-1},\eta_s])\geq e^{-3\mu_4}\mu_5^{\frac{3}{4}}\ell(\mathfrak{l}[u_{s-1},u_s])>\mu_5^{\frac{1}{2}}\ell(\mathfrak{l}[u_{s-1},u_s]).$$
	By taking $\omega_1=u_{s-1}$, $\varpi_1=u_s$ and $\varpi_0=u_{0,s}$, and the proof of Theorem \ref{9-15-1} is thus complete.
	\qed

	\section{Proofs of the main results}\label{sec-4}

	The purpose of this section is to prove the main results in this paper, i.e., Theorems \ref{thm1.1} $\sim$ \ref{thm1.3}, based on Theorems \ref{9-15-1} and \ref{22-11-16}.
	For the proof of Theorem \ref{thm1.1}, we need the following useful result bounding $\ell(\gamma)$ in terms of $\ell(\alpha)$. Recall the definition of $\mu_5$ given in Section \ref{subsec:notation}.
	\bprop\label{9-19-2}
	Under the \textbf{Standing Assumptions} in Section \ref{sec-3-1}, let $\alpha$ be an arc in $D$ joining $\omega$ and $\varpi$. Then we have
	$$\ell(\gamma)\leq \frac{4}{5}\mu_5\ell(\alpha).$$
	\eprop
	\bpf
	It suffices to consider the case when $\alpha$ is rectifiable. Suppose on the contrary that
	\be\label{24-4-10-1}
	\ell(\gamma)> \frac{4}{5}\mu_5 \ell(\alpha).
	\ee
	
	We claim $k_D(\omega,\varpi)> 1$. Otherwise, we know from \cite[Lemma 2.31]{Vai7}  that $\gamma\in \Lambda_{\omega\varpi}^{\lambda}$ is a  $c$-uniform arc with $c=18e^9$, which contradicts with
	the assumption in this theorem.
	
	To reach a contradiction, let $u_1=\omega$ and $v_1=\varpi$. Then it follows from Theorem \ref{9-15-1} and \eqref{24-4-10-1} that there are points
	$u_2\not=v_2$ on $\alpha$ such that $u_2\in \alpha[u_1, v_2]$, $v_2\in \alpha[u_2, v_1]$, and there is an element $\gamma_{2}\in\Lambda_{u_2v_2}^\lambda$ such that
	there exists some $u_0\in\gamma_2$ satisfying
	\be\label{24-4-10-2}
	d_D(u_0)\geq \mu_5^{\frac{1}{2}}\ell(\alpha[u_2,v_2]).
	\ee
	Combining Theorem \ref{22-11-16} with \eqref{24-4-10-2}, we get a contradiction
	Consequently, the proof is complete.
	\epf

	\begin{proof}[Proof of Theorem \ref{thm1.1}]
		Let $\vartheta$ be a $c_0$-quasigeodesic, $c_0\geq 1$, in $D$ with end points $\omega$ and $\varpi$, and let $L$ be an arc in $D$ with the same end points as $\vartheta$.
		By Proposition \ref{11-19-1}$(P_3)$, $D$ is $\delta$-Gromov hyperbolic, for some $\delta=\delta(c,C,M)$.
		
		We may assume that $L$ is rectifiable with $\ell(L)<\ell(\vartheta)$. Then
		as both $\vartheta$ and $\gamma$ are $(1+c_0)$-quasigeodesics, we obtain from Theorem \ref{24-4-11-12}  that
		\be\label{12-11-6}
		\ell(\vartheta)\leq \lambda_0\ell(\gamma),
		\ee
		where $\lambda_0=4e^{8(1+c_0)^2\mu_2}$. On the other hand, 	we obtain from Proposition \ref{9-19-2} that
		$$\ell(\gamma)\leq \frac{4}{5}\mu_5\ell(L).$$
		Combining the above two inequalities leads to our desired estiamte
		\be\label{12-11-2}\ell(\vartheta)\leq \frac{4}{5}\mu_5 \lambda_0\ell(L)=:b\ell(L).\ee
		Then, Theorem \ref{thm1.1} follows from the arbitrariness of $\omega$ and $\varpi$ in $D$.
	\end{proof}
	
	\begin{proof}[Proof of Theorem \ref{thm1.3}]
		By \cite[Theorem 4.8]{Vai3} and \cite[Theorem 4.11]{Vai4}, we know that the image of a $c_0$-quasigeodesic under an
		$M$-bi-Lipschitz homeomorphism (in the quasihyperbolic metric) is a $c'_0$-quasigeodesic with $c'_0=2c_0M^2(1+M^2)$. As the inverse of an $M$-bi-Lipschitz homeomorphism is again $M$-bi-Lipschitz, Theorem \ref{thm1.3} follows directly from Theorem \ref{thm1.1}.
	\end{proof}
	
%
	
	\section{Proof of Theorem \ref{22-11-16}}\label{sec-5}
	
	The purpose of this section is to prove Theorem \ref{22-11-16}, which has played a key role in the proof of Theorem \ref{thm1.1}. As in Section \ref{sec-3}, we fix $\lambda=\frac{1}{4}$.
	
	\subsection{Preparation for the proof}\label{subsec:step 1}

	Let $\alpha$ be a rectifiable curve in $D$ joining $x$ and $y$. Since $k_D(x,y)\leq 1$, it follows from Lemma \ref{lem-ne-1} that for each $w\in\alpha$,
	\be\label{22-11-16-1}d_D(w)\leq 2\ell(\alpha).\ee

	We shall prove the theorem by a contradiction. To reach a contradiction, let $x_1=x$ and $y_1=y$.
	To this end, suppose that there exists some point $x_0\in\gamma$ such that
	\be\label{22-11-16-2} 
	d_D(x_0)\geq \mu_5^{\frac{1}{2}}\ell(\alpha).
	\ee
	For each $z\in\alpha$ and $\gamma_{x_0z}\in\Lambda_{x_0z}^{\lambda}$, we infer from (\ref{22-11-16-1}) and (\ref{22-11-16-2}) that
	\be\label{10-3-2}
	\ell(\gamma_{x_0z})\geq d_D(x_0)-d_D(z)\geq d_D(x_0)-2\ell(\alpha)\geq (\mu_5^{\frac{1}{2}}-2)\ell(\alpha).
	\ee

	\subsection{$\zeta$-sequence on $\alpha[x_1,y_1]$}\label{sub-5.1}
	In this section, we shall construct a finite sequence of points on $\alpha$, called $\zeta$-sequence, with good controls on the length of $\lambda$-curves between successive points.
	\begin{prop}\label{24-5-7-3}
		There exist an integer $m\geq 2$ and a sequence $\{x_i\}_{i=1}^{m+1}$ of points on $\alpha$ such that the following conclusions hold.
		
		$(a)$ For each $i\in \{2,\cdots,m\}$ and for any $\lambda$-curve $\alpha_{i-1}\in \Lambda_{x_{i-1}x_i}^{\lambda}$, we have
		\beqq\label{lem2-eq-5}
		\mu_6^{\frac{1}{4}}\ell(\alpha)\leq \ell(\alpha_{i-1})\leq \mu_6^{\frac{1}{2}}\ell(\alpha).
		\eeqq

		$(b)$ For any $z\in \alpha[x_i,y_1]$, and for each $\lambda$-curve $\varsigma_{i-1}\in \Lambda_{x_{i-1}z}^{\lambda}$,
		\beqq\label{1-16-1}
		\ell(\varsigma_{i-1})>\mu_6^{\frac{1}{4}}\ell(\alpha).
		\eeqq
		
		$(c)$ There is an element  $\alpha_m\in \Lambda_{x_mx_{m+1}}^\lambda$ such that
		\beqq\label{1-16-01}
		0<\ell(\alpha_m)\leq \mu_6^{\frac{1}{2}}\ell(\alpha).
		\eeqq
	\end{prop}
	
	\begin{figure}[htbp]
		\begin{center}
			\begin{tikzpicture}[scale=2.5]
				\draw (-2,0) to [out=60,in=180] (-1,0.8) to [out=0,in=120] (0,0) to [out=300,in=180] (1.7,-0.7) to [out=0,in=-128] (3.304,0.1);
				\filldraw  (-2,0) circle (0.02);
				\node at(-1.85,-0.15) {\small $x\!\!=\!x_1\!\!=\!\zeta_1$};
				\filldraw  (-1.4,0.7) circle (0.02);
				\node[rotate=-60] at(-1.18,0.3) {\small $x_2=\zeta_2$};
				\filldraw  (-0.7,0.745)node[below, rotate=-15] {\small $x_3=\zeta_3$} circle (0.02);
				\filldraw  (-0.25,0.392)node[below left] {\small $x_{i-1}=\zeta_{i-1}$} circle (0.02);
				\filldraw  (0.2,-0.265)node[below, xshift=0.3cm] {\small $x_{i}\!\!=\!\zeta_{i}\;\;\;\;\;\;\;\;\;\;\;\;\;\;\;\;$} circle (0.02);
				\filldraw  (0.85,-0.622)node[below, rotate=15] {\small $x_{i+1}\!\!=\zeta_{i+1}\;\;\;\;\;\;\;\;$} circle (0.02);
				\filldraw  (1.45,-0.695)node[below, rotate=0] {\small $\;\;\;\;\;\;\;\;\;\;\;\;x_{m-1}\!\!=\zeta_{m-1}$} circle (0.02);
				\filldraw  (2.5,-0.555)node[below] {\small $\;\;\;x_m$} circle (0.02);
				\filldraw  (3.304,0.1)node[below, rotate=50] {\small $y_1\!\!=\!x_{m\!+\!1}\!\!=\!\zeta_m$} circle (0.02);
				
				\draw (-2,0) to [out=95,in=180] (-1.7,1) to [out=0,in=110] (-1.4,0.7);
				\filldraw  (-1.95,0.865)node[below right] {\tiny $w_1$} circle (0.02);
				\node[] at(-1.82,1.1) {\tiny $\alpha_1\in \Lambda_{x_1\!x_2}^{\lambda}$};
				
				\draw (-1.4,0.7) to [out=85,in=180] (-1.05,1.3) to [out=0,in=100] (-0.7,0.745);
				\filldraw  (-1.25,1.213)node[below right] {\tiny $w_2$} circle (0.02);
				\node[above] at(-1.05,1.3) {\tiny $\alpha_2\in \Lambda_{x_2\!x_3}^{\lambda}$};
				
				\draw (-0.25,0.392) to [out=85,in=180] (-0.1,0.6) to [out=0,in=90] (0.2,-0.265);
				\filldraw  (0.15,0.43)node[below left] {\tiny $w_{i-1}$} circle (0.02);
				\node[above, rotate=-10] at(0.037,0.55) {\tiny $\alpha_{i-1}\in \Lambda_{x_{i-1}\!x_i}^{\lambda}$};

				\draw (0.2,-0.265) to [out=80,in=185] (0.5,0.1) to [out=0,in=90] (0.85,-0.622);
				\filldraw  (0.65,0.06)node[below] {\tiny $w_{i}$} circle (0.02);
				\node[above, rotate=-10,xshift=0.225cm,yshift=-0.06cm] at(0.5,0.1) {\tiny $\alpha_i\in \Lambda_{x_i\!x_{i+1}}^{\lambda}$};

				\draw (1.45,-0.695) to [out=67,in=180] (2,-0.35) to [out=0,in=135] (2.5,-0.555);
				
				\filldraw[xshift=0.6cm,yshift=-0.085cm]  (1.5,-0.27)node[below] {$\eta_1$} circle (0.02);
				\node[rotate=5, xshift=1cm ,yshift=-0.4cm] at(1.45,-0.135) {\tiny $\alpha_{m-1}\in \Lambda_{x_{m\!-\!1}\!x_m}^{\lambda}$};
				
				\draw (2.5,-0.555) to [out=95,in=220] (2.72,0) to [out=40,in=170] (3.304,0.1);
				\filldraw[xshift=0.6cm]  (2.12,0) circle (0.02);
				\node[xshift=01.5cm] at(2.2,-0.1) {\small $\eta_2$};
				\node[xshift=1.6cm, rotate=10] at(2.1,0.15)  {\tiny $\alpha_{m}\in \Lambda_{x_m\!x_{\!m+\!1}}^{\lambda}$};
				
				\draw (1.45,-0.695) to [out=80,in=200] (2.2,0.4) to [out=20,in=135] (3.304,0.1);
				\filldraw[xshift=0.6cm,yshift=-0.145cm]  (1.18,0.265)node[above, rotate=45] {\tiny $w_{m\!-\!1}$} circle (0.02);
				\filldraw[xshift=0.6cm,yshift=-0.12cm]  (1.55,0.5)node[below] {$\eta_3$} circle (0.02);
				\node[xshift=1.3cm,yshift=-0.15cm, rotate=-5] at(2.1,0.65) {\tiny $\varsigma_{m}{\color{red}=\alpha_{m-1}}\in \Lambda_{x_{m-1}x_{m+1}}^{\lambda}$};
				
				\node at(0.38,-0.6) {\large $\alpha$};
				
				
				\filldraw  (0.5,4.5)node[above] { $x_0$} circle (0.02);
				\coordinate (Y) at (0.5,4.5);
				
				\node at(1.7,4.4) { $\varsigma\in \Lambda_{xy_1}^{\lambda}$};
				
				\draw (-2,0) to [out=110,in=175] (Y) to [out=-5,in=90] (3.304,0.1);
				\filldraw  (-1.94,2.6)node[below right] {\tiny $w_{1,1}$} circle (0.02);
				\node[left] at(-1.65,3.2) {\small $\gamma_1$};
				\filldraw  (-1.31,3.7)node[below right] {\tiny $w_{0,1}$} circle (0.02);
				\filldraw[color=red]  (-1,4.004)node[below] {\tiny $\xi_1$} circle (0.02);
				
				\draw (-1.4,0.7) to [out=100,in=200] (Y);
				\filldraw  (-1.445,1.89)node[right] {\tiny $w_{1,2}$} circle (0.02);
				\filldraw  (-1.355,2.4)node[right] {\tiny $w_{2,1}$} circle (0.02);
				\node[left] at(-1.13,3.02) {\small $\gamma_2$};
				\filldraw  (-0.825,3.53)node[below right] {\tiny $w_{0,2}$} circle (0.02);
				\filldraw[color=red]  (-0.55,3.85)node[below] {\tiny $\xi_2$} circle (0.02);
				
				\draw (-0.7,0.745) to [out=90,in=220] (Y);
				\filldraw  (-0.687,1.695)node[right] {\tiny $w_{2,2}$} circle (0.02);
				\filldraw  (-0.65,2.21)node[right] {\tiny $w_{3,1}$} circle (0.02);
				\node[left] at(-0.54,2.87) {\small $\gamma_3$};
				\filldraw  (-0.342,3.41)node[below] {\tiny $\;\;\;\;\;\;w_{0,3}$} circle (0.02);
				\filldraw[color=red]  (-0.15,3.77)node[below] {\tiny $\;\;\xi_3$} circle (0.02);
				
				\draw (-0.25,0.392) to [out=88,in=248] (Y);
				\filldraw  (-0.188,1.61)node[right] {\tiny $w_{i\!-\!2,2}$} circle (0.02);
				\filldraw  (-0.135,2.13)node[right] {\tiny $w_{i-1,1}$} circle (0.02);
				\node[left] at(-0.04,2.79) {\small $\gamma_{i-1}$};
				\filldraw  (0.105,3.35)node[below] {\tiny $\;\;\;\;\;\;\;\;\;\;\;w_{0,i-1}$} circle (0.02);
				\filldraw[color=red]  (0.21,3.703)node[below, xshift=0.1cm] {\tiny $\;\;\;\xi_{i-1}$} circle (0.02);
				
				\draw (0.2,-0.265) to [out=80,in=-90] (Y);
				\filldraw  (0.44,1.58)node[right] {\tiny $w_{i-1,2}$} circle (0.02);
				\filldraw  (0.47,2.1)node[right] {\tiny $w_{i,1}$} circle (0.02);
				\node[left] at(0.54,2.77) {\small $\gamma_{i}$};
				\filldraw  (0.5,3.33)node[below] {\tiny $\;\;\;\;\;\;\;\;\;w_{0,i}$} circle (0.02);
				\filldraw[color=red]  (0.5,3.7)node[below] {\tiny $\;\;\;\;\;\;\;\xi_{i}$} circle (0.02);
				
				\draw (0.85,-0.622) to [out=85,in=-70]coordinate[pos=0.41] (wi2) coordinate[pos=0.52] (wi+11) node [left, pos=0.656, xshift=0.1cm]{\small $\gamma_{i\!+\!1}$} coordinate[pos=0.78] (w0i+1) coordinate[pos=0.85] (xii+1)  (Y);
				\filldraw  (wi2)node[right, xshift=-0.25cm] {\tiny $\;\;\;w_{i,2}$} circle (0.02);
				\filldraw  (wi+11)node[right, xshift=-0.25cm] {\tiny $\;\;\;w_{i\!+\!1,1}$} circle (0.02);
				\filldraw  (w0i+1)node[below right, xshift=-0.3cm] {\tiny $\;\;\;w_{0,i\!+\!1}$} circle (0.02);
				\filldraw[red]  (xii+1)node[below right , xshift=-0.25cm] {\tiny $\;\;\;\xi_{i\!+\!1}$} circle (0.02);
				
				\draw (1.45,-0.695) to [out=85,in=-45] (Y);
				\filldraw[xshift=0.652cm]  (0.965,1.605)node[right] {\tiny $w_{m\!-\!2,2}$} circle (0.02);
				\filldraw[xshift=0.625cm]  (0.965,2.13)node[right] {\tiny $w_{m\!-\!1\!,1}$} circle (0.02);
				\node[left, xshift=1.43cm] at(0.98,2.75) {\small $\gamma_{m\!-\!1}$};
				\filldraw[xshift=0.465cm]  (0.825,3.36)node[below right] {\tiny $w_{0\!,m\!-\!1}$} circle (0.02);
				\filldraw[color=red]  (1.115,3.708)node[below right, xshift=0.08cm, yshift=0.1cm] {\tiny $\xi_{m\!-\!1}$} circle (0.02);
				
				\filldraw[xshift=0.439cm]  (2.67,1.8)node[left] {\tiny $w_{m\!-\!1,2}$} circle (0.02);
				\node[left,xshift=0.8cm] at(2.4,2.9) {\small $\gamma_{m}$};
				\filldraw[xshift=0.23cm]  (2.1,3.5)node[below left] {\tiny $w_{0,m}$} circle (0.02);
				\filldraw[color=red]  (1.96,3.89)node[below] {\tiny $\xi_m\;\;$} circle (0.02);
				
				\draw[color=red] (-1,4.004) to [out=30,in=110] (-0.55,3.85);
				\draw[color=red] (-0.55,3.85) to [out=30,in=110] (-0.15,3.77);
				\draw[dashed, color=red] (-0.15,3.77) to [out=30,in=110] (0.21,3.703);
				\draw[color=red] (0.21,3.703) to [out=45,in=118] (0.5,3.7);
				\draw[color=red] (0.5,3.7) to [out=45,in=118] (xii+1);
				\draw[dashed, color=red] (xii+1) to [out=45,in=118] (1.115,3.708);
				\draw[color=red] (1.115,3.708) to [out=60,in=145] (1.96,3.89);
				
				\node[color=red] at(-0.7,4.05) {$\beta_1$};
				\node[color=red] at(-0.3,3.95) {$\beta_2$};
				\node[color=red] at(0.35,3.85) {$\beta_{i-1}$};
				\node[color=red] at(0.6,3.85) {$\beta_{i}$};
				\node[color=red] at(1.35,4.05) {$\beta_0$};
			\end{tikzpicture}
		\end{center}
		\caption{Illustration for the construction of $\zeta$-sequence} \label{fig3}
	\end{figure}
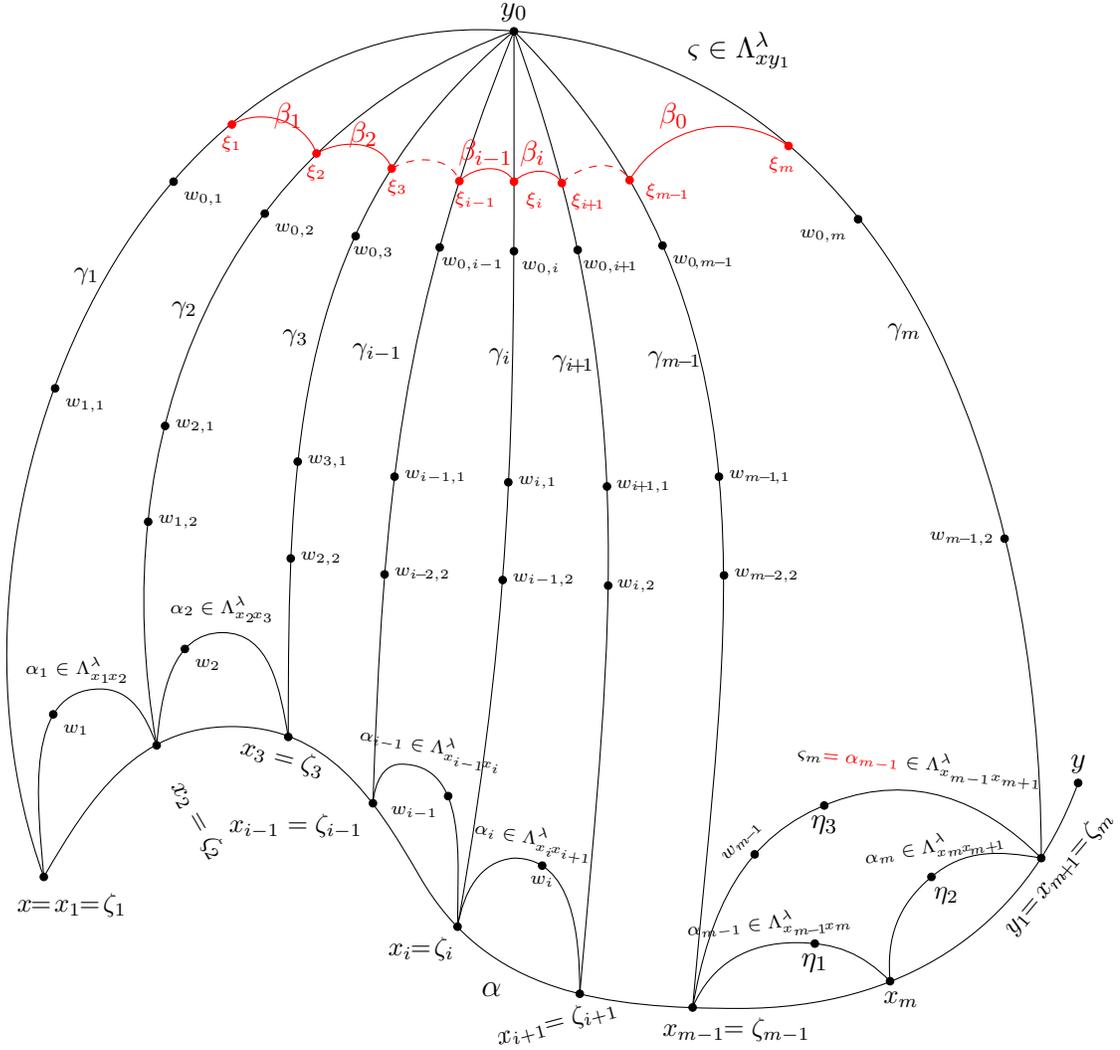

	\begin{proof}
		By Proposition \ref{11-19-1}$(P_3)$,  $D$ is $\delta$-Gromov hyperbolic.
		It follows from \eqref{22-11-16-1} and (\ref{10-3-2}) that the assumptions in Lemma \ref{24-4-22-3} are satisfied for any $\lambda$-curve $\varsigma\in \Lambda_{x_1y_1}^{\lambda}$ and $\alpha$. Then it follows that there is $x_2$ in $\alpha[x_1,y_1]$ such that for any $\lambda$-curve $\alpha_1\in \Lambda_{x_1x_2}^{\lambda}$, it holds
		\be\label{23-07-22-3}
		\mu_6^{\frac{1}{4}}\ell(\alpha)\leq \ell(\alpha_1)\leq \mu_6^{\frac{1}{2}}\ell(\alpha),
		\ee
		and for any $z\in \alpha[x_2,y_1]$, there is a $\lambda$-curve $\varsigma_1\in \Lambda_{x_1z}^{\lambda}$ such that
		$$\ell(\varsigma_1)\geq \mu_6^{\frac{1}{4}}\ell(\alpha).$$
		
		If there is an element in $\varsigma_{21}\in \Lambda_{x_2y_1}^\lambda$ such that
		$\ell(\varsigma_{21})\leq \mu_6^{\frac{1}{2}}\ell(\alpha),$ then we take $m=2$, $x_3=y_1$ and $\alpha_2=\varsigma_{21}$. Then $\{x_1,x_2,x_3\}$ and $\alpha_2$ are the desired points and $\lambda$-curve satisfying all the requirements of Proposition \ref{24-5-7-3}.
		
		If for any $\lambda$-curve $\varsigma_{21}\in \Lambda_{x_2y_1}^{\lambda}$,
		$\ell(\varsigma_{21})>\mu_6^{\frac{1}{2}}\ell(\alpha)$, then
		the similar reasoning as above ensures that there is $x_3\in\alpha[x_2,y_1]$ such that for any $\lambda$-curve $\alpha_2\in \Lambda_{x_2x_3}^{\lambda}$,
		$$\mu_6^{\frac{1}{4}}\ell(\alpha)\leq \ell(\alpha_2)\leq \mu_6^{\frac{1}{2}}\ell(\alpha),$$
		and for any $z\in \alpha[x_3,y_1]$, and for each $\lambda$-curve $\varsigma_2\in \Lambda_{x_2z}^{\lambda}$,
		$\ell(\varsigma_2)\geq \mu_6^{\frac{1}{4}}\ell(\alpha).$
		
		If there is an element $\varsigma_{31}\in \Lambda_{x_2y_1}^\lambda$ such that
		$\ell(\varsigma_{31})\leq \mu_6^{\frac{1}{2}}\ell(\alpha),$ then we take $m=3$, $x_4=y_1$ and $\alpha_3=\varsigma_{31}$. Then $\{x_1,x_2,x_3,x_4\}$ and $\alpha_3$ are the desired points and $\lambda$-curve satisfying all the requirements of Proposition \ref{24-5-7-3}.
		
		If for any $\lambda$-curve $\varsigma_{31}\in \Lambda_{x_3y_1}^{\lambda}$,
		$\ell(\varsigma_{31})>\mu_6^{\frac{1}{2}}\ell(\alpha)$, then
		the similar reasoning as above ensures that there is $x_4\in\alpha[x_3,y_1]$ such that for any $\lambda$-curve $\alpha_3\in \Lambda_{x_3x_4}^{\lambda}$,
		$$\mu_6^{\frac{1}{4}}\ell(\alpha)\leq \ell(\alpha_3)\leq \mu_6^{\frac{1}{2}}\ell(\alpha).$$
		
		As $\alpha$ is rectifiable, repeating the above procedure finitely many times, we obtain that there is a sequence of points $\{x_i\}_{i=1}^{m+1}$ on $\alpha$ and a $\lambda$-curve $\alpha_m\in \Lambda_{x_mx_{m+1}}^\lambda$ satisfying all the requirements of Proposition \ref{24-5-7-3}.
	\end{proof}
	
	\blem\label{lem2-eq-6} There is an element $\varsigma_m\in\Lambda_{x_{m-1}x_{m+1}}^\lambda$ such that
	$$\mu_6^{\frac{1}{4}}\ell(\alpha)<\ell(\varsigma_m)\leq \mu_6^{\frac{3}{4}}\ell(\alpha).$$
	\elem
	\bpf By Proposition \ref{24-5-7-3}$(b)$, there is an element $\varsigma_m\in\Lambda_{x_{m-1}x_{m+1}}^\lambda$ such that
	$\ell(\varsigma_m)>\mu_6^{\frac{1}{4}}\ell(\alpha).$
	Next, we prove that the $\lambda$-curve $\varsigma_m$ simultaneously satisfies
	\be\label{24-4-12-10}
	\ell(\varsigma_m)\leq \mu_6^{\frac{3}{4}}\ell(\alpha).
	\ee
	
	Fix a $\lambda$-curves $\alpha_{m-1}\in \Lambda_{x_{m-1}x_m}^{\lambda}$ and let $\alpha_m\in \gamma_{x_mx_{m+1}}^{\lambda}$ be given by Proposition \ref{24-5-7-3}$(c)$. Since $(D,k_D)$ is a $(\mu_0,\rho_0)$-Rips space by Proposition \ref{11-19-1}$(P_4)$,  we see from Lemma \ref{lem-8-12} that there exist $\eta_1\in \alpha_{m-1}$, $\eta_2\in \alpha_m$ and $\eta_3\in \varsigma_m$ such that for any $i\not=j\in\{1,2,3\}$ and  $\mu_1=3\mu_0$,
	\beq\label{lem2-eq-11}
	k_D(\eta_i,\eta_j)\leq \mu_1.
	\eeq 
	
	Next, we estimate $\ell(\varsigma_m[x_{m+1},\eta_3])$ and $\ell(\varsigma_m[x_{m-1},\eta_3])$.
	To this end, we divide our discussions into two cases:
	$\min\{k_D(x_{m+1},\eta_2),$ $k_D(x_{m+1},\eta_3)\}> 1$ and $\min\{k_D(x_{m+1},\eta_2),k_D(x_{m+1},\eta_3)\}\leq 1.$ In the former case, we consider the points $\eta_2,$ $\eta_3$ and $x_{m+1}$.
	Since Proposition \ref{11-19-1}$(P_3)$ implies that $D$ is $\delta$-Gromov hyperbolic, we infer from \eqref{lem2-eq-11},  Lemma \ref{lem-2.4} and Proposition \ref{24-5-7-3}$(c)$ that
	\beqq
	\ell(\varsigma_m[x_{m+1},\eta_3])\leq 4e^{25\mu_2}\ell(\alpha_m[x_{m+1},\eta_2])\leq 4\mu_6^{\frac{1}{2}}e^{25\mu_2}\ell(\alpha),\eeqq
	where in the first inequality, we used the fact $\mu_1\leq \mu_2$.
	
	In the latter case, that is, $\min\{k_D(x_{m+1},\eta_2),k_D(x_{m+1},\eta_3)\}\leq 1,$ it follows from \eqref{lem2-eq-11} that
	$$k_D(x_{m+1},\eta_3)\leq \min\{k_D(x_{m+1},\eta_2),k_D(x_{m+1},\eta_3)\}+k_D(\eta_2,\eta_3)\leq 1+\mu_1.$$
	Since $\varsigma_m[x_{m+1},\eta_3]$ is a $2\lambda$-curve by Lemma \ref{lem-2.3}, we infer from (\ref{(2.1)}) and Proposition \ref{11-19-1}$(P_1)$ that
	\beqq
		\log\Big(1+\frac{\ell(\varsigma_m[x_{m+1},\eta_3])}{d_D(x_{m+1})}\Big) \leq  \ell_{k}(\varsigma_m[x_{m+1},\eta_3])\leq \frac{9}{8}k_D(x_{m+1},\eta_3)
		\leq  \frac{9}{8}(1+\mu_1),
	\eeqq
	which, together with \eqref{22-11-16-1}, implies that
	\beqq\label{lem2-eq-ne-11}
	\ell(\varsigma_m[x_{m+1},\eta_3])\leq e^{\frac{9}{8}(1+\mu_1)}d_D(x_{m+1})\leq 2e^{\frac{9}{8}(1+\mu_1)}\ell(\alpha).
	\eeqq
	
	In either cases, we obtain
	\be\label{24-4-12-11} \ell(\varsigma_m[x_{m+1},\eta_3])\leq 4\mu_6^{\frac{1}{2}}e^{25\mu_2}\ell(\alpha).\ee
	Similarly, repeating the above arguments with $x_{m+1}$, $\eta_2$ replaced by $x_{m-1}$, $\eta_1$, we obtain
	\be\label{24-4-12-12}\ell(\varsigma_m[x_{m-1},\eta_3])\leq  4\mu_6^{\frac{1}{2}}e^{25\mu_2}\ell(\alpha).\ee
	
	Since $8\mu_6^{\frac{1}{2}}e^{25\mu_2}\leq \mu_6^{\frac{3}{4}}$, \eqref{24-4-12-10} follows from \eqref{24-4-12-11} and \eqref{24-4-12-12}, and hence, the proof of lemma is complete.
	\epf
	
	For convenience, we make the following notational conventions. 
	\begin{enumerate}
		\item\label{24-5-7-1}
		for each $i\in\{1,2,\cdots,m-1\}$, let $\zeta_i=x_i$ and $\zeta_m=x_{m+1}$.
		\item\label{24-5-7-2}
		For each $i\in\{1,2,\cdots,m-2\}$, we fix an element $\alpha_i\in\Lambda_{\zeta_i\zeta_{i+1}}^\lambda$. Moreover, let $\alpha_{m-1}=\varsigma_m\in \Lambda_{x_{m-1}x_{m+1}}^{\lambda}$ be given by Lemma \ref{lem2-eq-6}, and $\alpha_m\in \Lambda_{x_mx_{m+1}}^{\lambda}$ be given by Proposition \ref{24-5-7-3}$(c)$.
	\end{enumerate}
	Then we see from Proposition \ref{24-5-7-3}$(a)$ and Lemma \ref{lem2-eq-6} that for each $i\in\{1,2,\cdots,m-2\}$,
	\beq\label{lem2-eq-13}
	\mu_6^{\frac{1}{4}}\ell(\alpha)<\ell(\alpha_{i})\leq \mu_6^{\frac{1}{2}}\ell(\alpha)\;\;\mbox{and}\;\;\mu_6^{\frac{1}{4}}\ell(\alpha)<\ell(\alpha_{m-1})\leq \mu_6^{\frac{3}{4}}\ell(\alpha),\eeq
	and thus, we infer from \eqref{(2.1)} and \eqref{22-11-16-1} that for each $i\in\{1,2,\cdots,m-1\}$,
	\be\label{2022-10-24-2}\ell_{k}(\alpha_i)\geq \log\Big(1+\frac{\ell(\alpha_i)}{d_D(\zeta_i)}\Big)\geq \log\Big(1+\frac{\ell(\alpha_i)}{2\ell(\alpha)}\Big)> \frac{1}{4}\log \frac{\mu_6}{16}.\ee

	\subsection{$\lambda$-curves associated to $\zeta$-sequence}\label{sub-5.2}
	
	Fix any $\lambda$-curve $\varsigma\in \Lambda_{x_1y_1}^\lambda$ and write  $\gamma_1=\varsigma[\zeta_1,x_0]$ and $\gamma_{m}=\varsigma[y_1,x_0]$. For $i\in \{2,3,\cdots,$ $m-1\}$, let
	$\gamma_i\in \Lambda_{\zeta_ix_0}^{\lambda}$ be an element; see Figure \ref{fig3}.
	
	Since, by Proposition \ref{11-19-1}$(P_4)$, $(D,k_D)$ is a $(\mu_0,\rho_0)$-Rips space, it follows from Lemma \ref{lem-8-12} that for each $i\in \{1,2,\cdots,m-1\}$, there exist three points $w_i\in \alpha_i$, $w_{i,1}\in \gamma_i$, $w_{i,2}\in \gamma_{i+1}$ such that
	\beq\label{lem2-eq-9}
	k_D(w_i,w_{i,1})\leq \mu_1,\;\;k_D(w_{i,1},w_{i,2})\leq \mu_1\;\;{\rm and}\;\; k_D(w_i,w_{i,2})\leq \mu_1.
	\eeq
	
	\blem\label{2022-10-24-4}
	For each $i\in \{1,2,\cdots,m-1\}$, we have
	$$\max\{\ell(\gamma_i[\zeta_i,w_{i,1}]),\ell(\gamma_{i+1}[\zeta_{i+1},w_{i,2}])\}\leq 4\mu_6^{\frac{3}{4}}e^{25\mu_2}\ell(\alpha).$$
	\elem
	\bpf
	Without loss of generality, we assume $\max\{\ell(\gamma_i[\zeta_i,w_{i,1}]),$ $\ell(\gamma_{i+1}[\zeta_{i+1},w_{i,2}])\}=\ell(\gamma_i[\zeta_i,w_{i,1}]),$ since the proof of other case is similar.
	
	We first assume that $\min\{k_D(\zeta_i,w_i),$ $k_D(\zeta_i,w_{i,1})\}> 1.$
	It follows from \eqref{lem2-eq-9}, Proposition \ref{11-19-1}$(P_1)$, and Lemma \ref{lem-2.4} (applied to $\zeta_i$, $w_i$ and $w_{i,1}$) that
	$$\ell(\gamma_i[\zeta_i,w_{i,1}])\leq  4e^{25\mu_2}\ell(\alpha_i[\zeta_i,w_i])\stackrel{\eqref{lem2-eq-13}}{\leq} 4\mu_6^{\frac{3}{4}}e^{25\mu_2}\ell(\alpha).$$
	
	In the remaining case, that is, $\min\{k_D(\zeta_i,w_i),k_D(\zeta_i,w_{i,1})\}\leq 1,$ it follows from \eqref{lem2-eq-9}, 
	\eqref{(2.1)} and Proposition \ref{11-19-1}$(P_1)$ that
	$$
	\begin{aligned}
		\log\Big(1+\frac{\ell(\gamma_i[\zeta_i,w_{i,1}])}{d_D(\zeta_i)}\Big)&\leq\ell_{k}(\gamma_i[\zeta_i,w_{i,1}])\leq \frac{9}{8}k_D(\zeta_i,w_{i,1})\\
		&\leq \frac{9}{8}\min\{k_D(\zeta_i,w_i),k_D(\zeta_i,w_{i,1})\}+k_D(w_i,w_{i,1})
		\leq \frac{9}{8}(1+\mu_1),
	\end{aligned}
	$$
	which, together with \eqref{22-11-16-1}, implies that
	$$\ell(\gamma_i[\zeta_i,w_{i,1}])\leq e^{\frac{9}{8}(1+\mu_1)}d_D(\zeta_i)\leq 2e^{\frac{9}{8}(1+\mu_1)}\ell(\alpha)<4\mu_6^{\frac{3}{4}}e^{25\mu_2}\ell(\alpha).$$
	In either cases, the proof is complete.
	\epf
	
	Now, we know from \eqref{10-3-2} that there exists $w_{0,i}\in\gamma_i$ such that
	\be\label{2022-09-30-1}\ell(\gamma_i[\zeta_i,w_{0,i}])=\mu_6\ell(\alpha)>4\mu_6^{\frac{3}{4}}e^{25\mu_2}\ell(\alpha).\ee
	Furthermore, we infer from \eqref{2022-09-30-1} and Lemma \ref{2022-10-24-4} that
	$w_{0,1}\in\gamma_1[w_{1,1},x_0]$, $w_{0,m}\in\gamma_m[w_{m-1,2},x_0]$,
	and for each $i\in\{2,\cdots,m-1\}$,
	$$w_{0,i}\;\in \; \gamma_i[w_{i,1},x_0]\cap \gamma_i[w_{i-1,2},x_0].$$
	
	\blem\label{2022-09-30-2} For each $i\in\{1,\cdots,m-1\}$,
	$$\min\{k_D(w_{i,1},w_{0,i}),\;\;k_D(w_{i,2},w_{0,{i+1}})\} \geq \frac{1}{5}\log \mu_6.$$
	\elem
	\bpf
	Let $i\in\{1,\cdots,m-1\}$. Since Lemma \ref{2022-10-24-4} and \eqref{22-11-16-1} lead to
	$$d_D(w_{i,1})\leq \ell(\gamma_i[\zeta_i,w_{i,1}])+d_D(\zeta_i)< 2(1+2\mu_6^{\frac{3}{4}}e^{25\mu_2})\ell(\alpha),$$
	and since Lemma \ref{2022-10-24-4} and \eqref{2022-09-30-1} give
	$$\ell(\gamma_i[w_{i,1},w_{0,i}])= \ell(\gamma_i[\zeta_i,w_{0,i}])-\ell(\gamma_i[\zeta_i,w_{i,1}])\geq (\mu_6-4\mu_6^{\frac{3}{4}}e^{25\mu_2})\ell(\alpha),$$
	we get from Proposition \ref{11-19-1}$(P_1)$, \eqref{(2.1)} and the above two estimates that
	\beqq
	k_D(w_{i,1},w_{0,i})\geq \frac{8}{9}\ell_{k}(\gamma_i[w_{i,1},w_{0,i}])\geq\frac{8}{9}\log\Big(1+\frac{\ell(\gamma_i[w_{i,1},w_{0,i}])}{d_D(w_{i,1})}\Big)\geq \frac{1}{5}\log \mu_6.
	\eeqq
	
	Arguing similarly as above, replacing $\gamma_i$ by $\gamma_{i+1}$, we obtain $$k_D(w_{i,2},w_{0,{i+1}}) \geq \frac{1}{5}\log \mu_6.$$
	This completes the proof of lemma.
	\epf
	
	\blem\label{23-07-23-4} For each $i\in\{1,\cdots,m-1\}$,
	$$\min\{k_D(x_0,w_{i,1})-k_D(w_{0,i},x_0),\; k_D(x_0,w_{i,2})-k_D(w_{0,i+1},x_0)\}  \geq \frac{1}{5}\log \mu_6-\frac{1}{2}.$$
	\elem
	\bpf By Lemmas \ref{lem-2.2}, \ref{lem-2.3} and \ref{2022-09-30-2}, we see that for each $i\in\{1,\cdots,m-1\}$,
	\beqq
	k_D(x_0,w_{i,1})&\geq& \ell_k(\gamma_i[w_{i,1},x_0])-\frac{1}{2}\\ \nonumber&=&\ell_k(\gamma_i[w_{i,1},w_{0,i}])+\ell_k(\gamma_i[w_{0,i},x_0])-\frac{1}{2}
	\\ \nonumber&\geq& k_D(w_{0,i},x_0)+\frac{1}{5}\log \mu_6-\frac{1}{2}.
	\eeqq
	
	Similarly, we can show that for each $i\in\{1,\cdots,m-1\}$, $$k_D(x_0,w_{i,2})\geq k_D(w_{0,i+1},x_0)+\frac{1}{5}\log \mu_6-\frac{1}{2}.$$
	This completes the proof.
	\epf
	
	Let $t_1\in\{1,\cdots,m\}$ be such that
	\be\label{11-12-1} k_D(x_0,w_{0,t_1})=\min_{i\in\{1,\cdots,m\}}\{k_D(x_0,w_{0,i})\}.\ee
	Then we have the following assertion.
	\blem\label{cl-23-07-23}
	Suppose that there exist $i\in\{1,\cdots,m\}$ and $u_i\in \gamma_i[w_{0,i},x_0]$ such that  $k_D(x_0,u_i)=k_D(x_0,w_{0,t_1})$.
	\begin{enumerate}
		\item\label{23-07-23-6}
		If $i\in\{2,\cdots,m\}$, then there exists $u_{i-1}\in \gamma_{i-1}$ such that
		$$k_D(x_0,u_{i-1})=k_D(x_0,w_{0,t_1})\;\;\mbox{and}\;\; k_D(u_i,u_{i-1})\leq \frac{3}{2}+2\mu_0+2\mu_1.$$
		
		\item\label{23-07-23-7} If $i\in\{1,\cdots,m-1\}$, then there exists $u_{i+1}\in \gamma_{i+1}$ such that
		$$k_D(x_0,u_{i+1})=k_D(x_0,w_{0,t_1})\;\;\mbox{and}\;\;  k_D(u_i,u_{i+1})\leq \frac{3}{2}+2\mu_0+2\mu_1.$$
	\end{enumerate}
	\elem
	\bpf
	We only need to prove \eqref{23-07-23-6}  since the proof of \eqref{23-07-23-7} is similar. Let $i\in\{2,\cdots,m\}$. To find the asserted point $u_{i-1}$, we need some preparation.
	
	Firstly, we determine the position of $u_i$ on $\gamma_i$. Since Lemmas \ref{lem-2.2} and \ref{lem-2.3} lead to
	$$\ell_k(\gamma_i[x_0,u_i])\leq  k_D(x_0,u_i)+\frac{1}{2},$$
	and since Lemma \ref{23-07-23-4} and \eqref{11-12-1} give
	$$k_D(x_0,w_{i,1})\geq k_D(w_{0,i},x_0)+\frac{1}{5}\log \mu_6-\frac{1}{2}\geq k_D(x_0,u_i)+\frac{1}{5}\log \mu_6-\frac{1}{2},$$
	we obtain from the above two estimates that for $i\in\{2,\cdots,m-1\}$, it holds
	$$\ell_k(\gamma_{i}[x_0,w_{i,1}])\geq k_D(x_0,w_{i,1})\geq\ell_k(\gamma_i[x_0,u_i])+\frac{1}{5}\log \mu_6-1.$$
	
	Similarly, we know that for $i\in\{2,\cdots,m\}$,
	$$\ell_k(\gamma_{i}[x_0,w_{i-1,2}])\geq \ell_k(\gamma_i[x_0,u_i])+\frac{1}{5}\log \mu_6-1.$$
	These ensure that for $i\in\{2,\cdots,m-1\}$,
	\be\label{23-07-23-8}
	u_i\;\in \; \big(\gamma_i[w_{i,1},x_0]\setminus \{w_{i,1}\}\big)\cap \big(\gamma_i[w_{i-1,2},x_0]\setminus \{w_{i-1,2}\}\big)
	\ee
	and
	\be\label{24-5-15-1}
	u_m\;\in \; \gamma_i[w_{m-1,2},x_0]\setminus \{w_{m-1,2}\}.
	\ee
	
	Secondly, we find a point on $\gamma_{i-1}$. Since Proposition \ref{11-19-1}$(P_4)$ implies that $(D,k_D)$ is a $(\mu_0,\rho_0)$-Rips space,  it follows
	from \eqref{lem2-eq-9} and Lemma \ref{lem-2.4-1} (replacing $x$, $y$ and $z$ with $x_0$, $w_{i-1,2}$ and $w_{i-1,1}$), together with \eqref{23-07-23-8} and \eqref{24-5-15-1}, that
	there exists $v_{i-1}\in\gamma_{i-1}[w_{i-1,1},x_0]$ such that
	\beq\label{7-26-2}k_D(u_i,v_{i-1})\leq \frac{1}{2}+\mu_0+\mu_1.\eeq
	Then we have
	\beq\label{23-07-23-10}
	\begin{aligned}
		k_D(x_0,u_i)&-\frac{1}{2}-\mu_0-\mu_1\leq k_D(x_0,u_i)-k_D(u_i,v_{i-1})\\
		&\leq k_D(x_0,v_{i-1})\leq k_D(x_0,u_i)+k_D(u_i,v_{i-1})\\
		&\leq k_D(x_0,u_i)+\frac{1}{2}+\mu_0+\mu_1.
	\end{aligned}
	\eeq
	
	Now, we are ready to find the required point.
	If $k_D(x_0,v_{i-1})\leq k_D(x_0,u_i)$, then we infer from Lemma \ref{23-07-23-4} and \eqref{11-12-1} that
	$$ k_D(x_0,w_{i-1,1}) \geq k_D(x_0,w_{0,i-1})+\frac{1}{5}\log \mu_6-\frac{1}{2}\geq k_D(x_0,u_i)+\frac{1}{5}\log \mu_6-\frac{1}{2}> k_D(x_0,v_{i-1}),$$ from which we know that
	there exists $u_{i-1}\in\gamma_{i-1}[w_{i-1,1},v_{i-1}]$ such that
	$$k_D(x_0,u_{i-1})=k_D(x_0,u_i).$$ Then it follows from Lemmas \ref{lem-2.2} and \ref{lem-2.3}, and \eqref{23-07-23-10} that
	\begin{eqnarray*} k_D(u_{i-1},v_{i-1})&\leq& \ell_k(\gamma_{i-1}[x_0,u_{i-1}])-\ell_k(\gamma_{i-1}[x_0,v_{i-1}])\\ \nonumber&\leq& k_D(x_0,u_{i-1})-k_D(x_0,v_{i-1})+\frac{1}{2}\\ \nonumber
		&\leq& 1+\mu_0+\mu_1,
	\end{eqnarray*} which, together with \eqref{7-26-2}, gives
	$$k_D(u_{i-1},u_i)\leq k_D(u_{i-1},v_{i-1})+k_D(v_{i-1},u_i)\leq \frac{3}{2}+2\mu_0+2\mu_1.$$
	
	If $k_D(x_0,v_{i-1})>k_D(x_0,u_i)$, then select $u_{i-1}\in\gamma_{i-1}[x_0,v_{i-1}]$ such that
	$$k_D(x_0,u_{i-1})=k_D(x_0,u_i).$$ Then a similar computation as above gives
	$$k_D(u_{i-1},u_i)\leq k_D(u_{i-1},v_{i-1})+k_D(v_{i-1},u_i)\leq \frac{3}{2}+2\mu_0+2\mu_1.$$
	The proof of the lemma is thus complete.
	\epf

	\subsection{$\xi$-sequence on the associated $\lambda$-curves}\label{sub-5.3}
	In this section, we shall find a new sequence of points. By the choice of $w_{0,t_1}$ in \eqref{11-12-1} and Lemma \ref{cl-23-07-23}, we know that
	for each $i\in\{1,\cdots,m\}$, there exists $\xi_i\in\gamma_i[w_{0,i},x_0]$ (see Figure \ref{fig3}) such that
	\be\label{2022-09-30-3}
	k_D(x_0,\xi_i)=k_D(x_0,w_{0,t_1}),
	\ee
	and for every $i\in\{1,\cdots,m-1\}$,
	\be\label{2022-09-30-4}
	k_D(\xi_i,\xi_{i+1})\leq \frac{3}{2}+2\mu_0+2\mu_1.
	\ee

	\blem\label{cl3.2-1} $(i)$
	For each $i\in\{1,\cdots,m\}$ and all $z\in\gamma_i[\zeta_i,\xi_i]$, it holds
	\beqq\label{2022-10-11-1}
	k_D(x_0,z)\geq \log \frac{\mu_5-2}{2+\mu_6}-\frac{1}{2}.
	\eeqq
	
	$(ii)$ Let $i\in\{1,\cdots,m-1\}$. For any $\lambda$-curve $\beta_i\in \Lambda_{\xi_i\xi_{i+1}}^{\lambda}$, it holds
	\beqq\label{23-07-24-3}
	\ell_k(\beta_i)\leq  \frac{7}{4}+2\mu_0+2\mu_1.
	\eeqq
	
	$(iii)$ Let $i$, $j\in\{1,\cdots,m-1\}$. For any $\lambda$-curve $\beta_i\in \Lambda_{\xi_i\xi_{i+1}}^{\lambda}$ and any $u\in\beta_i$, it holds
	\beqq\label{23-12-01-1}
	k_D(x_0,\xi_j)-\frac{7}{4}-2\mu_0-2\mu_1\leq k_D(x_0,u)\leq k_D(x_0,\xi_j)+ \frac{7}{4}+2\mu_0+2\mu_1.
	\eeqq
	\elem
	\bpf
	(i). We get from \eqref{22-11-16-1} and \eqref{2022-09-30-1} that
	\be\label{2022-12-13-1}
	d_D(w_{0,t_1})\leq d_D(\zeta_{t_1})+
	\ell(\gamma_{t_1}[\zeta_{t_1},w_{0,t_1}])\leq (2+\mu_6)\ell(\alpha).\ee
	Let $z\in\gamma_i[\zeta_i,\xi_i]$. Then we obtain from Lemmas \ref{lem-2.2} and \ref{lem-2.3} that
	\beqq
	\begin{aligned}
		k_D(x_0,z) &\geq \ell_{k}(\gamma_i[x_0,z])-\frac{1}{2}\geq\ell_{k}(\gamma_i[x_0,\xi_i])-\frac{1}{2}\\
		&\geq  k_D(x_0,\xi_i)-\frac{1}{2}\stackrel{\eqref{2022-09-30-3}}{=}k_D(x_0,w_{0,t_1})-\frac{1}{2}
		\stackrel{\eqref{(2.2)}}{\geq}  \log \frac{d_D(x_0)}{d_D(w_{0,t_1})}-\frac{1}{2}\\
		&\stackrel{\eqref{10-3-2} + \eqref{2022-12-13-1}}{\geq}  \log \frac{\mu_5-2}{2+\mu_6}-\frac{1}{2}.
	\end{aligned}
	\eeqq
	This proves $(i)$.
	
	(ii). Let $i\in\{1,\cdots,m-1\}$. For any $\lambda$-curve $\beta_i\in \Lambda_{\xi_i\xi_{i+1}}^{\lambda}$,
	we have by \eqref{2022-09-30-4} and Lemma \ref{lem-2.2} that
	$$\ell_k(\beta_i)\leq k_D(\xi_i,\xi_{i+1})+\frac{1}{4}\leq \frac{7}{4}+2\mu_0+2\mu_1,$$
	which proves $(ii)$.
	
	(iii). Let $i\in\{1,\cdots,m-1\}$. For any $\lambda$-curve $\beta_i\in \Lambda_{\xi_i\xi_{i+1}}^{\lambda}$, let $u\in\beta_i$. By $(ii)$, we know that
	\begin{equation}\label{24-4-30-1}
		\begin{aligned}
			k_D(x_0,\xi_i)&-\frac{7}{4}-2\mu_0-2\mu_1 \leq k_D(x_0,\xi_i)-\ell_k(\beta_{i})\\
			&\leq k_D(x_0,u)
			\leq  k_D(x_0,\xi_i)+\ell_k(\beta_{i})
			\\
			&\leq  k_D(x_0,\xi_i)+\frac{7}{4}+2\mu_0+2\mu_1.
		\end{aligned}
	\end{equation}
	Since for any $i$, $j\in\{1,\cdots,m-1\}$,  $k_D(x_0,\xi_i)=k_D(x_0,\xi_j)$ by \eqref{2022-09-30-3}, statement $(iii)$ follows immediately from \eqref{24-4-30-1} . The proof of lemma is thus complete.
	\epf

	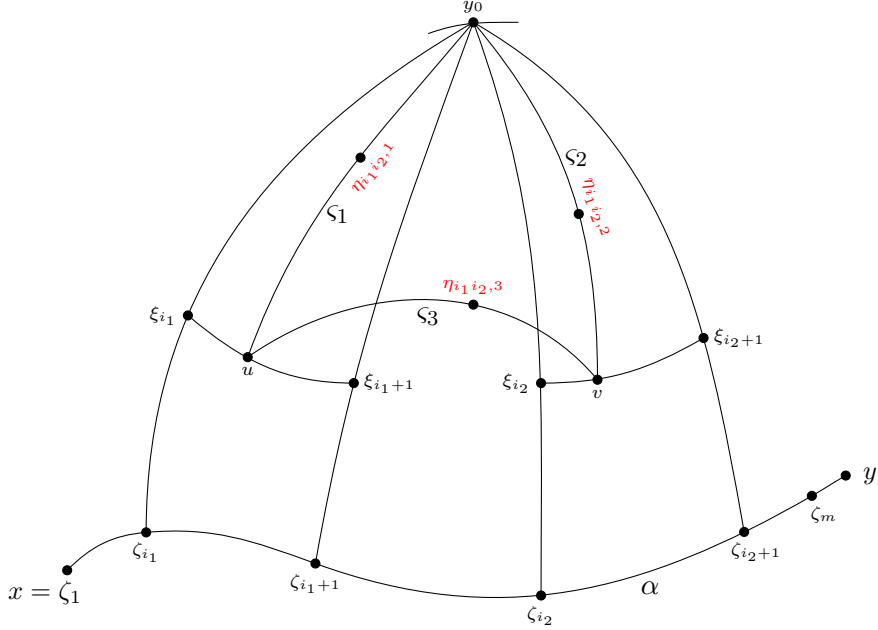
\begin{figure}[htbp]
		\begin{center}
			\begin{tikzpicture}[scale=3]
				\draw (-2.1,-0.03) to [out=45,in=185] (-1.75,0.138) to [out=5,in=160]  (-1,0) to [out=-20,in=210]  coordinate[pos=1] (zetam) node[below,pos=0.67] { $\alpha$} (1.2,0.3);
				\filldraw  (-2.1,-0.03)node[below] {\small $x=\zeta_{1}\;\;\;\;\;\;$} circle (0.02);
				\filldraw  (-1.75,0.138)node[below] {\tiny $\zeta_{i_1}$} circle (0.02);
				\filldraw  (-1,0)node[below] {\tiny $\zeta_{i_1+1}$} circle (0.02);
				\filldraw  (0,-0.142)node[below] {\tiny $\zeta_{i_2}$} circle (0.02);
				\filldraw  (0.9,0.14)node[below] {\tiny $\;\;\;\;\zeta_{i_2+1}$} circle (0.02);
				\filldraw  (zetam)node[below] {\tiny $\;\;\;\;\zeta_m$} circle (0.02);
				
				\draw (-0.5,2.35) to [out=20,in=180] (-0.1,2.4);
				\filldraw  (-0.3,2.4)node[above] {\tiny $x_0$} circle (0.02);
				
				\draw (-0.3,2.4) to [out=210,in=90] (-1.75,0.138);
				\draw (-0.3,2.4) to [out=250,in=80] (-1,0);
				\draw (-0.3,2.4) to [out=290,in=90] (0,-0.142);
				\draw (-0.3,2.4) to [out=330,in=100] (0.9,0.14);
				\filldraw  (-1.565,1.1)node[left] {\tiny $\xi_{i_1}$} circle (0.02);
				\filldraw  (-0.83,0.8)node[right] {\tiny $\xi_{i_1+1}$} circle (0.02);
				\filldraw  (0,0.8)node[left] {\tiny $\xi_{i_2}$} circle (0.02);
				\filldraw  (0.72,1)node[right] {\tiny $\xi_{i_2+1}$} circle (0.02);
				
				\draw (-1.565,1.1) to [out=-40,in=180] (-0.83,0.8);
				\draw (0,0.8) to [out=0,in=210] (0.72,1);
				
				\filldraw  (-1.3,0.915)node[below] {\tiny $u$} circle (0.02);
				\filldraw  (0.25,0.815)node[below] {\tiny $v$} circle (0.02);
				
				\draw (-1.3,0.915) to [out=35,in=130] (0.25,0.815);
				\filldraw  (-0.3,1.148)node[color=red,above] {\tiny ${\eta}_{i_1i_2,3}$} circle (0.02);
				
				\draw (-0.3,2.4) to [out=230,in=70] (-1.3,0.915);
				\draw (-0.3,2.4) to [out=310,in=90] (0.25,0.815);
				\filldraw  (-0.8,1.8)node[color=red,below,rotate=55] {\tiny ${\eta}_{i_1i_2,1}$} circle (0.02);
				\filldraw  (0.167,1.55)node[color=red,above,rotate=-70] {\tiny ${\eta}_{i_1i_2,2}$} circle (0.02);
				
				\node at(-0.90,1.55) {$\varsigma_{1}$};
				\node at(0.158,1.80) {$\varsigma_{2}$};
				\node[below] at(-0.5,1.18) {$\varsigma_{3}$};
				
			\end{tikzpicture}
		\end{center}
		\caption{Illustration for the proofs of Lemmas \ref{lem-3.4-0} and \ref{lem-22-09-20}} \label{fig5}
	\end{figure}

	In the following, for each $i\in \{1,\cdots,m\}$, let $\{\xi_i\}$ be the sequence given by \eqref{2022-09-30-3}. Then for each $i\in \{1,\cdots,m-1\}$, we fix a $\lambda$-curve $\beta_i\in \Lambda_{\xi_{i}\xi_{i+1}}^\lambda$ and set $\beta_0=\bigcup_{i=1}^{m-1}\beta_i$; see Figure \ref{fig3}. Fix $u$, $v\in \beta_0$. Then there are $i_1$, $i_2\in\{1,\cdots,m-1\}$ such that $u\in\beta_{i_1}$ and $v\in\beta_{i_2}$ (see Figure \ref{fig5}).
	The following result gives a sufficient condition for $i_1$ and $i_2$ to be different.
	
	\begin{lem}\label{lem-3.4-0} Let $u\in\beta_{i_1}$ and $v\in\beta_{i_2}$, where $i_1$, $i_2\in\{1,\cdots,m-1\}$. If
		$k_D(u,v)\geq 3\mu_1$, then
		$i_1\not= i_2$.
	\end{lem}
	\bpf Suppose on the contrary that $i_1=i_2$, that is, $u,v\in \beta_{i_1}$. Then we know from Lemma \ref{cl3.2-1}$(ii)$ that
	$$k_D(u,v)\leq \ell_k(\beta_{i_1})\leq 2+2\mu_0+2\mu_1< 3\mu_1.$$
	This contradiction implies that $i_1\neq i_2$.
	\epf
	
	\begin{lem}\label{lem-22-09-20} Suppose that $i_1$, $i_2\in\{1,\cdots,m-1\}$ with $i_1<i_2$, $u\in\beta_{i_1}$ and $v\in\beta_{i_2}$. Fix $\varsigma_1\in \Lambda_{ux_0}^{\lambda}$, $\varsigma_2\in\Lambda_{vx_0}^{\lambda}$ and $\varsigma_3\in\Lambda_{uv}^{\lambda}$.  Then for each $i\in \{1,2,3\}$, there exists a point ${\eta}_{i_1i_2,i}\in \varsigma_i$ such that
		\begin{enumerate}
			\item\label{lem-22-20-1}
			for any $r\not=s\in\{1,2,3\}$, $$k_D({\eta}_{i_1i_2,r},{\eta}_{i_1i_2,s})\leq \mu_1;$$
			\item\label{lem-22-20-2}
			$k_D(u,{\eta}_{i_1i_2,1})-\frac{15}{4}-\mu_0-6\mu_1\leq k_D(v,{\eta}_{i_1i_2,2})\leq k_D(u,{\eta}_{i_1i_2,1})+\frac{15}{4}+\mu_0+6\mu_1$;
			\item\label{lem-22-09-20-3}
			$\frac{1}{2}k_D(u,v)-\frac{15}{8}-\frac{1}{2}\mu_0-4\mu_1\leq k_D(u,{\eta}_{i_1i_2,3})\leq \frac{1}{2}k_D(u,v)+2+\frac{1}{2}\mu_0+4\mu_1$;
			\item\label{lem-22-09-20-4}
			$\frac{1}{2}k_D(u,v)-\frac{15}{8}-\frac{1}{2}\mu_0-4\mu_1\leq k_D(v,{\eta}_{i_1i_2,3})\leq \frac{1}{2}k_D(u,v)+2+\frac{1}{2}\mu_0+4\mu_1$.
		\end{enumerate}
	\end{lem}
	\bpf $(1)$. Since Proposition \ref{11-19-1}$(P_4)$ implies that $(D,k_D)$ is a $(\mu_0,\rho_0)$-Rips space, we see from Lemma \ref{lem-8-12} that for each $i\in \{1,2,3\}$, there exists ${\eta}_{i_1i_2,i}\in \varsigma_i$ so that statement (1) holds.
	\medskip
	
	$(2)$.
	By Lemma \ref{cl3.2-1}$(ii)$, we have
	\beqq
	\begin{aligned}
		k_D(x_0,\xi_{i_1+1})&-\frac{7}{4}-2\mu_0-2\mu_1 \leq k_D(x_0,\xi_{i_1+1})-\ell_k(\beta_{i_1})\\ &\leq  k_D(x_0,u)
		\leq
		k_D(x_0,\xi_{i_1+1})+\ell_k(\beta_{i_1})\\ &\leq  k_D(x_0,\xi_{i_1+1})+\frac{7}{4}+2\mu_0+2\mu_1,
	\end{aligned}
	\eeqq
	and
	\beqq 
	\begin{aligned}k_D(x_0,\xi_{i_2})&-\frac{7}{4}-2\mu_0-2\mu_1 \leq k_D(x_0,\xi_{i_2})-\ell_k(\beta_{i_2})\\ &\leq k_D(x_0,v)\leq
		k_D(x_0,\xi_{i_2})+\ell_k(\beta_{i_2})\\ &\leq k_D(x_0,\xi_{i_2})+\frac{7}{4}+2\mu_0+2\mu_1.\end{aligned}
	\eeqq	
	These, together with \eqref{2022-09-30-3}, give
	\be\label{2022-11-03-1}
	k_D(x_0,u)-\frac{7}{2}-\mu_0-5\mu_1\leq k_D(x_0,v)\leq k_D(x_0,u)+\frac{7}{2}+\mu_0+5\mu_1.\ee
	Moreover, we know from statement (1) of the lemma that
	\be\label{24-5-15-3}k_D(x_0,{\eta}_{i_1i_2,1})-\mu_1\leq k_D(x_0,{\eta}_{i_1i_2,2})\leq k_D(x_0,{\eta}_{i_1i_2,1})+\mu_1.\ee
	
	Since $\varsigma_1\in \Lambda_{ux_0}^{\lambda}$, we have
	\beqq
	\begin{aligned}k_D(u,{\eta}_{i_1i_2,1})&\leq \ell_{k}(\varsigma_1[u,{\eta}_{i_1i_2,1}]) =\ell_{k}(\varsigma_1)-\ell_{k}(\varsigma_1[x_0,{\eta}_{i_1i_2,1}])
		\\
		&\stackrel{\text{Lemma }\ref{lem-2.2} + \eqref{24-5-15-3}}{\leq} k_D(x_0,u)-k_D(x_0,{\eta}_{i_1i_2,2})+\frac{1}{4}+\mu_1\\
		&\stackrel{\eqref{2022-11-03-1}}{\leq} k_D(x_0,v)-k_D(x_0,{\eta}_{i_1i_2,2})+\frac{15}{4}+\mu_0+6\mu_1\\
		&\leq k_D(v,{\eta}_{i_1i_2,2})+\frac{15}{4}+\mu_0+6\mu_1.
	\end{aligned}
	\eeqq
	Similarly, we have
	$$ k_D(v,{\eta}_{i_1i_2,2})\leq k_D(u,{\eta}_{i_1i_2,1})+\frac{15}{4}+\mu_0+6\mu_1.$$
	These prove the second statement \eqref{lem-22-20-2}.
	\medskip
	
	$(3)$. To prove the third statement of the lemma, we note first that
	\beqq
	\begin{aligned}
		k_D(u,{\eta}_{i_1i_2,3})
		&\leq  k_D(u,{\eta}_{i_1i_2,1})+\mu_1& \;\;\;\;\;\;\;\;&\mbox{(by $\eqref{lem-22-20-1}$ of the lemma)}    \\
		&\leq
		k_D(v,{\eta}_{i_1i_2,2})+\frac{15}{4}+\mu_0+7\mu_1& \;\;\;\;\;\;\;&\mbox{(by $\eqref{lem-22-20-2}$ of the lemma)}    \\ &\leq
		k_D(v,{\eta}_{i_1i_2,3})+\frac{15}{4}+\mu_0+8\mu_1&\;\;\;\;\;\;\;&\mbox{(by $\eqref{lem-22-20-1}$ of the lemma)}
	\end{aligned}
	\eeqq
	Similarly, we may prove
	\beqq
	k_D(u,{\eta}_{i_1i_2,3})\geq k_D(v,{\eta}_{i_1i_2,3})-\frac{15}{4}-\mu_0-8\mu_1.
	\eeqq
	These ensure that
	\beqq
	\begin{aligned}
		k_D(u,v)-k_D&(u,{\eta}_{i_1i_2,3})-\frac{15}{4}-\mu_0-8\mu_1\leq k_D(u,v)-k_D(v,{\eta}_{i_1i_2,3})\\
		&\leq k_D(u,{\eta}_{i_1i_2,3})\leq \ell_k(\varsigma_3)-\ell_k(\varsigma_3[v, {\eta}_{i_1i_2,3}])\\
		&\leq k_D(u,v)-k_D(v,{\eta}_{i_1i_2,3})+\frac{1}{4}\\
		&\leq k_D(u,v)-k_D(u,{\eta}_{i_1i_2,3})+4+\mu_0+8\mu_1.
	\end{aligned}
	\eeqq
	Consequently, we get
	$$\frac{1}{2}k_D(u,v)-\frac{15}{8}-\frac{1}{2}\mu_0-4\mu_1\leq k_D(u,{\eta}_{i_1i_2,3})\leq \frac{1}{2}k_D(u,v)+2+\frac{1}{2}\mu_0+4\mu_1.$$
	This proves \eqref{lem-22-09-20-3}.
	\medskip
	
	$(4)$. The proof of \eqref{lem-22-09-20-4} is similar to the one given above and thus it is omitted here.
	
	\epf
	
	\subsection{Six-tuples with Property \ref{Property-A}}\label{sub-5.4}
	In this section, we use $\mathfrak{c}$ to denote a constant with $\mathfrak{c}\geq 64\mu_1$. Given $\beta_0$, $\xi$-sequence as in previous sections, we first introduce the following useful definition.
	
\begin{figure}[htbp]
		\begin{center}
			\begin{tikzpicture}[scale=0.9]
				\draw (-5.5,-0.5) to [out=65,in=155] coordinate[pos=0] (xi1) coordinate[pos=0.3] (w_2) coordinate[pos=0.95] (w_4) (0,0) to [out=-25,in=210] coordinate[pos=0.5] (varpi) node [above, pos=0.75]{$\beta_0$} coordinate[pos=1] (xi_m) (5,0.44) ;
				\filldraw  (xi1)node[below,xshift=0cm] { $\xi_1$} circle (0.04);
				\filldraw  (w_2)node[below,xshift=0.2cm] { $w_2$} circle (0.04);
				\filldraw  (w_4)node[below,xshift=0cm] { $w_4$} circle (0.04);
				\filldraw  (varpi)node[below,xshift=0cm] { $\varpi$} circle (0.04);
				\filldraw  (xi_m)node[right,xshift=0cm] { $\xi_m$} circle (0.04);
				
				\draw (w_2) to [out=85,in=190]  coordinate[pos=0.4] (w_3) coordinate[pos=0.7] (w_1) node[below,pos=0.88,yshift=-0.1cm] { $\gamma_{12}\in \Lambda_{w_1w_2}^{\lambda}$} (2,7);
				\filldraw  (w_3)node[above,xshift=-0.3cm] { $w_3$} circle (0.04);
				\filldraw  (2,7)node[right,xshift=0cm] { $w_1\in D$} circle (0.04);
				
				
			\end{tikzpicture}
		\end{center}
		\caption{Illustration for the definition of a six-tuple $\Omega=[w_1,w_2,\gamma_{12},\varpi,w_3,w_4]$} \label{fig-25-5-15}
	\end{figure}
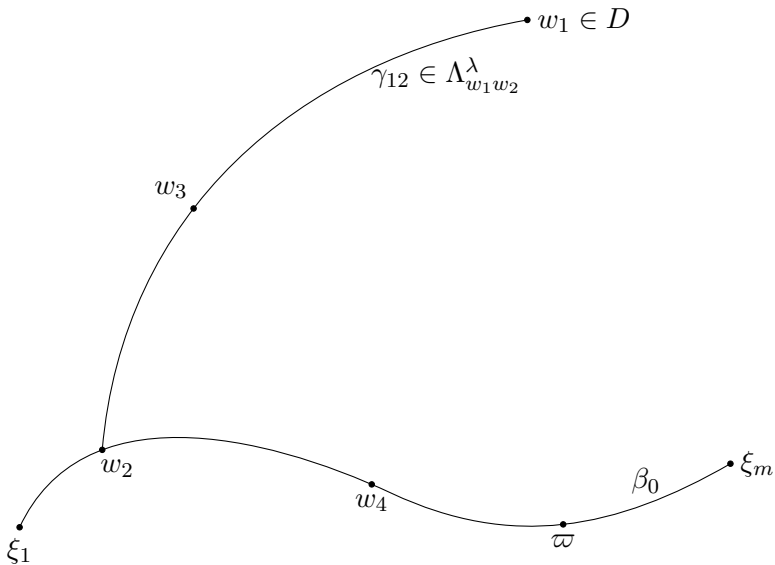

	\begin{defn}[Six-tuple]\label{def:six-tuple}
		We call $[w_1,w_2,\gamma_{12},\varpi,w_3,w_4]$ a {\it six-tuple} (see Figure \ref{fig-25-5-15}) if the following conditions are satisfied:
		\begin{itemize}
			\item$w_1\in D$ and $w_2\in \beta_0$;
			\item $\varpi\in \beta_0[w_2,\xi_m]$ and $\gamma_{12}\in \Lambda_{w_1w_2}^\lambda$;
			\item $w_3\in \gamma_{12}$ and $w_4\in \beta_0[w_2,\varpi]$.
		\end{itemize}
	\end{defn}
	
	We usually use the symbol $\Omega=[w_1,w_2,\gamma_{12},\varpi,w_3,w_4]$ to represent a six-tuple.
	
	\begin{Prop}\label{Property-A}
		Given $\mathfrak{c}\geq 64\mu_1$, a six-tuple $\Omega=[w_1,w_2,\gamma_{12},\varpi,w_3,w_4]$ is said to satisfy {\it Property \ref{Property-A}} with constant $\mathfrak{c}$ if
		\begin{enumerate}
			\item\label{3-3.1}
			$k_D(w_1, w_2)\leq \mathfrak{c}$ and $k_D(w_3, w_4)=\mathfrak{c}$;
			\item\label{3-3.2}
			for any $u\in \gamma_{12}$ and any $v\in \beta_0[w_4,\varpi]$,
			$k_D(u,v)\geq \mathfrak{c}-\frac{1}{4}.$
		\end{enumerate}
		If a six-tuple $\Omega$ satisfies Property \ref{Property-A} with constant $\mathfrak{c}$, then we shall write $\Omega^{\mathfrak{c}}=\Omega$ to indicate its dependence on $\mathfrak{c}$.
	\end{Prop}
	In most situations, we shall simply say that $\Omega^{\mathfrak{c}}=[w_1,w_2,\gamma_{12},\varpi,w_3,w_4]$ is a six-tuple with Property \ref{Property-A} to indicate that the six-tuple $[w_1,w_2,\gamma_{12},\varpi,w_3,w_4]$ satisfies Property \ref{Property-A} with constant $\mathfrak{c}$.

	The following result provides a sufficient condition for the existence of a six-tuples with Property \ref{Property-A}.
	\begin{lem}\label{lem-3.8}
		Let $w_1\in D$, $w_2\in\beta_0$ and $\gamma_{12}\in \Lambda_{w_1w_2}^\lambda$. Suppose that $0<k_D(w_1,w_2)\leq \mathfrak{c}$ and there is $\varpi\in\beta_0[w_2,\xi_m]$ such that for any $w\in \gamma_{12}$,
		$k_D(w,\varpi)>\mathfrak{c}$.
		Then there exist two points $w_3\in \gamma_{12}\backslash\{w_2\}$ and $w_4\in\beta_0[w_2,\varpi]$
		such that $\Omega^{\mathfrak{c}}=[w_1,w_2,\gamma_{12},\varpi,w_3,w_4]$ is a six-tuple with Property \ref{Property-A}.
	\end{lem}
	\bpf
	For the proof, it is helpful to revisit Figure \ref{fig-25-5-15} for a couple of times. For convenience, we write $\gamma=\gamma_{12}$ and set
	\be\label{24-2-18-1}
	\ell_k(\beta_0[w_2,\varpi])=\nu.
	\ee
	Then $\nu>\mathfrak{c}$.
	
	We start the proof with the following claim.
	\bcl\label{8-24-1} There are $u_1\in\gamma\backslash\{w_2\}$ and $v_1\in\beta_0[w_2,\varpi]$ such that
	\begin{enumerate}
		\item\label{8-24-2}
		$\ell_k(\beta_0[v_1,\varpi])\leq \nu-\mathfrak{c}+\frac{1}{8}$, and
		\item\label{8-24-3}
		one of the following statements is true:
		\begin{enumerate}
			\item\label{24-5-16-1}
			The six-tuple $[w_1,w_2,\gamma,\varpi,u_1,v_1]$ satisfies Property \ref{Property-A} with constant $\mathfrak{c}$.
			\item\label{24-5-16-2}
			There exist $u_2\in\gamma[u_1,w_1]$ and $q_1\in\beta_0[v_1,\varpi]$ such that
			$$k_D(q_1,u_2)<\mathfrak{c}-\frac{1}{4}.$$
	\end{enumerate}\end{enumerate}
	\ecl
	
	To prove the claim, select $u_1\in\gamma$ such that
	\be\label{23-08-22-8}\ell_k(\gamma[u_1,w_2])=\frac{1}{8}\min\big\{1,k_D(w_1,w_2)\big\}\leq \frac{\mathfrak{c}}{8}.\ee
	Then we see from the assumption $k_D(u_1,\varpi)>\mathfrak{c}$ and the choice of $u_1$ in \eqref{23-08-22-8} that there exists $v_1\in\beta_0[w_2,\varpi]$ such that
	\be\label{8-25-1}
	k_D(u_1,v_1)=\mathfrak{c},
	\ee and for any $w\in \beta_0[v_1,\varpi]\setminus \{v_1\}$,
	$$k_D(u_1,w)>\mathfrak{c}.$$
	
	Since the choice of $v_1$ ensures that $k_D(u_1, v)\geq \mathfrak{c}$, we infer from  \eqref{23-08-22-8} that for any $u\in \gamma[u_1,w_2]$ and any $v\in \beta_0[v_1,\varpi]$, it holds
	\beq\label{23-08-22-9}
	\quad\quad\quad k_D(u,v)\geq k_D(v,u_1)-k_D(u_1,u)\geq k_D(v,u_1)-\ell_k(\gamma[u_1,w_2])\geq \mathfrak{c}-\frac{1}{8}.
	\eeq
	This implies that $$\ell_k(\beta_0[v_1,w_2])\geq k_D(v_1,w_2)\geq \mathfrak{c}-\frac{1}{8},$$ and so, \eqref{24-2-18-1} leads to
	$$\ell_k(\beta_0[v_1,\varpi])=\ell_k(\beta_0[w_2,\varpi])-\ell_k(\beta_0[w_2,v_1])\leq \nu-\mathfrak{c}+\frac{1}{8}.$$
	
	If for any $u\in\gamma[u_1,w_1]$ and any $v\in\beta_0[v_1,\varpi]$, we have
	$$k_D(u,v)\geq \mathfrak{c}-\frac{1}{4},$$
	then it follows from \eqref{8-25-1}, \eqref{23-08-22-9} and the assumption  $0<k_D(w_1,w_2)\leq \mathfrak{c}$ that $[w_1,w_2,\gamma,\varpi,u_1, v_1]$ is a six-tuple with Property \ref{Property-A}.
	Otherwise, there must exist $u_2\in\gamma[u_1,w_1]$ and $q_1\in\beta_0[v_1,\varpi]$ such that
	$$k_D(q_1,u_2)<\mathfrak{c}-\frac{1}{4}.$$
	This proves the claim.\medskip
	
	If Claim \ref{8-24-1}\eqref{24-5-16-1} occurs, by setting $w_3=u_1$ and $w_4=v_1$, we know that $[w_1,w_2,\gamma,\varpi,w_3,w_4]$ is the required six-tuple with Property \ref{Property-A}, and the proof of lemma is complete.\medskip
	
	If Claim \ref{8-24-1}\eqref{24-5-16-2} occurs, then we make the following claim.
	\bcl\label{8-24-5} There is $v_2\in\beta_0[q_1,\varpi]$ such that
	\begin{enumerate}
		\item\label{8-24-6}
		$\ell_k(\beta_0[v_2,\varpi])\leq \nu-\mathfrak{c}-\frac{1}{8}$, and
		\item\label{8-24-7}
		one of the following statements is true:
		\begin{enumerate}
			\item\label{24-5-16-3}
			The six-tuple $[w_1,w_2,\gamma,\varpi,u_2,v_2]$ satisfies Property \ref{Property-A} with constant $\mathfrak{c}$.
			\item\label{24-5-16-4}
			There exist $u_3\in\gamma[u_1,w_1]$ and $q_2\in\beta_0[v_2,\varpi]$ such that
			$$k_D(q_2,u_3)<\mathfrak{c}-\frac{1}{4}.$$
	\end{enumerate}\end{enumerate}
	\ecl
	
	Since $k_D(u_2,\varpi)>\mathfrak{c}$ and $k_D(q_1,u_2)<\mathfrak{c}-\frac{1}{4}$, we see from Claim \ref{8-24-1}\eqref{24-5-16-2} that there is $v_2\in\beta_0[q_1,\varpi]$ such that
	\be\label{8-25-2}
	k_D(u_2,v_2)=\mathfrak{c},
	\ee and for any $w\in \beta_0[v_2,\varpi]\setminus \{v_2\}$,
	$$k_D(u_2,w)>\mathfrak{c}.$$
	Then it follows from Claim \ref{8-24-1}\eqref{24-5-16-2} that
	\beqq\label{23-08-21-1}
	\ell_k(\beta_0[v_1,v_2])\geq k_D(q_1,v_2)\geq k_D(v_2, u_2)-k_D(u_2, q_1)>\frac{1}{4},
	\eeqq
	and thus, we infer from Claim \ref{8-24-1}\eqref{8-24-2} that
	$$\ell_k(\beta_0[v_2,\varpi])=\ell_k(\beta_0[v_1,\varpi])-\ell_k(\beta_0[v_1,v_2])\leq \nu-\mathfrak{c}-\frac{1}{8}.$$
	
	If for any $u\in\gamma[u_1,w_1]$ and any $v\in\beta_0[v_2,\varpi]$, $$k_D(u,v)\geq \mathfrak{c}-\frac{1}{4},$$
	we see from
	\eqref{23-08-22-9}, \eqref{8-25-2} and the assumption $0<k_D(w_1,w_2)\leq \mathfrak{c}$ that $[w_1,w_2,\gamma,\varpi,u_2,v_2]$ is a six-tuple with Property \ref{Property-A}.
	Otherwise, there are $u_3\in\gamma[u_1,w_1]$ and $q_2\in\beta_0[v_2,\varpi]$ such that
	$$k_D(q_2,u_3)<\mathfrak{c}-\frac{1}{4}.$$
	This prove the claim.\medskip
	
	If Claim \ref{8-24-5}\eqref{24-5-16-3} occurs, by choosing $w_3=u_2$ and $w_4=v_2$, we know that $[w_1,w_2,\gamma,\varpi,w_3, w_4]$ is our required six-tuple satisfying Property \ref{Property-A} with constant $\mathfrak{c}$, and thus the proof is complete.\medskip
	
	If Claim \ref{8-24-5}\eqref{24-5-16-4} occurs, then we make the following claim.
	\bcl\label{8-24-9} There is $v_3\in\beta_0[q_2,\varpi]$ such that
	\begin{enumerate}
		\item\label{8-24-10}
		$\ell_k(\beta_0[v_3,\varpi])\leq \nu-\mathfrak{c}-\frac{3}{8}$, and
		\item\label{8-24-11}
		one of the following statements is true:
		\begin{enumerate}
			\item\label{24-5-16-5}
			The six-tuple $[w_1,w_2,\gamma,\varpi,u_3,v_3]$ satisfies Property \ref{Property-A} with constant $\mathfrak{c}$.
			\item\label{24-5-16-6}
			There exist $u_4\in\gamma[u_1,w_1]$ and $q_3\in\beta_0[v_3,\varpi]$ such that
			$$k_D(q_3,u_4)<\mathfrak{c}-\frac{1}{4}.$$
	\end{enumerate}\end{enumerate}
	\ecl
	
	Continuing this procedure for finitely many times, we arrive at the following conclusion.
	\begin{conc}\label{8-24-9} Let $a=4\nu-4\mathfrak{c}+\frac{1}{2}$ and $b=4\nu-4\mathfrak{c}+\frac{3}{2}$. Then there are an integer $K$ with $K\in [a, b)$ and $v_K\in\beta_0[q_{K-1},\varpi]$ such that
		\begin{enumerate}
			\item\label{8-24-10}
			$\ell_k(\beta_0[v_K,\varpi])\leq \nu-\mathfrak{c}+\frac{1}{8}-\frac{K-1}{4}\leq \frac{1}{4}$, and
			\item\label{8-24-11}
			one of the following statements is true:
			\begin{enumerate}
				\item\label{24-5-16-7}
				The six-tuple $[w_1,w_2,\gamma,\varpi,u_K,v_K]$ satisfies Property \ref{Property-A} with constant $\mathfrak{c}$.
				\item\label{24-5-16-8}
				There exist $u_{K+1}\in\gamma[u_1,w_1]$ and $q_K\in\beta_0[v_K,\varpi]$ such that
				$$k_D(q_K,u_{K+1})<\mathfrak{c}-\frac{1}{4}.$$
		\end{enumerate}\end{enumerate}
	\end{conc}
	
	We assert that the statement \eqref{24-5-16-8} in Conclusion \ref{8-24-9} cannot occur. Otherwise, since $k_D(u_{K+1},\varpi)>\mathfrak{c}$, we get that there is $v_{K+1}\in\beta_0[q_K,\varpi]$ such that
	$$k_D(u_{K+1},v_{K+1})=\mathfrak{c},$$
	which leads to
	\beqq\label{23-08-21-1}
	\ell_k(\beta_0[v_K,v_{K+1}])\geq k_D(q_K,v_{K+1})\geq k_D(v_{K+1},u_{K+1})-k_D(u_{K+1},q_K)>\frac{1}{4}.
	\eeqq
	This contradicts with Conclusion \ref{8-24-9}\eqref{8-24-10} since $v_{K+1}\in \beta_0[v_K,\varpi]$.
	
	Set $w_3=u_K$ and $w_4=v_K$. Then Conclusion \ref{8-24-9}\eqref{24-5-16-7} implies that the six-tuple $[w_1,w_2,\gamma,\varpi,w_3,w_4]$ is our required six-tuple satisfying Property \ref{Property-A}, and thus, the proof is complete.
	\epf

	The following three lemmas concern basic properties of six-tuples with Property \ref{Property-A}.
	
	\begin{lem}\label{lem23-01}
		Suppose that $[w_1,w_2,\gamma_{12},\varpi,w_3,w_4]$ is a six-tuple satisfying {Property \ref{Property-A}} with constant $\mathfrak{c}$ $($see Figure \ref{fig7-1}$)$. For each integer $r\in \{1,2,3\}$, let $\gamma_{rs}\in \Lambda_{w_rw_s}^\lambda$ with $s=r+1$, where $\gamma_{13}=\gamma_{12}[w_1,w_3]$ and $\gamma_{23}=\gamma_{12}[w_2,w_3]$. Then
		\begin{enumerate}
			\item\label{2022-08-20}
			\begin{enumerate}
				\item\label{23-07-26-4}
				there exist $x_1\in \gamma_{23}$, $x_2\in\gamma_{24}$ ($\gamma_{24}\in \Lambda_{w_2w_4}^{\lambda}(D)$) and $x_3\in\gamma_{34}$ such that for any $r\not=s\in\{1,2,3\}$, $$k_D(x_r,x_s)\leq \mu_1;$$
				\item\label{23-07-26-5}
				$k_D(w_3,x_3)\leq \frac{1}{2}+\mu_1$.
			\end{enumerate}
			
			\item\label{2022-08-30}
			\begin{enumerate} \item\label{23-07-26-6} there exist $p_1\in \gamma_{13}$, $p_2\in\gamma_{14}$ and $p_3\in\gamma_{34}$ such that for any $r\not=s\in\{1,2,3\}$, $$k_D(p_r,p_s)\leq \mu_1;$$
				\item\label{23-07-26-6} $k_D(w_3,p_3)\leq \frac{1}{2}+\mu_1$.
			\end{enumerate}
			
			\item\label{23-01-10-0} $k_D(w_2,w_4)\geq k_D(w_2,w_3)-\frac{5}{4}-5\mu_1+\mathfrak{c}$.
			
			\item\label{23-12-1-Add} for any $w\in\gamma_{34}$, it holds
			$$k_D(x_0,w)\leq k_D(x_0,\xi_1)-\min\{k_D(w_2,w_3)+k_D(w_3,w), k_D(w_4,w)\}+\frac{15}{2}+13\mu_1,$$
			and for each $w\in\gamma_{23}$,
			$$k_D(x_0,w)\leq k_D(x_0,\xi_1)-k_D(w_2,w)+8+16\mu_1.$$
		\end{enumerate}
	\end{lem}

	\begin{figure}[htbp]
		\begin{center}
			\begin{tikzpicture}[scale=4]
				\draw (-1.6,0.3) to [out=50,in=180] (-1,0.6)node[below] { \small $w_2$} to [out=0,in=120] (0,0)node[below left] { \small $w_4$} to [out=300,in=180] (1,-0.33);
				\filldraw  (-1.6,0.3)node[below] { \small $\xi_1$} circle (0.015);
				\filldraw  (-1,0.6) circle (0.015);
				\filldraw  (0,0) circle (0.015);
				\node at(0.8,-0.2) { \small $\beta_0$};
				\filldraw  (0.55,-0.32) circle (0.015);
				\node[below] at(0.55,-0.34) { \small $\varpi$};
				\filldraw  (1,-0.33)node[below] { \small $\xi_m$} circle (0.015);
				
				\filldraw[red]  (-1.3,0.545)node[below,xshift=0.2cm] { \small ${\color{red} \xi_{r_1}}$} circle (0.015);
				\filldraw[red]  (-0.7,0.565)node[below] { \small ${\color{red} \xi_{r_1\!+\!1}}$} circle (0.015);
				\filldraw[red]  (-0.2,0.269)node[below left] { \small ${\color{red} \xi_{r_2}}$} circle (0.015);
				\filldraw[red] (0.2,-0.21)node[below left] { \small ${\color{red} \xi_{r_2\!+\!1}}$} circle (0.015);

				\draw (-1,0.6) to [out=30,in=130] (-0.2,0.6) to [out=310,in=100] (0,0);
				\filldraw  (-0.4,0.733)node[above] { \small $x_2$} circle (0.015);
				\filldraw[red]  (-0.65,0.738)node[above] { \small ${\color{red} \eta_{r_1\!r_2,3}}$} circle (0.015);
				\filldraw[blue]  (-0.2,0.6)node[above right] { \small ${\color{blue} x_5}$} circle (0.015);
				
				\draw (-1,0.6) to [out=80,in=160] (0.5,1.5);
				\filldraw  (-0.75,1.125)node[left] { \small $x_1$} circle (0.015);
				\filldraw  (-0.4,1.415)node[above left] { \small $w_3$} circle (0.015);
				\filldraw  (0.075,1.556)node[above] { \small $p_1$} circle (0.015);
				\filldraw  (0.5,1.5)node[above right] { \small $w_1$} circle (0.015);
				
				\draw (-0.4,1.415) to [out=-40,in=70] (0,0);
				\filldraw  (-0.15,1.15)node[right] { \small $p_3$} circle (0.015);
				\filldraw  (0.04,0.77)node[right] { \small $x_3$} circle (0.015);
				\filldraw[red]  (-0.27,1.3)node[right] { \small ${\color{red} w}$} circle (0.015);
				\filldraw[blue]  (0.083,0.5)node[right] { \small ${\color{blue} w}$} circle (0.015);
				
				\draw (0.5,1.5) to [out=-65,in=40] (0,0);
				\filldraw  (0.597,1)node[right] { \small $p_2$} circle (0.015);
				
				\draw[red] (-1,0.6) to [out=100,in=160] (0.8,2)node[above right] { \small ${\color{red} x_0}$};
				\filldraw[red]  (0.2,2.042)node[below] { \small ${\color{red} \eta_{r_1\!r_2,1}}$} circle (0.015);
				
				\draw[red] (0,0) to [out=0,in=-60] (0.8,2);
				\filldraw[red]  (0.945,1.6)node[right] { \small ${\color{red} \eta_{r_1\!r_2,2}}$} circle (0.015);
				\filldraw[red]  (0.945,1)node[right] { \small ${\color{red} x_4}$} circle (0.015);
				\filldraw[blue]  (0.705,0.4)node[below right] { \small ${\color{blue} x_6}$} circle (0.015);
			\end{tikzpicture}
		\end{center}
		\caption{Illustration for the proof of Lemma \ref{lem23-01} I} \label{fig7-1}
	\end{figure}

	\bpf $(1)$ By Proposition \ref{11-19-1}$(P_4)$, $(D,k_D)$ is a $(\mu_0,\rho_0)$-Rips space, and by Lemma \ref{lem-2.3}, $\gamma_{23}$ is a $2\lambda$-curve. Then Lemma \ref{lem-8-12} guarantees the existence of the points $x_1\in\gamma_{23}$, $x_2\in\gamma_{24}$ and $x_3\in\gamma_{34}$ satisfying statement \eqref{23-07-26-4}.
	
	For \eqref{23-07-26-5}, we may estimate as follows:
	\[
	\begin{aligned}
		k_D(w_3,x_3)+k_D(x_3,w_4)&\leq \ell_k(\gamma_{34})\leq k_D(w_3,w_4)+\frac{1}{4}& \qquad &\text{(by Lemma \ref{lem-2.2})} \\
		&\leq k_D(w_4,x_1)+\frac{1}{2}&\quad &\text{(by Property \ref{Property-A})}\\
		&\leq k_D(w_4, x_3)+k_D(x_3,x_1)+\frac{1}{2}&\quad &\\
		&\leq k_D(w_4, x_3)+\frac{1}{2}+\mu_1& \quad &\text{(by statement \eqref{23-07-26-4})}
	\end{aligned}
	\]
	\medskip
	
	$(2)$ The proof of \eqref{2022-08-30} is similar to that of \eqref{2022-08-20} and thus is omitted.\medskip
	
	$(3)$ Statement \eqref{2022-08-20}  ensures that
	\beqq k_D(w_2,x_2)\geq  k_D(w_2,w_3)-k_D(w_3,x_3)-k_D(x_3,x_1)-k_D(x_1,x_2)
	\geq k_D(w_2,w_3)-\frac{1}{2}-3\mu_1\eeqq
	and that
	\beqq k_D(w_4,x_2)\geq k_D(w_4,w_3)-k_D(w_3,x_3)-k_D(x_3,x_2)\geq  k_D(w_3,w_4)-\frac{1}{2}-2\mu_1. \eeqq
	These two estimates, together with Lemma \ref{lem-2.2} and Property \ref{Property-A}\eqref{3-3.1}, imply that
	\begin{eqnarray*}
		k_D(w_2,w_4)&\geq& \ell_k(\gamma_{24})-\frac{1}{4}\geq
		k_D(w_2,x_2)+k_D(x_2,w_4)-\frac{1}{4}
		\\ \nonumber&\geq& k_D(w_2,w_3)-\frac{5}{4}-5\mu_1+\mathfrak{c}.\end{eqnarray*}
	This proves the statement \eqref{23-01-10-0}.

	$(4)$ To prove the statement \eqref{23-12-1-Add}, we first come to discuss the following. For each $w\in\gamma_{34}$,
	\be\label{H26-02-23-1}k_D(x_0,w)\leq k_D(x_0,\xi_1)-\min\{k_D(w_2,w_3)+k_D(w_3,w), k_D(w_4,w)\}+\frac{15}{2}+13\mu_1.\ee
	Next we shall construct several auxiliary points.
	By the statement \eqref{23-01-10-0}, we see that $k_D(w_2,w_4)>58\mu_1$. It follows from Lemma \ref{lem-3.4-0} that there exist two integers $r_1<r_2\in\{1,\cdots,m\}$ such that
	$w_2\in\beta_{r_1}$ and $w_4\in\beta_{r_2}$ (see Figure \ref{fig7-1}), where $\beta_i$ are defined in Section \ref{sub-5.3}. Let $\varsigma_1\in \Lambda_{w_2x_0}^\lambda$ and $\varsigma_2\in \Lambda_{w_4x_0}^\lambda$.
	
	Applying
	Lemma \ref{lem-22-09-20} (with $u$ and $v$ replaced by $w_2$ and $w_4$), we infer that there are $\eta_{r_1r_2,1}\in\varsigma_1$, $\eta_{r_1r_2,2}\in\varsigma_2$ and $\eta_{r_1r_2,3}\in\gamma_{24}$ (see Figure \ref{fig7-1}) such that
	for any $r\not=s\in\{1,2,3\}$,
	\be\label{Add-12-04-5}k_D(\eta_{r_1r_2,r},\eta_{r_1r_2,s})\leq \mu_1.\ee
	
	Let $w\in\gamma_{34}$. To establish the desired upper bound for $k_D(x_0,w)$, we consider two cases.
	\medskip
	
	\noindent
	{\bf Case 1}: $\eta_{r_1r_2,3}\in\gamma_{24}[w_2,x_2]$.  
	\medskip
	
	Under this assumption, we separate the arguments into two subcases: $w\in\gamma_{34}[w_3,x_3]$ and $w\in\gamma_{34}[x_3,w_4]$. In the former case $($see the red $w$ in Figure \ref{fig7-1}$)$,
	we know from the statement \eqref{23-07-26-5}  and Lemma \ref{lem-2.3} that
	\be\label{Add-12-04-1}\max\{k_D(w_3,w),k_D(w,x_3)\}\leq \ell_k(\gamma_{34}[w_3,x_3])\leq \frac{1}{2}+k_D(w_3,x_3)\leq 1+\mu_1,\ee
	and thus, we get from the statement \eqref{23-07-26-4} that
	\be\label{Add-12-04-2}k_D(w,x_2)\leq k_D(w,x_3)+k_D(x_3,x_2)\leq  1+2\mu_1.\ee
	
	Since $(D,k_D)$ is a $(\mu_0,\rho_0)$-Rips space, it follows
	from \eqref{Add-12-04-5} and Lemma \ref{lem-2.4-1} (replacing $x$, $y$ and $z$ by $w_4$, $\eta_{r_1r_2,2}$ and $\eta_{r_1r_2,3}$) that
	there exists $x_4\in\varsigma_2[w_4,\eta_{r_1r_2,2}]$ $($see Figure \ref{fig7-1}$)$ such that
	\beqq\label{24-1-25-5}
	k_D(x_2,x_4)\leq 1+\mu_0+\mu_1.
	\eeqq
	This, together with \eqref{Add-12-04-2}, implies
	\beqq
	k_D(w,x_4)\leq k_D(w,x_2)+k_D(x_2,x_4)\leq 2+\mu_0+3\mu_1.
	\eeqq
	Then $$k_D(w_4,w)\leq k_D(w_4,x_4)+k_D(x_4, w)\leq k_D(w_4, x_4)+ 2+\mu_0+3\mu_1,$$
	and thus, we obtain from the above two estimates and Lemma \ref{lem-2.2} that
	\beqq
	\begin{aligned}
		k_D(x_0,w)&\leq k_D(x_0, x_4)+k_D(x_4, w)\leq \ell_k(\varsigma_2)-\ell_k(\varsigma_2[w_4,x_4])+k_D(x_4, w)
		\\ &\leq k_D(x_0,w_4)-k_D(w_4,x_4)+k_D(x_4, w)+\frac{1}{4}\\ &\leq
		k_D(x_0,w_4)-k_D(w_4,w)+\frac{17}{4}+2\mu_0+6\mu_1\\
		&\stackrel{\text{Lemma }\ref{cl3.2-1}(iii)}{\leq} k_D(x_0,\xi_1)-k_D(w_4,w)+6+\mu_0+9\mu_1.
	\end{aligned}
	\eeqq
	This proves \eqref{H26-02-23-1}.
	
	In the latter case, that is, $w\in\gamma_{34}[w_4,x_3]$ $($see the blue $w$ in Figure \ref{fig7-1}$)$, since $(D,k_D)$ is a $(\mu_0,\rho_0)$-Rips space, it follows
	from statement \eqref{23-07-26-4}  and Lemma \ref{lem-2.4-1} (replacing $x$, $y$ and $z$ by $w_4$, $x_2$ and $x_3$) that
	there exists $x_5\in\gamma_{24}[w_4,x_2]$ such that
	\be\label{24-4-16-10}k_D(w,x_5)\leq 1+\mu_0+\mu_1.\ee
	
	Let us consider $w_4$, $\eta_{r_1r_2,2}$ and $\eta_{r_1r_2,3}$. As $\eta_{r_1r_2,3}\in\gamma_{24}[w_2,x_2]$, it follows
	from \eqref{Add-12-04-5} and Lemma \ref{lem-2.4-1} that
	there exists $x_6\in\varsigma_2[w_4,\eta_{r_1r_2,2}]$ such that
	$$k_D(x_5,x_6)\leq 1+\mu_0+\mu_1.$$ 
	This, together with the choice of $x_5$ in \eqref{24-4-16-10}, implies $$k_D(w,x_6)\leq k_D(w,x_5)+k_D(x_5,x_6)\leq 2(1+\mu_0+\mu_1).$$
	Combining the above estimate with Lemma \ref{lem-2.2} gives
	\beqq
	\begin{aligned}k_D(x_0,w)&\leq k_D(x_0, x_6)+k_D(x_6, w)\leq \ell_k(\varsigma_2)-\ell_k(\varsigma_2[w_4,x_6])+k_D(x_6,w)
		\\ &\leq k_D(x_0,w_4)-k_D(w_4,x_6)+k_D(x_6,w)+\frac{1}{4}\\ &\leq
		k_D(x_0,w_4)-k_D(w_4,w)+2k_D(x_6,w)+\frac{1}{4}\\ &\leq
		k_D(x_0,w_4)-k_D(w_4,w)+\frac{17}{4}+\mu_0+5\mu_1\\
		&\stackrel{\text{Lemma }\ref{cl3.2-1}(iii)}{\leq}k_D(x_0,\xi_1)-k_D(w_4,w)+6+8\mu_1,
	\end{aligned}
	\eeqq
	which again proves \eqref{H26-02-23-1}.
	\smallskip
	
	\begin{figure}[htbp]
		\begin{center}
			\begin{tikzpicture}[scale=4]
				\draw (-1.6,0.3) to [out=50,in=180] (-1,0.6)node[below] { \small $w_2$} to [out=0,in=120] (0,0)node[below left] { \small $w_4$} to [out=300,in=180] (1,-0.33);
				\filldraw  (-1.6,0.3)node[below] { \small $\xi_1$} circle (0.015);
				\filldraw  (-1,0.6) circle (0.015);
				\filldraw  (0,0) circle (0.015);
				\node at(0.8,-0.2) { \small $\beta_0$};
				\filldraw  (0.55,-0.32) circle (0.015);
				\node[below] at(0.55,-0.34) { \small $\varpi$};
				\filldraw  (1,-0.33)node[below] { \small $\xi_m$} circle (0.015);
				
				\filldraw[red]  (-1.3,0.545)node[below,xshift=0.2cm] { \small ${\color{red} \xi_{r_1}}$} circle (0.015);
				\filldraw[red]  (-0.7,0.565)node[below] { \small ${\color{red} \xi_{r_1\!+\!1}}$} circle (0.015);
				\filldraw[red]  (-0.2,0.269)node[below left] { \small ${\color{red} \xi_{r_2}}$} circle (0.015);
				\filldraw[red] (0.2,-0.21)node[below left] { \small ${\color{red} \xi_{r_2\!+\!1}}$} circle (0.015);
				
				\draw (-1,0.6) to [out=30,in=130] (-0.2,0.6) to [out=310,in=100] (0,0);
				\filldraw  (-0.4,0.733)node[above] {\small $x_2$} circle (0.015);
				\filldraw[red]  (-0.13,0.5)node[right] {\tiny ${\color{red} \eta_{r_1\!r_2\!,3}}$} circle (0.015);
				\filldraw[blue]  (-0.25,0.65)node[above right] {\small ${\color{blue} x_{8}}$} circle (0.015);
				\filldraw[purple]  (-0.05,0.3)node[right] {\small ${\color{purple} x_{8}}$} circle (0.015);
				
				\draw (-1,0.6) to [out=80,in=160] (0.5,1.5);
				\filldraw  (-0.75,1.125)node[left] {\small $x_1$} circle (0.015);
				\filldraw  (-0.4,1.415)node[above left] {\small $w_3$} circle (0.015);
				\filldraw  (0.5,1.5)node[above right] {\small $w_1$} circle (0.015);
				
				\draw (-0.4,1.415) to [out=-35,in=65] (0,0);
				\filldraw  (0.07,0.77)node[right] {\small $x_3$} circle (0.015);
				\filldraw[red]  (-0.128,1.157)node[right] { \small ${\color{red} w}$} circle (0.015);
				\filldraw[blue]  (0.107,0.4)node[right] { \small ${\color{blue} w}$} circle (0.015);
				
				\draw (0.5,1.5) to [out=-65,in=40] (0,0);
				
				\draw[red] (-1,0.6) to [out=100,in=160] (0.8,2)node[above right] {\small ${\color{red} x_0}$};
				\filldraw[red]  (0.2,2.042)node[below] {\small ${\color{red} \eta_{r_1\!r_2,1}}$} circle (0.015);
				\filldraw[red]  (-0.4,1.795)node[below right] {\small ${\color{red} x_7}$} circle (0.015);
				\filldraw[blue]  (-0.8,1.4)node[below right] {\small ${\color{blue} x_9}$} circle (0.015);
				
				\draw[red] (0,0) to [out=0,in=-60] (0.8,2);
				\filldraw[red]  (0.945,1.6)node[right] {\small ${\color{red} \eta_{r_1\!r_2,2}}$} circle (0.015);
				\filldraw[purple]  (0.8,0.56)node[below right] {\small ${\color{purple} x_{10}}$} circle (0.015);
				
			\end{tikzpicture}
		\end{center}
		\caption{Illustration for the proof of Lemma \ref{lem23-01} II} \label{fig7-2}
	\end{figure}

	\noindent
	{\bf Case 2}:  $\eta_{r_1r_2,3}\in\gamma_{24}[w_4,x_2]$. 
	
	\smallskip
	
	Under this assumption, we again separate the arguments into two subcases: $w\in\gamma_{34}[w_3,x_3]$ and $w\in\gamma_{34}[w_4,x_3]$.
	For the former $($see the red $w$ in Figure \ref{fig7-2}$)$, it follows
	from \eqref{Add-12-04-5} and Lemma \ref{lem-2.4-1} (replacing $x$, $y$ and $z$ by $w_2$, $\eta_{r_1r_2,1}$ and $\eta_{r_1r_2,3}$) that
	there exists $x_7\in\varsigma_1[w_2,\eta_{r_1r_2,1}]$ $($see Figure \ref{fig7-2}$)$ such that
	\beqq
	k_D(x_2,x_7)\leq 1+\mu_0+\mu_1.
	\eeqq
	This, together with \eqref{Add-12-04-2}, gives
	\beqq
	k_D(w,x_7)\leq k_D(w,x_2)+k_D(x_2,x_7)\leq 2+2\mu_1.
	\eeqq
	By Lemma \ref{lem-2.2}, statement \eqref{23-07-26-4} and the above two estimates, we have
	\beq\label{24-5-17-2}
	\begin{aligned}
		k_D(x_0,w)&\leq k_D(x_0, x_7)+k_D(x_7, w)\leq \ell_k(\varsigma_1)-\ell_k(\varsigma_1[w_2,x_7])+k_D(x_7, w)
		\\ &\leq k_D(x_0,w_2)-k_D(w_2,x_7)+k_D(x_7, w)+\frac{1}{4}\\
		&\leq k_D(x_0,w_2)-k_D(w_2,x_1)+k_D(x_1,x_2)+k_D(x_2,x_7)+k_D(x_7, w)+\frac{1}{4}\\
		&\leq  k_D(x_0,w_2)-k_D(w_2,x_1)+\frac{13}{4}+\mu_0+4\mu_1\\
		&\stackrel{\text{Lemma }\ref{cl3.2-1}(iii)}{\leq}k_D(x_0,\xi_1)-k_D(w_2,x_1)+5+7\mu_1.
	\end{aligned}
	\eeq
	
	It remains to estimate $k_D(w_2,x_1)$.  Note that statements \eqref{23-07-26-4} and \eqref{23-07-26-5} give
	\be\label{Add-12-04-3}k_D(x_1,w_3)\leq k_D(x_1,x_3)+k_D(x_3,w_3)\leq \frac{1}{2}+2\mu_1.\ee
	It follows from this and \eqref{Add-12-04-1} that
	$k_D(x_1,w_3)+k_D(w_3,w)\leq \frac{3}{2}+3\mu_1$,
	and thus, we have
	\beqq
	\begin{aligned}
		k_D(w_2,x_1) & \geq  k_D(w_2,x_1)+k_D(x_1,w_3)+k_D(w_3,w)-\frac{3}{2}-3\mu_1\\  &\geq k_D(w_2,w_3)+k_D(w_3,w)-\frac{3}{2}-3\mu_1.
	\end{aligned}
	\eeqq
	Substituting this to \eqref{24-5-17-2}, we infer that
	\beqq
	k_D(x_0,w)\leq  k_D(x_0,\xi_1)-\big(k_D(w_2,w_3)+k_D(w_3,w)\big)+\frac{13}{2}+10\mu_1.
	\eeqq
	This proves \eqref{H26-02-23-1}.

	In the latter case, that is, $w\in\gamma_{34}[w_4,x_3]$ $($see the blue $w$ in Figure \ref{fig7-2}$)$, it follows from statement \eqref{23-07-26-4}  and Lemma \ref{lem-2.4-1} (replacing $x$, $y$ and $z$ by $w_4$, $x_2$ and $x_3$) that
	there exists $x_8\in\gamma_{24}[w_4,x_2]$ such that
	\be\label{Add-12-04-7}k_D(w,x_8)\leq 1+\mu_0+\mu_1.\ee
	
	If $x_8\in\gamma_{24}[x_2,\eta_{r_1r_2,3}]$ $($see the blue $x_8$ in Figure \ref{fig7-2}$)$, then \eqref{Add-12-04-5} and Lemma \ref{lem-2.4-1} imply that there exists $x_9\in\varsigma_1[w_2,\eta_{r_1r_2,1}]$ such that
	\beqq
	k_D(x_8,x_9)\leq 1+\mu_0+\mu_1,
	\eeqq
	and by \eqref{Add-12-04-7}, it holds
	\beqq
	k_D(w,x_9)\leq k_D(w, x_8)+k_D(x_8,x_9)\leq 2(1+\mu_0+\mu_1).
	\eeqq
	Then we obtain from Lemma \ref{lem-2.2} and the above two estimates that
	\beq\label{24-5-18-1}
	\begin{aligned}
		k_D(x_0,w)&\leq k_D(x_0, x_9)+k_D(x_9,w)\leq \ell_k(\varsigma_1)-\ell_k(\varsigma_1[w_2,x_9])+k_D(w,x_9)
		\\ &\leq k_D(x_0,w_2)-k_D(w_2,x_9)+k_D(w,x_9)+\frac{1}{4}\\ &\leq
		k_D(x_0,w_2)-k_D(w_2,x_8)+k_D(x_8,x_9)+k_D(w,x_9)+\frac{1}{4}\\
		&\stackrel{\text{Lemma }\ref{cl3.2-1}(iii)}{\leq}k_D(x_0,\xi_1)-k_D(w_2,x_8)+5+2\mu_0+6\mu_1.
	\end{aligned}
	\eeq

	Next, we derive a lower bound for $k_D(w_2,x_8)$. Note that by statement \eqref{23-07-26-4} and \eqref{Add-12-04-7}, we have
	\beqq
	\begin{aligned}
		k_D(x_2,x_8)&\geq k_D(x_3,w)-k_D(x_2,x_3)-k_D(w,x_8)\geq
		k_D(x_3,w)-(1+\mu_0+2\mu_1),
	\end{aligned}
	\eeqq
	and that by \eqref{Add-12-04-3} and statement \eqref{23-07-26-5}, it holds
	$$k_D(x_1,w_3)+k_D(w_3,x_3)\leq 1+3\mu_1.$$
	Combining the above two estimates with Lemma \ref{lem-2.3} gives
	\beqq
	\begin{aligned}
		k_D(w_2,x_8)&\geq \ell_k(\gamma_{24}[w_2,x_8])-\frac{1}{2}\geq k_D(w_2,x_2)+k_D(x_2,x_8)-\frac{1}{2}
		\\ &\geq k_D(w_2,x_1)-k_D(x_1,x_2)+k_D(x_2,x_8)-\frac{1}{2}\\ &\geq
		k_D(w_2,x_1)+k_D(x_3,w)-\big(\frac{3}{2}+\mu_0+3\mu_1\big)\\
		&\geq k_D(w_2,w_3)-k_D(w_3,x_1)+k_D(w,w_3)-k_D(w_3,x_3)-\big(\frac{3}{2}+\mu_0+3\mu_1\big)\\
		&\geq k_D(w_2,w_3)+k_D(w_3,w)-\big(\frac{5}{2}+\mu_0+6\mu_1\big),
	\end{aligned}
	\eeqq

	Substituting the above estimate to \eqref{24-5-18-1} yields
	\begin{eqnarray*}
		k_D(x_0,w) \leq k_D(x_0,\xi_1)-\big(k_D(w_2,w_3)+k_D(w_3,w)\big)+\frac{15}{2}+13\mu_1,
	\end{eqnarray*}
	which proves \eqref{H26-02-23-1}.
	
	If $x_8\in\gamma_{24}[\eta_{r_1r_2,3},w_4]$ $($see the red $x_8$ in Figure \ref{fig7-2}$)$,
	then Lemma \ref{lem-2.4-1} and \eqref{Add-12-04-5} show that there exists $x_{10}\in\varsigma_2[w_4,\eta_{r_1r_2,2}]$ $($see Figure \ref{fig7-2}$)$
	such that
	\beqq\label{24-1-26-10}k_D(x_8,x_{10})\leq 1+\mu_0+\mu_1.\eeqq
	
	Since (\ref{Add-12-04-7}) ensures
	\beqq\label{24-1-26-11}k_D(w,x_{10})\leq k_D(w,x_8)+k_D(x_8,x_{10})\leq 2(1+\mu_0+\mu_1),\eeqq
	we derive from Lemma \ref{lem-2.2} and the above estimate that
	\beqq
	\begin{aligned}
		k_D(x_0,w)&\leq k_D(x_0, x_{10})+k_D(x_{10},w)\leq \ell_k(\varsigma_2)-\ell_k(\varsigma_2[w_4,x_{10}])+k_D(x_{10},w)
		\\ &\leq k_D(x_0,w_4)-k_D(w_4,x_{10})+k_D(x_{10},w)+\frac{1}{4}
		\\ &\leq k_D(x_0,w_4)-k_D(w_4,w)+2k_D(x_{10},w)+\frac{1}{4}
		\\ &\leq
		k_D(x_0,w_4)-k_D(w_4,w)+\frac{17}{4}+\mu_0+5\mu_1\\
		&\stackrel{\text{Lemma }\ref{cl3.2-1}(iii)}{\leq} k_D(x_0,\xi_1)-k_D(w_4,w)+6+8\mu_1,
	\end{aligned}
	\eeqq
	which again proves \eqref{H26-02-23-1}.
	
	Since $(D,k_D)$ is a $(\mu_0,\rho_0)$-Rips space, we know from Lemmas \ref{lem-2.3} and \ref{lem-8-12} that there exist three points
	$z_{1,1}\in\gamma_{w_3x_0}$, $z_{1,2}\in\gamma_{23}[w,w_3]$ and $z_{1,3}\in\gamma_{wx_0}$ such that
	for each $i\not=j\in\{1,2,3\}$,
	\be\label{H26-02-23-3}k_D(z_{1,i},z_{1,j})\leq \mu_1.\ee
	Combining the above estimate with Lemma \ref{lem-2.2} gives
	\beqq
	\begin{aligned}k_D(x_0,w)&\leq k_D(w,z_{1,3})+k_D(z_{1,3},x_0)\leq k_D(w,z_{1,2})+k_D(z_{1,1},x_0)+k_D(z_{1,2},z_{1,3})+k_D(z_{1,3},z_{1,1})
		\\ &\leq k_D(w,z_{1,2})+k_D(z_{1,1},x_0)+2\mu_1
		\\ &\leq\ell_k(\gamma_{w_3x_0})-k_D(w_3,z_{1,1})+k_D(w,z_{1,2})+2\mu_1
		\\ &\leq k_D(w_3,x_0)-k_D(w_3,z_{1,2})+k_D(z_{1,1},z_{1,2})+k_D(w,z_{1,2})+\frac{1}{4}+2\mu_1\\ &\leq
		k_D(w_3,x_0)-k_D(w_3,z_{1,2})+k_D(w,z_{1,2})+\frac{1}{4}+3\mu_1,
	\end{aligned}
	\eeqq
	and again by Lemma \ref{lem-2.2} and \eqref{H26-02-23-1},
	we get \beqq
	\begin{aligned}k_D(x_0,w)&\leq k_D(\xi_1,x_0)-k_D(w_2,w_3)-k_D(w_3,z_{1,2})+k_D(w,z_{1,2})+\frac{31}{4}+16\mu_1
		\\ &\leq k_D(\xi_1,x_0)-\ell_k(\gamma_{23})+k_D(w,z_{1,2})-k_D(w_3,z_{1,2})+8+16\mu_1
		\\ &\leq k_D(\xi_1,x_0)-k_D(w_2,w)+8+16\mu_1.
	\end{aligned}
	\eeqq
	The proof of lemma is thus complete.
	\epf

	\begin{lem}\label{lem-22-12-1}
		Let $w_1\in D$, $w_2\in D$ and $\gamma_{12}=\gamma_1[w_1,w_2]$, $w_3\in \gamma_{12}$, $\varpi\in\beta_0[w_2,\xi_m]$ and $w_4\in \beta_0[w_2,\varpi]$. Suppose that the six-tuple $[w_1,w_2,\gamma_{12},\varpi,w_3,w_4]$ satisfies Property \ref{Property-A} with constant $\mathfrak{c}$. If $k_D(w,x_0)\leq k_D(\xi_1,x_0)-k_D(w_2,w)+8+16\mu_1$, then
		$k_D(w_1,w_3)\leq 22\mu_1$.
	\end{lem}
	\bpf
	Fix two $\lambda$-curves $\varsigma_{3}\in \Lambda_{w_3w_4}^\lambda$, $\varsigma_{4}\in \Lambda_{w_4x_0}^\lambda$ and $\varsigma_{34}\in \Lambda_{w_3x_0}^\lambda$.
	Since $(D,k_D)$ is a $(\mu_0,\rho_0)$-Rips space by Proposition \ref{11-19-1}$(P_4)$, we know from Lemma \ref{lem-8-12} that there exist $u_1\in\varsigma_{34}$, $u_2\in\varsigma_3$ and
	$u_3\in\varsigma_4$ such that for any $r\not=s\in\{1,2,3\}$,
	\be\label{23-01-09-1} k_D(u_r,u_s)\leq \mu_1.\ee

	We prove the lemma by contradiction. Suppose that
	\be\label{23-01-10-2}
	k_D(w_1,w_3)> 22\mu_1.
	\ee
	
	To reach a contradiction, we will claim \be\label{23-01-10-1}k_D(u_2,w_3)\geq 6\mu_1.\ee
	Otherwise, $k_D(u_2,w_3)<6\mu_1$. Then
	$$k_D(x_0,u_2)\geq k_D(x_0,w_2)-k_D(w_2,w_3)-k_D(w_3,u_2)\geq  k_D(x_0,w_2)-k_D(w_2,w_3)-6\mu_1,$$
	and since Lemma \ref{lem-2.3} and \eqref{23-01-10-2} imply
	$$k_D(w_2,w_3)\leq \ell_k(\gamma_{12})-\ell_k(\gamma_{12}[w_1,w_3])\leq k_D(w_1,w_2)-k_D(w_1,w_3)+\frac{1}{4}< k_D(w_1,w_2)+\frac{1}{4}-21\mu_1,$$
	we see from Lemma \ref{cl3.2-1}(iii) that
	\beq\label{H26-02-24-1}
	\begin{aligned}
		k_D(x_0,u_2) &>k_D(x_0,w_2)-k_D(w_1,w_2)-\frac{1}{4}+16\mu_1\\ &\geq k_D(x_0,\xi_1)-k_D(w_1,w_2)-2+14\mu_1-2\mu_0.
	\end{aligned}
	\eeq
	
	On the other hand, since
	$$k_D(w_4,u_2)\geq k_D(w_3,w_4)-k_D(w_3,u_2)\geq k_D(w_3,w_4)-6\mu_1,$$
	we have by Lemma \ref{lem-2.3}, (\ref{23-01-09-1})  and (\ref{H26-02-24-1}) that
	\beqq
	\begin{aligned}k_D(x_0,w_4)&\geq \ell_k(\varsigma_4)-\frac{1}{4}\geq k_D(w_4,u_3)+k_D(x_0,u_3)-\frac{1}{4}
		\\ &\geq k_D(w_4,u_2)-k_D(u_2,u_3)+k_D(x_0,u_2)-k_D(u_2,u_3)-\frac{1}{4}
		\\ &\geq k_D(w_3,w_4)+k_D(x_0,u_2)-8\mu_1-\frac{1}{4}>k_D(x_0,\xi_1)-\frac{9}{4}+9\mu_1-2\mu_0.
	\end{aligned}
	\eeqq
	Further, by Lemma \ref{cl3.2-1}(iii),
	$$k_D(x_0,w_4)\leq k_D(x_0,\xi_1)+\frac{7}{4}+2(\mu_0+\mu_1).$$
	This contradicts with the above inequalities, and so (\ref{23-01-10-1}) holds.
	
	Let $u_4\in \varsigma_3[w_3,u_2]$ be such that
	\be\label{23-08-01-3}
	k_D(u_4,w_3)=3\mu_1.
	\ee
	Then we see from \eqref{23-01-09-1} and Lemma \ref{lem-2.4-1} (with $x$, $y$ and $z$ replaced by $w_3$, $u_1$ and $u_2$) that there exists $u_5\in \varsigma_{34}[w_3,u_1]$ such that
	\be\label{23-08-01-4}
	k_D(u_4,u_5)\leq 1+\mu_0+\mu_1,
	\ee
	which, together with (\ref{23-08-01-3}), ensures that
	\beqq\label{H26-02-25-1}
	\begin{aligned}
	2\mu_1-1-\mu_0&\leq k_D(w_3,u_4)-k_D(u_4,u_5)\leq k_D(w_3,u_5)\\ &\leq k_D(w_3,u_4)+k_D(u_4,u_5)\leq 1+\mu_0+4\mu_1.
	\end{aligned}
	\eeqq
	
	Since $(D,k_D)$ is a $(\mu_0,\rho_0)$-Rips space by Proposition \ref{11-19-1}$(P_4)$, there exists some $u_{1,5}\in\gamma_{12}[w_3,w_1]\cup\gamma_{w_1x_0}$ satisfying
	\be\label{8-12-1}
	k_D(u_5,u_{1,5})\leq \mu_0.
	\ee
	and thus, we know from (\ref{23-08-01-3}) and (\ref{23-08-01-4}) that
	\be\label{23-08-02-1}
	k_D(w_3,u_{1,5})\leq k_D(w_3,u_4)+k_D(u_4,u_5)+k_D(u_5,u_{1,5})\leq 1+2(\mu_0+2\mu_1).
	\ee
	Now we will claim $u_{1,5}\in\gamma_{12}[w_1,w_3]$. Otherwise, it follows from Lemma \ref{lem-2.3} and \eqref{23-08-02-1} that
	\beqq 
	\begin{aligned}
	k_D(x_0,w_1)\nonumber&\geq \ell_k(\varsigma_{34})-\frac{1}{4}\geq k_D(u_{1,5},x_0)-\frac{1}{4}
		\\ \nonumber&\geq k_D(w_3,x_0)-k_D(w_3,u_{1,5})-\frac{1}{4}\geq k_D(w_3,x_0)-\frac{5}{4}-2(\mu_0+2\mu_1)
		\\ \nonumber&\geq k_D(w_2,x_0)-k_D(w_3,w_2)-\frac{5}{4}-2(\mu_0+2\mu_1)
		\\ \nonumber&\geq k_D(w_2,x_0)-(\ell_k(\gamma_{12})-\ell_k(\gamma_{12})[w_1,w_3])-\frac{5}{4}-2(\mu_0+2\mu_1)
		\\ \nonumber&\geq k_D(w_2,x_0)-k_D(w_1,w_2)+k_D(w_1,w_3)-\frac{5}{4}-2(\mu_0+2\mu_1),
		\end{aligned}
			\eeqq 
	and so Lemma \ref{cl3.2-1}(iii) and (\ref{23-01-10-2}) yields
	\beqq
	\begin{aligned}
	k_D(x_0,w_1)\nonumber&\geq k_D(w_2,x_0)-k_D(w_1,w_2)-\frac{3}{2}-2(\mu_0-10\mu_1)
		\\ \nonumber&\geq k_D(w_2,\xi_1)-k_D(w_1,w_2)-\frac{13}{4}-2(2\mu_0-9\mu_1).
		\end{aligned}
			\eeqq 
	But we infer from the assumption in this lemma that
	$$k_D(x_0,w_1)\leq k_D(w_2,\xi_1)-k_D(w_1,w_2)+6+16\mu_1.$$ This contradicts with the above inequalities, and hence, $u_{1,5}\in\gamma_{12}[w_1,w_3]$.
	
	Then we get from Lemma \ref{lem-2.3}, \eqref{23-08-01-4} and \eqref{8-12-1} that
	\beq\label{26-02-25-3}
	\begin{aligned}
		k_D(w_4,u_{1,5})&\leq k_D(w_4,u_4)+k_D(u_4,u_5)+ k_D(u_5,u_{1,5})\leq k_D(w_4,u_4)+1+2\mu_0+\mu_1
		\\ &\leq \ell_k(\varsigma_{3})-\ell_k(\varsigma_{3}[w_3,u_4])+1+2\mu_0+\mu_1
		\\ &\leq k_D(w_3,w_4)-k_D(w_3,u_4)+\frac{5}{4}+2\mu_0+\mu_1,
	\end{aligned}
	\eeq
	which, together with \eqref{23-08-01-3}, shows that
	$$k_D(w_4,u_{1,5})\leq k_D(w_3,w_4)+\frac{5}{4}+2\mu_0-2\mu_1<k_D(w_3,w_4)-\frac{1}{4},$$
	and since the six-tuple $[w_1,w_2,\gamma_{12},\varpi,w_3,w_4]$ satisfies Property \ref{Property-A} with constant $\mathfrak{c}$, we get a contradiction, and hence the proof is complete.
	\epf

	\begin{figure}[htbp]
		\begin{center}
			\begin{tikzpicture}[scale=1]
				\draw (-5.5,-0.5) to [out=65,in=185] coordinate[pos=0] (xi1) coordinate[pos=0.3] (w12) coordinate[pos=0.75] (w14) (0,0) to [out=5,in=140] coordinate[pos=0.35] (w24) coordinate[pos=0.8] (varpi) coordinate[pos=1] (xim) (5,0.6);
				\filldraw  (xi1)node[left,xshift=0cm] {\small $\xi_1$} circle (0.04);
				\filldraw  (w12)node[below,xshift=0.1cm] {\small $w_{1,2}$} circle (0.04);
				\filldraw  (w14)node[below,rotate=-5,xshift=0.1cm] {\small $w_{1,4}=w_{2,2}$} circle (0.04);
				\filldraw  (w24)node[below,xshift=0.3cm] {\small $w_{2,4}$} circle (0.04);
				\filldraw  (varpi)node[below,xshift=0cm] {\small $\varpi$} circle (0.04);
				\filldraw  (xim)node[right,yshift=-0.1cm] {\small $\xi_m$} circle (0.04);
				
				\draw (w12) to [out=85,in=220] coordinate[pos=0.6] (w13) coordinate[pos=1] (w11) (-1.2,7);
				\filldraw  (w13)node[left,xshift=0cm] {\small $w_{1,3}=w_{2,1}$} circle (0.04);
				\filldraw  (w11)node[right,xshift=0cm] {\small $w_{1,1}$} circle (0.04);
				
				\draw (w14) to [out=90,in=-50] coordinate[pos=0.6] (w23) coordinate[pos=0.7] (z12) coordinate[pos=0.8] (u12) (w13);
				\filldraw  (w23)node[left,xshift=0cm] {\small $w_{2,3}$} circle (0.04);
				\filldraw  (z12)node[right,xshift=0cm] {\small $z_{1,2}$} circle (0.04);
				\filldraw  (u12)node[right,xshift=0cm] {\small $u_{1,2}$} circle (0.04);
				
				\draw (w24) to [out=100,in=5] coordinate[pos=0.75] (z13) coordinate[pos=0.6] (u2) (w23);
				\filldraw  (u2)node[below,xshift=-0.1cm] {\small $u_2$} circle (0.04);
				\filldraw  (z13)node[below,xshift=-0.1cm] {\small $z_{1,3}$} circle (0.04);
				
				\draw (w24) to [out=100,in=-120] coordinate[pos=0.5] (u3) coordinate[pos=1] (x0) (2.8,6.5);
				\filldraw  (u3)node[right,xshift=0cm] {\small $u_3$} circle (0.04);
				\filldraw  (x0)node[right,xshift=0cm] {\small $x_0$} circle (0.04);
				
				\draw (w23) to [out=10,in=-150] coordinate[pos=0.3] (z11) coordinate[pos=0.4] (u1) coordinate[pos=0.5] (u11) (x0);
				\filldraw  (z11)node[above,xshift=-0.2cm] {\small $z_{1,1}$} circle (0.04);
				\filldraw  (u1)node[below,xshift=0.2cm] {\small $u_1$} circle (0.04);
				\filldraw  (u11)node[above,xshift=-0.2cm] {\small $u_{1,1}$} circle (0.04);
				
				\draw (w13) to [out=0,in=170] coordinate[pos=0.3] (u13) (x0);
				\filldraw  (u13)node[below,xshift=0.2cm] {\small $u_{1,3}$} circle (0.04);
			\end{tikzpicture}
		\end{center}
		\caption{Illustration for the proof of Lemma \ref{lem-22-12-0} } \label{Fig-10-10}
	\end{figure}
	
	\begin{lem}\label{lem-22-12-0} If both $\Omega_1$ and $\Omega_2$ satisfy Property \ref{Property-A} with constant $\mathfrak{c}$, then
		$$k_D(w_{1,3},w_{2,3})=k_D(w_{2,1},w_{2,3})\leq \frac{1}{2}\big(\mathfrak{c}-k_D(w_{1,2},w_{1,3})\big)+28\mu_1.$$
	\end{lem}

	\bpf Suppose on the contrary that
	\be\label{24-1-28-2}
	k_D(w_{1,3},w_{2,3})=k_D(w_{2,1},w_{2,3})> \frac{1}{2}\big(\mathfrak{c}-k_D(w_{1,2},w_{1,3})\big)+28\mu_1.
	\ee
	
	Under this contrary assumption, we do some preparation. First of all, since for each $i\in\{1,2\}$, the six-tuple $\Omega_i=[w_{i,1},w_{i,2},\gamma_{12}^i,\varpi, w_{i,3},w_{i,4}]$ satisfies Property \ref{Property-A} with constant $\mathfrak{c}$, we know from Property \ref{Property-A}\eqref{3-3.1} that
	\be\label{23-07-27-1} k_D(w_{i,1},w_{i,2})\leq \mathfrak{c}\;\;\mbox{and}\;\; k_D(w_{i,3},w_{i,4})=\mathfrak{c}.\ee
	We take $\gamma_{w_{1,3}x_0}\in \Lambda_{w_{1,3}x_0}^{\lambda}$, $\gamma_{w_{2,3}x_0}\in \Lambda_{w_{2,3}x_0}^{\lambda}$, $\gamma_{w_{2,4}x_0}\in \Lambda_{w_{2,4}x_0}^{\lambda}$ and $\gamma_{w_{2,3}w_{2,4}}\in \Lambda_{w_{2,3}w_{2,4}}^{\lambda}$.
	Since $(D,k_D)$ is a $(\mu_0,\rho_0)$-Rips space by Proposition \ref{11-19-1}$(P_4)$, we know from Lemma \ref{lem-8-12} that there exist $u_1\in\gamma_{w_{2,3}x_0}$, $u_2\in\gamma_{w_{2,3}w_{2,4}}$ and
	$u_3\in\gamma_{w_{2,4}x_0}$ such that for any $r\not=s\in\{1,2,3\}$,
	\be\label{24-2-1-1} k_D(u_r,u_s)\leq \mu_1.\ee
	
	Since $\Omega_2$ satisfies Property \ref{Property-A} with constant $\mathfrak{c}$, by (\ref{24-1-28-2}), a similar discussion as in (\ref{23-01-10-1}) shows that
	\be\label{23-07-28-7}k_D(u_2,w_{2,3})\geq 6\mu_1.\ee
	
	
	Since $(D,k_D)$ is a $(\mu_0,\rho_0)$-Rips space by Proposition \ref{11-19-1}$(P_4)$, there exist points $u_{1,1}\in\gamma_{w_{2,3}x_0}$,
	$u_{1,2}\in\gamma_{12}^2[w_{2,1},w_{2,3}]$ and $u_{1,3}\in \gamma_{w_{1,3}x_0}$ such that for
	any $i\not=j\in\{1,2,3\}$,
	\be\label{24-1-27-1}
	k_D(u_{1,i},u_{1,j})\leq \mu_1.
	\ee
	Now we will claim 
	\be\label{clm-08-20-2}
	k_D(w_{2,3},u_{1,2})\geq 6\mu_1.
	\ee
	Otherwise, (\ref{24-1-28-2}) yields
	\be\label{26-02-25-7}k_D(w_{2,1},u_{1,2})\geq k_D(w_{2,1},w_{2,3})-k_D(w_{2,3},u_{1,2})> \frac{1}{2}\big(\mathfrak{c}-k_D(w_{1,2},w_{1,3})\big)+22\mu_1.\ee
	Then it follows from Lemma \ref{lem-2.3} and (\ref{24-1-27-1}) that
	\beqq 
	\begin{aligned}
	k_D(x_0,w_{2,1})\nonumber&\geq \ell_k(\gamma_{w_{1,3}x_0})-\frac{1}{4}\geq k_D(w_{2,1},u_{1,3})+k_D(x_0,u_{1,3})-\frac{1}{4}
	\\ \nonumber&\geq k_D(w_{2,1},u_{1,2})-k_D(u_{1,2},u_{1,3})+k_D(x_0,u_{1,2})-k_D(u_{1,2},u_{1,3})-\frac{1}{4}
	\\ \nonumber&\geq k_D(x_0,u_{1,2})+\frac{1}{2}\big(\mathfrak{c}-k_D(w_{1,2},w_{1,3})\big)-\frac{1}{4}+21\mu_1,
	\end{aligned}
	\eeqq
	and since Lemma \ref{lem-2.3} and (\ref{26-02-25-7}) imply that
	\beqq 
	\begin{aligned}
	k_D(x_0,u_{1,2})\nonumber&\geq k_D(x_0,w_{2,2})-k_D(w_{2,2},u_{1,2})\geq k_D(x_0,w_{2,2})-(\ell_k(\gamma_{12}^2)-\ell_k(\gamma_{12}^2[w_{2,1},u_{1,2}]))
	\\ \nonumber&\geq k_D(x_0,w_{2,2})-k_D(w_{2,1},w_{2,2})+k_D(u_{1,2},w_{2,1})-\frac{1}{4}
	\\ \nonumber&\geq k_D(x_0,w_{2,2})-\frac{1}{2}\big(\mathfrak{c}+k_D(w_{1,2},w_{1,3})\big)-\frac{1}{4},
	\end{aligned}
	\eeqq
	we have by Lemma \ref{cl3.2-1}(iii) that
	$$k_D(x_0,w_{2,1})\geq k_D(x_0,w_{2,2})-k_D(w_{1,2},w_{1,3})-\frac{1}{2}+21\mu_1\geq k_D(x_0,\xi_1)-k_D(w_{1,2},w_{1,3})-\frac{9}{4}+18\mu_1-2\mu_0.$$
	But Lemmas \ref{lem-2.3} and \ref{lem23-01}(\ref{23-12-1-Add}) yield that
	\beq k_D(x_0,w_{2,1})\nonumber&\leq& k_D(x_0,\xi_1)-\min\{k_D(w_{1,2},w_{1,3}),k_D(w_{2,2},w_{1,3})\}
	\\ \nonumber&\leq& k_D(x_0,\xi_1)-k_D(w_{1,2},w_{1,3})+\frac{31}{4}+13\mu_1.\eeq
	This contradicts with the above inequalities, and so (\ref{clm-08-20-2}) holds.
	
	Based on (\ref{24-2-1-1}) $\sim$ (\ref{clm-08-20-2}), we get
	\be\label{26-02-25-8}k_D(w_{2,3},u_1)\geq k_D(w_{2,3},u_2)-k_D(u_2,u_1)\geq 5\mu_1\ee
	and \be\label{26-02-25-9}k_D(w_{2,3},u_{1,1})\geq k_D(w_{2,3},u_{1,2})-k_D(u_{1,1},u_{1,2})\geq 5\mu_1.\ee
	We take $z_{1,1}\in\gamma_{w_{2,3}x_0}$ with $k_D(w_{2,3},z_{1,1})=5\mu_1$. Thus, Lemma \ref{lem-2.2}, (\ref{26-02-25-8}) and (\ref{26-02-25-9})
	imply $z_{1,1}\in \gamma_{w_{2,3}x_0}[w_{2,3},u_1]\cap \gamma_{w_{2,3}x_0}[w_{2,3},u_{1,1}]$.
	Hence we see from \eqref{24-2-1-1}, \eqref{24-1-27-1} and Lemma \ref{lem-2.4-1} (with $x$, $y$ and $z$ replaced by $w_{2,3}$, $u_1$ and $u_2$, and by $w_{2,3}$, $u_{1,1}$ and $u_{1,2}$, respectively) that there exist $z_{1,2}\in \gamma_{12}[w_{2,1},w_{2,3}]$ and $z_{1,3}\in \gamma_{w_{2,3}w_{2,4}}[w_{2,3},u_2]$ such that
	$$
	\max\{k_D(z_{1,1},z_{1,2}),k_D(z_{1,1},z_{1,3})\}\leq 1+\mu_0+\mu_1.
	$$
	This ensures that
	$$k_D(w_{2,3},z_{1,3})\geq k_D(w_{2,3},z_{1,1})-k_D(z_{1,1},z_{1,3})\geq 5\mu_1-1-\mu_0$$
	and $$k_D(z_{1,3},z_{1,2})\leq k_D(z_{1,1},z_{1,3})+k_D(z_{1,1},z_{1,2})\leq 2(1+\mu_0+\mu_1)$$
	Indeed, it follows from Lemma \ref{lem-2.2} and the above two inequalities that
	\beqq 
	\begin{aligned}
		k_D(w_{2,4},z_{1,2})&\leq k_D(w_{2,4},z_{1,3})+k_D(z_{1,3},z_{1,2})\leq \ell_k(\gamma_{w_{2,3}w_{2,4}})-\ell_k(\gamma_{w_{2,3}w_{2,4}}[w_{2,3},z_{1,3}])
	\\ &\leq k_D(w_{2,3},w_{2,4})-k_D(w_{2,3},z_{1,3})+\frac{9}{4}+2(\mu_0+\mu_1)\\ &\leq k_D(w_{2,3},w_{2,4})+\frac{13}{4}+3\mu_0-2\mu_1,
	\end{aligned}
	\eeqq
	which, together with Property \ref{Property-A}, shows $2\mu_1\leq 3\mu_0+\frac{7}{2}$, which is impossible as $\mu_1=3\mu_0$. Therefore, the proof is complete.
	\epf

	\subsection{A sequence of six-tuples}\label{sub-5.5}
	In the rest of this paper, we use
	$\tau_1$, $\tau_2$ and $\tau_3$ to denote constants defined in Section \ref{subsec:notation}.

	let $\Omega_1=[w_{1,1},w_{1,2},\gamma_{12}^1,\varpi,w_{1,3},w_{1,4}]$ be a six-tuple, where $w_{1,1}\in D$, $w_{1,2}\in \beta_0$, $\gamma_{12}^1\in \Lambda_{w_{1,1}w_{1,2}}^{\lambda}$,
	$\varpi\in\beta_0[w_{1,2},\xi_m]$, $w_{1,3}\in \gamma_{12}^1$ and $w_{1,4}\in \beta_0[w_{1,2},\varpi]$ (see Figure \ref{fig-8-16-2}).
	
	In the rest of this section, we assume that
	\be\label{24-1-29-1}
	k_D(w_{1,1},w_{1,2})=\tau_3\;\;\mbox{and}\;\;k_D(w_{1,1},w_{1,3})\leq\tau_1.
	\ee
	Then
	\be\label{24-2-13-4} 
	k_D(w_{1,2},w_{1,3})\geq k_D(w_{1,2},w_{1,1})-k_D(w_{1,1},w_{1,3})\geq -\tau_1+\tau_3\geq \frac{15}{4}\tau_2.
	\ee
	This implies that there is $w_{2,1}\in\gamma_{12}^{1}[w_{1,2},w_{1,3}]$ (see Figure \ref{fig-8-16-2}) such that
	\be\label{23-07-29-10}
	k_D(w_{1,3},w_{2,1})=\tau_2.
	\ee
	
	Based on $\Omega_1$ and $w_{2,1}$, we obtain the second six-tuple $\Omega_2=[w_{2,1},w_{2,2},\gamma_{12}^2,\varpi_2,w_{2,3},w_{2,4}]$ by setting $w_{2,2}=w_{1,2}$, $\gamma_{12}^2=\gamma_{12}^1[w_{2,1},w_{2,2}]$, $\varpi_2=w_{1,4}\in \beta_0[w_{2,2},\xi_m]$, $w_{2,3}\in \gamma_{12}^2$ and $w_{2,4}\in \beta_0[w_{2,2},\varpi_2]$ (see Figure \ref{fig-8-16-2}).
	
	
	\begin{figure}[htbp]
		\begin{center}
			\begin{tikzpicture}[scale=1]
				\draw (-6.6,-1) to [out=65,in=155] coordinate[pos=0] (xi0) (0,0) to [out=-25,in=210] coordinate[pos=0.97] (varpi) (7,1.44) to [out=30,in=210]  node[above,pos=0,xshift=0cm,yshift=0cm] {$\beta_0$} coordinate[pos=1] (xim) (7.5,1.74);
				\filldraw  (xi0)node[below] {\tiny $\xi_1$} circle (0.04);
				\filldraw  (-6.02,-0.15)node[below] {\tiny $\;\;\;\;\;w_{2,2}$} circle (0.04);
				\node[left] at(-6.02,0.05) {\tiny $w_{1,2}$};
				\node[rotate=-55] at(-6.2,-0.17) {\tiny $=$};

				\filldraw  (-3,0.79)node[below] {\tiny $w_{3\!,2}\;\;\;\;\;$} circle (0.04);
				\node[above] at(-3,0.74) {\tiny $w_{2\!,\!4}\;\;\;\;\;\;\,$};
				\node[rotate=-85] at(-3.25,0.78) {\tiny $=$};
				
				\filldraw  (-2.4,0.73)node[below] {\tiny $w_{4\!,2}\;\;\;\;$} circle (0.04);
				\node[above] at(-2.4,0.68) {\tiny $w_{3\!,\!4}\;\;\;\;\;\;$};
				\node[rotate=-85] at(-2.65,0.72) {\tiny $=$};
				
				\filldraw  (-1.8,0.616)node[below] {\tiny $w_{5\!,2}\;\;\;$} circle (0.04);
				\node[above] at(-1.8,0.57) {\tiny $w_{4\!,\!4}\;\;\;\;\;\,$};
				\node[rotate=-85] at(-2.02,0.616) {\tiny $=$};
				
				\filldraw  (-1.2,0.472)node[below] {\tiny $w_{6\!,2}\;\;$} circle (0.04);
				\node[above] at(-1.2,0.422) {\tiny $w_{5\!,\!4}\;\;\;\;\;$};
				\node[rotate=-85] at(-1.41,0.472) {\tiny $=$};
				
				\filldraw  (-0.6,0.25)node[below] {\tiny $w_{m_{\!1}\!-\!2\!,2}\;\;\;\;\;\;\;\;\;$} circle (0.04);
				\node[rotate=70] at(-0.75,0.27) {\tiny $=$};
				\node[rotate=20] at(-0.44,0.63) {\tiny $w_{\!m_{1}\!-3\!,4}$};
				
				\filldraw  (0,0)node[below] {\tiny $w_{m_{\!1}\!-\!1\!,2}\;\;\;\;\;\;\;\;\;$} circle (0.04);
				\node[rotate=70] at(-0.15,0) {\tiny $=$};
				\node[rotate=20] at(0.16,0.35) {\tiny $w_{\!m_{1}\!-2\!,4}$};
				
				\filldraw  (0.6,-0.24)node[below] {\tiny $w_{m_{\!1}\!,2}\;\;\;\;\;\;$} circle (0.04);
				\node[rotate=70] at(0.45,-0.23) {\tiny $=$};
				\node[rotate=20] at(0.77,0.1) {\tiny $w_{\!m_{1}\!-1\!,4}$};
				
				\filldraw  (1.8,-0.465) circle (0.04);
				\node[rotate=12,xshift=-0.25cm,yshift=0.08cm] at(1.59,-0.27) {\tiny $\;\;w_{m_{\!1}\:\!,4}$};

				\filldraw  (6,0.87)node[below] {\tiny $\;\;\;\;\;\;\;w_{1,4}$} circle (0.04);
				
				\filldraw  (varpi)node[below,xshift=0.1cm] {\tiny $\varpi$} circle (0.04);
				\filldraw  (xim)node[below,xshift=0.2cm] {\tiny $\xi_m$} circle (0.04);

				\filldraw  (1,10)node[right] {\tiny $w_{1,1}$} circle (0.04);
				\draw (-6.02,-0.15)  to [out=110,in=185] (1,10);
				\filldraw  (0,9.84)node[above] {\tiny $w_{1,3}$} circle (0.04);
				\filldraw  (-1,9.53)node[above left] {\tiny $y_1$} circle (0.04);
				\filldraw  (-2,9.07)node[above left] {\tiny $w_{2,1}$} circle (0.04);
				
				\filldraw  (-4,7.58) circle (0.04);
				\node[below, rotate=-4] at(-4.3,7.65) {\tiny $w_{2,3}=$};
				\node[below] at(-3.55,7.6) {\tiny $w_{3,1}$};
				
				\draw (-4,7.58)  to [out=-60,in=90] (-3,0.79);
				\filldraw  (-3.64,6.9)node[below] {\tiny $w_{3,3}\;\;\;\;\;\;$} circle (0.04);
				\node[rotate=60] at(-3.78,6.89) {\tiny $=$};
				\node[above] at(-3.64,6.8) {\tiny $\;\;\;\;\;w_{4,1}$};
				
				\draw (-3.64,6.9) to [out=-5,in=105] (-2.75,6.5)  to [out=-75,in=85] (-2.4,0.73);
				\filldraw  (-2.78,6.6)node[below] {\tiny $w_{4,3}\;\;\;\;\;\;$} circle (0.04);
				\node[rotate=60] at(-2.9,6.58) {\tiny $=$};
				\node[above] at(-2.78,6.5) {\tiny $\;\;\;\;\;w_{5,1}$};
				
				\draw (-2.78,6.6) to [out=-5,in=100] (-1.9,6)  to [out=-80,in=83] (-1.8,0.616);
				\filldraw  (-1.94,6.15)node[below] {\tiny $w_{5,3}\;\;\;\;\;\;\,$} circle (0.04);
				\node[rotate=60] at(-2.05,6.13) {\tiny $=$};
				\node[above] at(-1.94,6.05) {\tiny $\;\;\;\;\;w_{6,1}$};
				
				\draw (-1.94,6.15) to [out=-5,in=100] (-1.1,5.5)  to [out=-80,in=83] (-1.2,0.472);
				
				\draw[dashed] (-1.125,5.6) to [out=-5,in=105] (-0.31,5.2);
				\draw (-0.31,5.2)  to [out=-85,in=80] (-0.6,0.25);
				\filldraw  (-0.31,5)node[left] {\tiny $w_{m_{1}-2,3}$} circle (0.04);
				\node[rotate=30] at(-0.32,5.1) {\tiny $=$};
				\node[above] at(-0.31,4.9) {\tiny $\;\;\;\;\;\;\;\;\;\;\;\;\;\;w_{\!m_{\!1}\!-\!1,1}$};
				
				\draw (-0.31,5) to [out=-5,in=95] (0.42,4.5)  to [out=-92,in=80] (0,0);
				\filldraw  (0.405,4.6)node[below] {\tiny $w_{m_{1}\!-\!1,3}\;\;\;\;\;\;\;\;\;\;\;\;\;$} circle (0.04);
				\node[rotate=45] at(0.25,4.6) {\tiny $=$};
				\node[above] at(0.405,4.5) {\tiny $\;\;\;\;\;\;\;w_{\!m_{\!1}\!,1}$};
				
				\draw (0.405,4.6) to [out=-5,in=90] (1.1,4)  to [out=-90,in=77] (0.6,-0.24);
				\filldraw  (1.09,4.1) circle (0.04);
				\node[below] at(0.767,4.19) {\tiny $w_{m_{\!1}\!,\!3}$};
				
				\draw (1.09,4.1) to [out=-40,in=86] coordinate[pos=0.2] (w) coordinate[pos=0.3] (eta12) coordinate[pos=0.5] (eta6) (1.8,-0.465);
				\filldraw  (w)node[above,xshift=0.2cm] {\tiny$w$} circle (0.04);
				\filldraw  (eta12)node[below,xshift=-0.25cm] {\tiny$\eta_{1,2}$} circle (0.04);
				\filldraw  (eta6)node[left,xshift=0cm] {\tiny$\eta_{6}$} circle (0.04);
				
				\draw (-4,7.58) to [out=-5,in=130] (1.6,5.3) to [out=-50,in=80] (1.8,-0.465);
				
				\draw (-4,7.58) to [out=5,in=150] (2.7,6.2) to [out=-30,in=105] (6,0.87);
				\filldraw  (1.5,6.82)node[above] {\tiny $y_2$} circle (0.04);
				
				\draw (0,9.84) to [out=-10,in=90] coordinate[pos=0.6] (w4) coordinate[pos=0.425] (eta2) coordinate[pos=0.25] (y3) (6,0.87);
				\filldraw  (y3)node[above,xshift=0.3cm] {\tiny $y_3$} circle (0.04);
				
				\draw (1.8,-0.465) to [out=40,in=160] coordinate[pos=0.34] (w2) coordinate[pos=0.5] (eta4) coordinate[pos=0.65] (eta3)  (6,0.87);
				\filldraw  (w2)node[above] {\tiny $\omega_2$} circle (0.04);
				\filldraw  (eta4)node[above] {\tiny $\eta_4$} circle (0.04);
				\filldraw  (eta3)node[above] {\tiny $\eta_3$} circle (0.04);

				\draw (1.09,4.1) to [out=30,in=180] coordinate[pos=0.31] (eta13) coordinate[pos=1] (x0) (7,6);
				\filldraw  (eta13)node[below,xshift=0.2cm] {\tiny $\eta_{1,3}$} circle (0.04);
				\filldraw  (x0)node[right] {\tiny $x_0$} circle (0.04);
				
				\draw (6,0.87) to [out=70,in=-90] coordinate[pos=0.4] (eta2) (7,6);
				\filldraw  (eta2)node[right] {\tiny $\eta_2$} circle (0.04);
				
				\draw (1.8,-0.465) to [out=60,in=220] coordinate[pos=0.25] (eta5) coordinate[pos=0.3] (eta1) coordinate[pos=0.4] (eta11) (7,6);
				\filldraw  (eta5)node[right,yshift=-0.1cm] {\tiny $\eta_{5}$} circle (0.04);
				\filldraw  (eta1)node[right,yshift=-0.1cm] {\tiny $\eta_{1}$} circle (0.04);
				\filldraw  (eta11)node[right,yshift=-0.1cm] {\tiny $\eta_{1,1}$} circle (0.04);
				
			\end{tikzpicture}
		\end{center}
		\caption{Illustration for the proof of Lemma \ref{lem-22-09-30}} \label{fig-8-16-2}
	\end{figure}
	
	
	Then we obtain the third six-tuple $\Omega_3=[w_{3,1},w_{3,2},\gamma_{12}^3,\varpi_3,w_{3,3},w_{3,4}]$ by setting $w_{3,1}=w_{2,3}$, $w_{3,2}=w_{2,4}$, $\gamma_{12}^3\in \Lambda_{w_{2,3}w_{2,4}}^\lambda$, $\varpi_3=w_{1,4}\in \beta_0[w_{32},\xi_m]$, $w_{3,3}\in \gamma_{12}^3$ and $w_{3,4}\in \beta_0[w_{3,2},\varpi_3]$
	and the fourth six-tuple $\Omega_4=[w_{4,1},w_{4,2},\gamma_{12}^4,\varpi_4,w_{4,3},w_{4,4}]$ by setting $w_{4,1}=w_{3,3}$, $w_{4,2}=w_{3,4}$, $\gamma_{12}^4\in \Lambda_{w_{3,3}w_{3,4}}^\lambda$, $\varpi_4=w_{1,4}\in \beta_0[w_{4,2},\xi_m]$,
	$w_{4,3}\in \gamma_{12}^{4}$ and $w_{4,4}\in \beta_0[w_{4,2},\varpi_4]$ (see Figure \ref{fig-8-16-2}).
	
	By repeating this procedure, we find a finite sequence of six-tuples $\{\Omega_i\}_{i=1}^{m_1}$ with $\Omega_i=[w_{i,1},w_{i,2},\gamma_{12}^i,\varpi_i,w_{i,3},w_{i,4}]$ obtained by setting
	for each $i\in\{3,\cdots,m_1\}$, $w_{i,1}=w_{i-1,3}$, $w_{i,2}=w_{i-1,4}$, $\gamma_{12}^i\in \Lambda_{w_{i-1,3}w_{i-1,4}}^\lambda$, $\varpi_i=w_{1,4}\in \beta_0[w_{i,2},\xi_m]$, $w_{i,3}\in \gamma_{12}^i$ and $w_{i,4}\in \beta_0[w_{i,2},\varpi_i]$.
	
	Then we have the following useful estimates.
	\begin{lem}\label{lem-22-09-30}
		$(1)$
		Suppose that $m_1\geq 2$ and $\Omega_1$ satisfies Property \ref{Property-A} with constant $\tau_3$. Then
		\beqq
		k_D(w_{1,4},w_{2,3})\geq k_D(w_{1,3},w_{1,4})+k_D(w_{1,3},w_{2,3})-\frac{5}{4}-4\mu_1.\eeqq
		
		\noindent $(2)$
		Suppose $m_1=2$. If both $\Omega_1$ and $\Omega_2$ satisfy Property \ref{Property-A} with constant $\tau_3$, then
		
		$(2a)$
		$k_D(w_{2,3},w_{2,4})\leq k_D(w_{1,4},w_{2,3})+\frac{7}{4}+4\mu_1-\tau_2$, and
		
		$(2b)$
		for any $w\in \gamma_{12}^2$, it holds
		$
		k_D(w_{1,4},w)> \tau_3.$
		
		\noindent $(3)$
		Suppose $m_1\geq 3$ and for each $i\in\{1,\cdots,m_1\}$, $\Omega_i$ satisfies Property \ref{Property-A} with constant $\tau_3$. If $$k_D(w_{2,3},w_{m_1,3})\leq -56\mu_1-\frac{\tau_1+\tau_2}{2^{m_1-2}}+\tau_2,$$ then
		for any $\gamma_{12}^{m_1}\in\Lambda_{w_{m_1,3}w_{m_1,4}}^{\lambda}$ and $w\in \gamma_{12}^{m_1}$, it holds
		$
		k_D(w_{1,4},w)> \tau_3.$
	\end{lem}
	\begin{proof}
		(1). Fix $\varsigma_{23}\in \Lambda_{w_{1,4}w_{2,3}}^{\lambda}$ and
		$\gamma_{34}^1\in\Lambda_{w_{1,3}w_{1,4}}^{\lambda}$. 		
		Since $(D,k_D)$ is a $(\mu_0,\rho_0)$-Rips space by Proposition \ref{11-19-1}$(P_4)$, we obtain from Lemma \ref{lem-8-12} that there exist $y_1\in \gamma_{12}^{1}[w_{1,3},w_{2,3}]$, $y_2\in\varsigma_{23}$ and
		$y_3\in\gamma_{34}^{1}$ (see Figure \ref{fig-8-16-2}) such that for any $r\not=s\in\{1,2,3\}$,
		\be\label{23-07-29-5}
		k_D(y_r,y_s)\leq \mu_1.
		\ee
		
		Since $\Omega_1$ satisfies Property \ref{Property-A} with constant $\tau_3$, it follows that
		$$ k_D(w_{1,3},w_{1,4})=\tau_3\leq k_D(y_1,w_{1,4})+\frac{1}{4},$$
		and thus, we get from the above estimate, Lemma \ref{lem-2.2} and \eqref{23-07-29-5} that
		\beqq\label{23-07-29-6}
		\begin{aligned}
		k_D(w_{1,3},y_3)&\leq  \ell_{k}(\gamma_{34}^{1})-\ell_{k}(\gamma_{34}^{1}[w_{1,4},y_3])
		\leq
		k_D(w_{1,3},w_{1,4})-k_D(w_{1,4},y_3)+\frac{1}{4}\\ &\leq
		k_D(y_{1},w_{1,4})-k_D(w_{1,4},y_3)+\frac{1}{2}\leq k_D(y_1,y_3)+\frac{1}{2} \leq \frac{1}{2}+\mu_1.
		\end{aligned}
		\eeqq
		This implies
		\be\label{Add-12-06-1}\quad \quad k_D(w_{1,4},y_3) \geq k_D(w_{1,4},w_{1,3})-k_D(w_{1,3},y_3)
		\geq k_D(w_{1,3},w_{1,4})-\frac{1}{2}-\mu_1
		\ee
		and 
		$$
		k_D(w_{2,3},y_2)\geq k_D(w_{2,3},w_{1,3})-k_D(w_{1,3},y_3)-k_D(y_3,y_2)\geq k_D(w_{2,3},w_{1,3})-\frac{1}{2}-2\mu_1.
		$$
		Moreover, \eqref{23-07-29-5} and \eqref{Add-12-06-1} ensure that
		\beqq
		k_D(w_{1,4},y_2)\geq k_D(w_{1,4},y_3)-k_D(y_3, y_2)\geq
		k_D(w_{1,3},w_{1,4})-\frac{1}{2}-2\mu_1.
		\eeqq
		Therefore, it follows from Lemma \ref{lem-2.2} and the above two estimates that
		\beqq
		\begin{aligned}
		k_D(w_{1,4},w_{2,3})&\geq  \ell_{k}(\varsigma_{23})-\frac{1}{4}\geq
		k_D(w_{2,3},y_2)+k_D(y_2,w_{1,4})-\frac{1}{4}\\
		&\geq k_D(w_{1,3},w_{1,4}) + k_D(w_{1,3},w_{2,3})-\frac{5}{4}-4\mu_1.
		\end{aligned}
		\eeqq
		This proves (1).

		(2). The proofs are similar to those used in (3) below and thus will be omitted here.
		
		(3). Since the choice of $w_{2,1}$ in \eqref{23-07-29-10}, the assumption $k_D(w_{2,3},w_{m_1,3})\leq -36\mu_1-\tau_1+\tau_2$ and Lemma \ref{lem-2.3} guarantee that
		$$ k_D(w_{1,3},w_{2,3}) \geq \ell_{k}(\gamma_{12}^{1}[w_{1,3},w_{2,3}])-\frac{1}{2}\geq \ell_{k}(\gamma_{12}^{1}[w_{1,3},w_{2,1}])-\frac{1}{2}\geq\tau_2-\frac{1}{2},
		$$ which, together with Lemma \ref{lem-22-09-30}$(1)$, shows
		$$k_D(w_{1,4},w_{2,3})\geq k_D(w_{1,3},w_{1,4}) + k_D(w_{1,3},w_{2,3})-\frac{5}{4}-4\mu_1\geq k_D(w_{2,3},w_{2,4})+\tau_2-\frac{7}{4}-4\mu_1.$$
		Hence, Lemma \ref{lem-22-09-30}$(2a)$ holds.
		
		Since Property \ref{Property-A}(\ref{3-3.1}), we have by Lemmas \ref{lem23-01}(\ref{23-12-1-Add}) and Lemma \ref{lem-22-12-1},
		$$k_D(w_{2,1},w_{2,3})\leq 22\mu_1.$$
		Then we get that
		\beq
		k_D(w_{2,2},w_{2,3})\nonumber&\geq& \ell_k(\gamma_{12}^1)-\ell_k(\gamma_{12}^1[w_{2,3},w_{1,1}])\\ \nonumber&\stackrel{Lemma \ref{lem-2.2}}{\geq}&
		k_D(w_{1,2},w_{1,1})-(k_D(w_{1,1},w_{1,3})+k_D(w_{1,3},w_{2,1})+k_D(w_{2,1},w_{2,3}))-\frac{1}{4}\\
		\nonumber&\stackrel{(\ref{23-07-29-10})}{\geq}& \tau_3-\frac{1}{4}-22\mu_1-\tau_2-k_D(w_{1,1},w_{1,3}).
		\eeq
		Moreover, Lemma \ref{lem-22-12-0} implies that
		\be\label{26-02-26-1}k_D(w_{2,3},w_{3,3})\leq \frac{1}{2}(\tau_3-k_D(w_{2,2},w_{2,3}))+28\mu_1\leq \frac{1}{2}(\frac{1}{2}+\tau_2+k_D(w_{1,1},w_{1,3}))+39\mu_1.\ee
		For each $i\in\{2,\cdots,m_1-1\}$, since $\Omega_i$ satisfies Property \ref{Property-A} with constant $\tau_3$, we have
		$$k_D(w_{i,3},w_{i+1,3})\geq \tau_3-k_D(w_{i+1,2},w_{i+1,3}).$$
		This, together with Lemma \ref{lem-22-12-0}, shows that for $i\in\{3,\cdots,m_1-1\}$,
		$$k_D(w_{i,3},w_{i+1,3})\leq \frac{1}{2}k_D(w_{i,3},w_{i-1,3})+28\mu_1\leq \frac{1}{2^{i-2}}k_D(w_{2,3},w_{3,3})+28\mu_1\sum_{t=3}^{i}\frac{1}{2^{t-3}},$$
		and so by (\ref{26-02-26-1}),
		\be\label{26-02-26-2}k_D(w_{i,3},w_{i+1,3})<\frac{1}{2^{i-1}}(\tau_1+\tau_2+22\mu_1)+56\mu_1.\ee
		
		We take $\varsigma_1\in \Lambda_{w_{m_1,4}x_0}^{\lambda}$, $\varsigma_2\in \Lambda_{w_{1,4}x_0}^{\lambda}$ $\varsigma_3\in \Lambda_{w_{m_1,3}x_0}^{\lambda}$ and  $\varsigma_4\in \Lambda_{w_{m_1,4}w_{1,4}}^{\lambda}$.
		It follows from Lemma \ref{lem-22-09-20} that there exist $\eta_1\in\varsigma_1$, $\eta_2\in \varsigma_2$ and $\eta_3\in\varsigma_3$
		such that for each $i\not=j\in\{1,2,3\}$,\be\label{26-02-26-3}k_D(\eta_i,\eta_j)\leq\mu_1\ee
		and \be\label{26-02-26-4}\frac{1}{2}(k_D(w_{m_1,4},w_{1,4})-\mu_0)-\frac{15}{8}-4\mu_1\leq k_D(\eta_3,w_{m_1,4})\leq \frac{1}{2}(k_D(w_{m_1,4},w_{1,4})+\mu_0)+2+4\mu_1.\ee
		Since $(D,k_D)$ is a $(\mu_0,\rho_0)$-Rips space by Proposition \ref{11-19-1}$(P_4)$, we obtain from Lemma \ref{lem-8-12} that there exist $\eta_{1,1}\in \varsigma_1$, $\eta_{1,2}\in\gamma_{w_{m_1,3}w_{m_1,4}}$ and
		$\eta_{1,3}\in\varsigma_4$  such that for any $r\not=s\in\{1,2,3\}$,
		\be\label{26-02-26-5}
		k_D(\eta_{1,r},\eta_{1,s})\leq \mu_1.
		\ee
		Now we will claim  \be\label{26-02-26-6}k_D(w_{m_1,4},\eta_{1,2})\geq \tau_3-42\mu_1-\frac{1}{2^{m_1-2}}(\tau_1+\tau_2).\ee
		Otherwise, $k_D(w_{m_1,4},\eta_{1,2})< \tau_3-42\mu_1-\frac{1}{2^{m_1-2}}(\tau_1+\tau_2)$. Then we see from Lemma \ref{cl3.2-1}(iii) that
		$$k_D(x_0,\eta_{1,2})\geq k_D(x_0,w_{m_1,4})-k_D(w_{m_1,4},\eta_{1,2})\geq  k_D(x_0,\xi_1)+40\mu_1-\tau_3-\frac{7}{4}-2\mu_0+\frac{1}{2^{m_1-2}}(\tau_1+\tau_2),$$
		and further, by \eqref{26-02-26-5},
		$$k_D(x_0,\eta_{1,3})\geq k_D(x_0,\eta_{1,2})-k_D(\eta_{1,3},\eta_{1,2})\geq k_D(x_0,\xi_1)+39\mu_1-\tau_3-\frac{7}{4}-2\mu_0+\frac{1}{2^{m_1-2}}(\tau_1+\tau_2),$$
		and Lemma \ref{lem-2.2} and \eqref{26-02-26-5} yield
		\beqq 
		\begin{aligned}
		k_D(x_0,&w_{m_1,3})\geq \ell_k(\varsigma_3)-\frac{1}{4}\geq k_D(x_0,\eta_{1,3})+k_D(w_{m_1,3},\eta_{1,3})-\frac{1}{4}
		\\ &\geq k_D(x_0,\xi_1)+k_D(w_{m_1,3},\eta_{1,2})-k_D(\eta_{1,3},\eta_{1,2})+39\mu_1-2-\tau_3-2\mu_0+\frac{1}{2^{m_1-2}}(\tau_1+\tau_2)
		\\ &\geq k_D(x_0,\xi_1)+80\mu_1-2-2\mu_0-\tau_3+\frac{1}{2^{m_1-2}}(\tau_1+\tau_2).
	\end{aligned}
		\eeqq
		Again since Lemma \ref{lem23-01} and (\ref{26-02-26-2}), we have
		\beqq
		\begin{aligned}
		 k_D(x_0,w_{m_1,3})&\leq k_D(x_0,\xi_1)-k_D(w_{m_1,3},w_{m_1,2})+8+16\mu_1
		\\ &\leq k_D(x_0,\xi_1)+k_D(w_{m_1,3},w_{m_1-1,3})+8+16\mu_1-\tau_3\\ &< k_D(x_0,\xi_1)+8+\frac{1}{2^{m_1-2}}(\tau_1+\tau_2)+78\mu_1-\tau_3,
		\end{aligned}
		\eeqq
		which, contradicts with the above inequalities, and (\ref{26-02-26-2}) holds.
		
		Next, we start the proof of (3b) with the following estimate:
		\be\label{23-08-29-1}k_D(w_{1,4},w_{m_1,4})\geq -48\mu_1-\frac{\tau_1+\tau_2}{2^{m_1-2}}+2\tau_3.\ee
		Suppose (\ref{23-08-29-1}) doesn't hold. Then by (\ref{26-02-26-4}), we let $\eta_4\in \varsigma_4[w_{m_1,4},\eta_3]$ be such that
		$$k_D(w_{m_1,4},\eta_4)=\frac{1}{2}k_D(w_{m_1,4},w_{1,4})-24\mu_1-\frac{\tau_1+\tau_2}{2^{m_1-1}}.$$
		Considering $w_{m_1,4}$, $\eta_3$ and $\eta_1$, it follows
		from \eqref{26-02-26-3} and Lemma \ref{lem-2.4-1} that
		there exists $\eta_5\in\varsigma_1[w_{m_1,4},\eta_1]$ (see Figure \ref{fig-8-16-2}) such that
		\be\label{23-07-30-0} k_D(\eta_4,\eta_5)\leq 1+\mu_0+\mu_1.\ee
		Then, we have 
		$$k_D(w_{m_1,4},\eta_5)\leq k_D(w_{m_1,4},\eta_4)+k_D(\eta_4,\eta_5)\leq \frac{1}{2}k_D(w_{m_1,4},w_{1,4})-\frac{1}{2^{m_1-2}}(\tau_1+\tau_2)+1+\mu_0-23\mu_1,$$
		and based on (\ref{26-02-26-5}) and (\ref{26-02-26-6}),
		$$k_D(w_{m_1,4},\eta_{1,1})\geq k_D(w_{m_1,4},\eta_{1,2})-k_D(\eta_{1,1},\eta_{1,2})\geq \tau_3-\frac{\tau_1+\tau_2}{2^{m_1-2}}-43\mu_1,$$
		which, together with \ref{lem-2.3}, shows $\eta_5\in\in\varsigma_1[w_{m_1,4},\eta_{1,1}]$. Again by using (\ref{26-02-26-5}) and Lemma \ref{lem-2.4-1}, there
		exists $\eta_6\in\gamma_{w_{m_1,3}w_{m_1,4}}[w_{m_1,4},\eta_{1,2}]$ satisfies
		\be\label{26-02-28-1}
		k_D(\eta_5,\eta_6)\leq 1+\mu_0+\mu_1.
		\ee
		Then it follows from  Lemma \ref{lem-2.2} that
		\[
		\begin{aligned}
			k_D(w_{1,4},\eta_{1,4})&\leq \ell_k(\varsigma_4)-\ell_k(\varsigma_4[w_{m_1,4},\eta_4])\leq k_D(w_{m_1,4},w_{1,4})-k_D(w_{m_1,4},\eta_4)+\frac{1}{4}\\
			&\leq \frac{1}{2}k_D(w_{m_1,4},w_{1,4})+\frac{1}{4}+24\mu_1+\frac{\tau_1+\tau_2}{2^{m_1-1}}.
		\end{aligned}
		\]
		And based on  Lemma \ref{lem-2.2}, (\ref{23-07-30-0}) and (\ref{26-02-28-1}) that
		\beqq 
		\begin{aligned}
		k_D(w_{m_1,3},\eta_6)&\leq \ell_k(\gamma_{12}^{m_1})- \ell_k(\gamma_{12}^{m_1}[w_{m_1,4},\eta_6])\leq \tau_3-k_D(w_{m_1,4},\eta_6)+\frac{1}{4}
		\\ &\leq \tau_3-(k_D(w_{m_1,4},\eta_5)-k_D(\eta_5,\eta_6))+\frac{1}{4}
		\\ &\leq \tau_3-(k_D(w_{m_1,4},\eta_4)-k_D(\eta_5,\eta_4)-k_D(\eta_5,\eta_6))+\frac{1}{4}
		\\ &\leq \tau_3-\frac{1}{2}k_D(w_{m_1,4},w_{1,4})+\frac{9}{4}+26\mu_1+2\mu_0+\frac{\tau_1+\tau_2}{2^{m_1-1}}
		\end{aligned}
		\eeqq
		and
		further, by the assumption as in this Lemma, (\ref{23-07-30-0}) and (\ref{26-02-28-1}),
		we have
		\beqq
\begin{aligned}
	k_D(w_{2,3},w_{1,4})&\leq k_D(w_{2,3},w_{m_1,3})+k_D(w_{m_1,3},\eta_6)+k_D(w_{1,4},\eta_{1,4})+k_D(\eta_5,\eta_4)+k_D(\eta_5,\eta_6)
	\\ &\leq k_D(w_{2,3},w_{m_1,3})+\frac{5}{2}+2\mu_0+50\mu_1+\frac{\tau_1+\tau_2}{2^{m_1-2}}
	\\ &\leq \tau_3+2\mu_0+\tau_2-\frac{9}{4}-6\mu_1.
\end{aligned}
\eeqq
		This contradicts with Lemma \ref{lem-22-09-30}$(2a)$, and so (\ref{23-08-29-1}) holds.
		
		Now, we may finish the proof of $(3)$ based on \eqref{23-08-29-1}. For this, let $w\in\gamma_{12}^{m_1}$ (see Figure \ref{fig-8-16-2}). We consider two cases: $k_D(w_{m_1,3},w)\leq 50\mu_1+\frac{\tau_1+\tau_2}{2^{m_1-2}}$ and $k_D(w_{m_1,3},w)> 50\mu_1+\frac{\tau_1+\tau_2}{2^{m_1-2}}$.
		
		In the former case, it follows from the assumption $k_D(w_{2,3},w_{m_1,3})\leq -56\mu_1-\frac{\tau_1+\tau_2}{2^{m_1-2}}+\tau_2$ that
		\beq\label{24-2-3-1}
		\begin{aligned}
			k_D(w_{1,4},w) &\geq  k_D(w_{1,4},w_{2,3})-k_D(w_{2,3},w_{m_1,3})-k_D(w_{m_1,3},w)
			\\  &\geq  k_D(w_{1,4},w_{2,3})-k_D(w_{2,3},w_{m_1,3})-50\mu_1-\frac{\tau_1+\tau_2}{2^{m_1-2}}\\
			&\geq k_D(w_{1,4},w_{2,3})+6\mu_1-\tau_2.
		\end{aligned}
		\eeq
		On the other hand, statement $(1)$ and Lemma \ref{lem-2.3} lead to
		\[
		\begin{aligned}
			k_D(w_{1,4},w_{2,3})& \geq  k_D(w_{1,3},w_{1,4})+\ell_k(\gamma_{12}^{1}[w_{1,3},w_{2,3}])-\frac{7}{4}-4\mu_1 \\ \nonumber
			& \geq  k_D(w_{1,3},w_{1,4})+\ell_k(\gamma_{12}^{1}[w_{1,3},w_{2,1}])-\frac{7}{4}-4\mu_1\\
			\nonumber &\stackrel{\eqref{23-07-29-10}+{\rm Property\ } \ref{Property-A}(\ref{3-3.1}) }{\geq} -\frac{7}{4}-4\mu_1+\tau_2+\tau_3.
		\end{aligned}
		\]
		Combining this with \eqref{24-2-3-1} gives
		\beqq
		k_D(w_{1,4},w)\geq -\frac{7}{4}+2\mu_1+\tau_3>\tau_3.
		\eeqq
		
		In the latter case, that is, $k_D(w_{m_1,3},w)>50\mu_1+\frac{\tau_1+\tau_2}{2^{m_1-2}}$,  it follows from Lemma \ref{lem-2.2} that
		\beqq
		\begin{aligned}k_D(w,w_{m_1,4})&\leq \ell_k(\gamma_{34}^{m_1})-\ell_k(\gamma_{12}^{m_1}[w_{m_1,3},w])\\  & \leq 
			k_D(w_{m_1,3},w_{m_1,4})-k_D(w_{m_1,3},w)+\frac{1}{4}\\  & \leq  \frac{1}{4}-50\mu_1-\frac{\tau_1+\tau_2}{2^{m_1-2}}+\tau_3,\end{aligned}
			\eeqq 
		and thus, we infer from \eqref{23-08-29-1} that
		\beqq
		k_D(w_{1,4},w)\geq k_D(w_{1,4},w_{m_1,4})-k_D(w_{m_1,4},w)\geq -\frac{1}{4}+2\mu_1+\tau_3>\tau_3.
		\eeqq
		In either cases, we have verified (3) and hence the proof is complete.
	\end{proof}

	\subsection{Iterative construction of new six-tuples}
	In this section, we shall prove an useful lemma, which gives a way to iteratively construct new six-tuples based on old ones. To this end, for each integer $n\geq 1$, we set
	\be\label{24-4-20-1} \varrho_{n}=\tau_2\sum_{i=1}^{n}\frac{1}{4^i}+64\mu_1\sum_{i=1}^{n}\frac{1}{2^i}.\ee
	Then, it is clear that
	\be\label{24-2-28-4}
	32\mu_1+\frac{1}{4}\tau_2< \varrho_{n}< 64\mu_1+\frac{1}{3}\tau_2.
	\ee
	
	
	\begin{lem}\label{23-08-20} Let $w_{1,1}\in D$, $w_{1,2}\in\beta_0$, $\varpi\in\beta_0[w_{1,2},\xi_{m}]$, $\gamma_{12}^1\in \Lambda_{w_{1,1}w_{1,2}}^{\lambda}$, $w_{1,3}\in \gamma_{12}^1$ and $w_{1,4}\in\beta_0[w_{1,2},\varpi]$ $($see Figure \ref{fig-11-12}$)$. Suppose that the six-tuple $\Omega_1=[w_{1,1},w_{1,2},\gamma_{12}^1,\varpi,w_{1,3},w_{1,4}]$ satisfies Property \ref{Property-A} with constant $\tau_3$,
		\be\label{24-2-15-1}
		k_D(w_{1,1},w_{1,2})=\tau_3,
		\ee and there is an integer $n\geq 2$ such that
		\be\label{24-2-26-10}
		k_D(w_{1,1},w_{1,3})\leq\varrho_{n}.
		\ee   Then the following statements hold:
		\begin{enumerate}
			\item\label{23-11-1} there exist $u_1\in D$, $u_2\in\beta_0[w_{1,2},w_{1,4}]$, $\alpha_{12}\in \Lambda_{u_1u_2}^{\lambda}$, $\varpi_1=w_{1,4}\in \beta_0[u_2,\varpi]$, $u_3\in\alpha_{12}$, and $u_4\in\beta_0[u_2,\varpi_1]$ such that the six-tuple $\Omega=[u_1,u_2,\alpha_{12},\varpi_1,u_3,u_4]$ satisfies Property \ref{Property-A} with constant $\tau_3$; 
			\item\label{23-02-07-1}
			$k_D(u_1,u_2)= \tau_3$,
			$k_D(w_{1,2},u_2)\geq \tau_3-\frac{1}{4}$, and
			\item\label{24-2-26-9}
			$k_D(u_1,u_3)\leq \varrho_{n+1}$.
		\end{enumerate}
	\end{lem}
	
	\bpf To find the required points in the lemma, we need to construct several auxiliary six-tuples based on $\Omega_1$. The construction of the first new six-tuple is included in the following claim.
	
	\begin{figure}[htbp]
		\begin{center}
			\begin{tikzpicture}[scale=2.5]
				\draw (-2,-0.2) to [out=75,in=160] coordinate[pos=0] (xi1) coordinate[pos=0.45] (w12) coordinate[pos=0.85] (w24)  (0,0) to [out=-20,in=210] coordinate[pos=0.09] (w34) coordinate[pos=0.33] (w44) coordinate[pos=0.6] (w14) coordinate[pos=0.75] (varpi) coordinate[pos=1] (xim)node [above, pos=0.88]{$\beta_0$} (2.5,0.2);
				\filldraw  (xi1)node[below] {\small $\xi_1$} circle (0.02);
				\filldraw  (w12)node[below] {\small $w_{2,2}\;\;\;\;$} circle (0.02);
				\node[above,yshift=-0.05cm,rotate=95] at  (w12) {\small $\;\;=$};
				\node[above,xshift=-0.36cm,yshift=0.05cm] at (w12) {\small $w_{1,2}$};
				
				\filldraw  (w24)node[below, rotate=0] {\small $w_{3,2}\;\;\;\;$} circle (0.02);
				\node[above,rotate=90] at  (w24) {\small $\;\;=$};
				\node[above,xshift=-0.36cm,yshift=0.15cm] at (w24) {\small $w_{2,4}$};
				
				\filldraw  (w34)node[below,xshift=0.3cm,yshift=-0.1cm,rotate=-15] {\small $w_{4,2}\!\!=\!{\color{red}u_2}\;\;\;\;$} circle (0.02);
				\node[above,rotate=90] at  (w34) {\small $\;\;=$};
				\node[above,xshift=-0.36cm,yshift=0.15cm] at (w34) {\small $w_{3,4}$};
				
				\filldraw  (w44)node[below] {\small $\;\;\;w_{4,4}\!\!=\!{\color{red}u_4}$} circle (0.02);
				\filldraw  (w14)node[below] {\small $w_{1,4}$} circle (0.02);
				\filldraw  (varpi)node[below] {\small $\;\varpi$} circle (0.02);
				\filldraw  (xim)node[below] {\small $\;\;\;\xi_{m}$} circle (0.02);
				
				\draw (w12) to [out=80,in=180] coordinate[pos=0.4] (w23) coordinate[pos=0.64] (w21)  coordinate[pos=0.86] (w13) coordinate[pos=1] (w11) (1.5,2.5);
				\filldraw  (w23)node[above,rotate=45] {\small $w_{2,3}=w_{3,1}$} circle (0.02);
				\filldraw  (w21)node[above] {\small $w_{2,1}$} circle (0.02);
				\filldraw  (w13)node[above] {\small $w_{1,3}$} circle (0.02);
				\filldraw  (w11)node[above] {\small $\;\;\;\;w_{1,1}$} circle (0.02);
				
				\draw (w24) to [out=85,in=-45] coordinate[pos=0.8] (w33) (w23);
				\filldraw  (w33)node[left] {\small $w_{3\!,3}$} circle (0.02);
				\node[above,rotate=45] at  (w33) {\small $\;\;=$};
				\node[above,xshift=0.65cm,yshift=0.15cm] at (w33) {\small $w_{4,1}\!\!=\!{\color{red}u_1}$};
				
				\draw (w33) to [out=45,in=110]  (0.2,1.3) to [out=-70,in=78] coordinate[pos=0.05] (w43)  (w34);
				\filldraw  (w43)node[left] {\small $w_{4\!,3}$} circle (0.02);
				\node[above,rotate=45] at  (w43) {\small $\;\;=$};
				\node[above,xshift=0.26cm,yshift=0.12cm] at (w43) {\small ${\color{red}u_3}$};
				
				\draw (w43) to [out=45,in=115]  (0.8,1.1) to [out=-65,in=80]  (w44);
				
				\draw (w14) to [out=88,in=-55] (w13);
			\end{tikzpicture}
		\end{center}
		\caption{Illustration for the proof of Lemma \ref{23-08-20}} \label{fig-11-12}
	\end{figure}
	\medskip
	
	\noindent {\bf Claim \ref{23-08-20}-{\rm I}}: {\it
		Let $w_{2,2}=w_{1,2}$. Then there are $w_{2,1}\in\gamma_{12}^{1}[w_{1,3},w_{2,2}]$, $\gamma_{12}^2=\gamma_{12}^{1}[w_{2,1},w_{2,2}]$, $\varpi_2=w_{1,4}$, $w_{2,3}\in \gamma_{12}^{2}$ and $w_{2,4}\in \beta_0[w_{2,2},\varpi_2]$ $($see Figure \ref{fig-11-12}$)$ such that}
\begin{enumerate}
	\item\label{23-08-29-3} {\it
		the six-tuple $\Omega_2=[w_{2,1},w_{2,2},\gamma_{12}^2,\varpi_2,w_{2,3},w_{2,4}]$ satisfies Property \ref{Property-A} with constant $\tau_3$;  }
	\item\label{23-08-29-4} {\it
		$k_D(w_{1,3},w_{2,1})=\tau_2$, and   }
	\item\label{23-08-29-5} {\it
		for any $u\in \gamma_{12}^{2}$ and any $v\in \beta_0[w_{2,4},\varpi_2]$,
		$k_D(u,v)\geq \tau_3-\frac{1}{4}$.}
\end{enumerate}

We first find the point $w_{2,1}$. Note that 
\be\label{8-27-1}
k_D(w_{1,3},w_{2,2})\geq k_D(w_{2,2},w_{1,1})-k_D(w_{1,1},w_{1,3}) \stackrel{\eqref{24-2-15-1}+\eqref{24-2-26-10}}{\geq} \tau_3-\varrho_{n}\stackrel{\eqref{24-2-28-4}}{>}3\tau_2\ee
and so there exists $w_{2,1}\in\gamma_{12}^{1}[w_{1,3},w_{2,2}]$ $($see Figure \ref{fig-11-12}$)$ such that
\be\label{23-08-09-1}k_D(w_{1,3},w_{2,1})=\tau_2.\ee
Then, we know from Lemma \ref{lem-2.2} that
\beq\label{8-26-1}
\begin{aligned}
	k_D(w_{2,1},w_{2,2})  &\leq  \ell_k(\gamma_{12}^{1}[w_{2,1},w_{2,2}] = \ell_k(\gamma_{12}^{1})- \ell_k(\gamma_{12}^{1}[w_{1,3},w_{2,1}])\\  &\leq
	k_D(w_{1,1},w_{1,2})-k_D(w_{1,3},w_{2,1})+\frac{1}{4}\stackrel{\eqref{24-2-15-1}+\eqref{23-08-09-1}}{=}
	\frac{1}{4}-\tau_2 +\tau_3.
\end{aligned}
\eeq

Set $\gamma_{12}^2=\gamma_{12}^1[w_{2,1},w_{2,2}]$ and $\varpi_2=w_{1,4}$. We shall apply Lemma \ref{lem-3.8} to find the remaining two points $w_{2,3}$ and $w_{2,4}$. For this, we need to show that $\tau_3$ is a lower bound for the quantity $k_D(v,w_{1,4})$ for any $v\in \gamma_{12}^{2}$. 

Since the six-tuple $\Omega_1$ satisfies Property \ref{Property-A} with constant $\tau_3$, we see from
Lemma \ref{lem23-01}(\ref{23-01-10-0}) that
\be\label{23-08-20-2}
\begin{aligned}
	k_D(w_{1,4},w_{2,2})&=  k_D(w_{1,2},w_{1,4})\geq k_D(w_{1,2},w_{1,3})-\frac{5}{4}-5\mu_1+\tau_3\\
	&\stackrel{\eqref{8-27-1}}{\geq} -\frac{5}{4}-5\mu_1-\varrho_{n}+2\tau_3.
\end{aligned}
\ee
For any $v\in\gamma_{12}^{2}$, we know from Lemma \ref{lem-2.3} that
\be\label{11-7}
\begin{aligned}k_D(v,w_{1,4})&\geq k_D(w_{1,4},w_{2,2})-k_D(w_{2,2},v)\geq  k_D(w_{1,4},w_{2,2})-\ell_k(\gamma_{12}^{1}[w_{2,2},v])\\ &\geq k_D(w_{1,4},w_{2,2})-\ell_k(\gamma_{12}^{1}[w_{2,2},w_{2,1}])
	\\&\geq k_D(w_{1,4},w_{2,2})-k_D(w_{2,2},w_{2,1})-\frac{1}{2}\\
	&\stackrel{\eqref{8-26-1}+\eqref{23-08-20-2}}{\geq}-2-5\mu_1-\varrho_{n}+\tau_2+\tau_3\stackrel{\eqref{24-2-28-4} }{>}\tau_3
\end{aligned}
\ee

Now, we see from \eqref{8-26-1} and \eqref{11-7} that $w_{2,1}$, $w_{2,2}$, $\gamma_{12}^2$  and $\varpi_2=w_{1,4}$  satisfy all the assumptions in Lemma \ref{lem-3.8}. It follows that
there are $w_{2,3}\in \gamma_{12}^{2}$ and $w_{2,4}\in \beta_0[w_{2,2},w_{1,4}]$ such that
the six-tuple $\Omega_2=[w_{2,1},w_{2,2},\gamma_{12}^2,w_{1,4},w_{2,3},w_{2,4}]$ satisfies Property \ref{Property-A} with constant $\tau_3$ $($see Figure \ref{fig-11-12}$)$.
Then Property \ref{Property-A}(\ref{3-3.2}) guarantees that for any $u\in \gamma_{12}^{2}$ and any $v\in \beta_0[w_{2,4},w_{1,4}]$,
$$k_D(u,v)\geq \tau_3-\frac{1}{4}.$$ This proves the last statement of Claim \ref{23-08-20}-{\rm I}, and hence, the proof of claim is complete.\medskip


The construction of next six-tuple is formulated in the following claim.\medskip

\noindent {\bf Claim \ref{23-08-20}-{\rm II}}: {\it Let $w_{3,1}=w_{2,3}$, $w_{3,2}=w_{2,4}$, $\gamma_{12}^3\in \Lambda_{w_{3,1}w_{3,2}}^\lambda$ and $\varpi_3=w_{1,4}$. Then   }
\begin{enumerate}
	\item\label{24-2-26-1}
	$k_D(w_{3,1},w_{3,2})=\tau_3$;
	\item\label{24-2-26-2}  {\it
		there are $w_{3,3}\in\gamma_{12}^{3}$
		and $w_{3,4}\in\beta_0[w_{3,2},\varpi_3]$ such that   }
	\begin{enumerate}
		\item\label{23-08-29-6}   {\it
			the six-tuple $\Omega_3=[w_{3,1},w_{3,2},\gamma_{12}^3,\varpi_3,w_{3,3},w_{3,4}]$ satisfies Property \ref{Property-A} with constant $\tau_3$;   }
		\item\label{23-08-29-7}   {\it
			$k_D(w_{1,2},w_{3,2})=k_D(w_{2,2},w_{3,2})\geq \tau_3-\frac{1}{4}$, and    }
		\item\label{24-2-26-5}
		$k_D(w_{2,3},w_{3,3})=k_D(w_{3,1},w_{3,3})\leq 39\mu_1+\frac{1}{2}(\varrho_{n}+\tau_2)$.
\end{enumerate}\end{enumerate}

Since $\Omega_2$ satisfies Property \ref{Property-A} with constant $\tau_3$, we infer from Property \ref{Property-A}\eqref{3-3.1} that
\beqq\label{9-1-2}
k_D(w_{3,1},w_{3,2})=k_D(w_{2,3},w_{2,4})=\tau_3.
\eeqq 
This implies statement \eqref{24-2-26-1} of Claim \ref{23-08-20}-{\rm II}.

Moreover, 
we know from Lemma \ref{lem-22-09-30}$(2b)$ that for any $w\in\gamma_{12}^{3}$,
\be\label{24-2-13-3}
k_D(w,w_{1,4})>\tau_3.
\ee
Thus, we know from statement  \eqref{24-2-26-1} of Claim \ref{23-08-20}-{\rm II}, \eqref{24-2-13-3} and Lemma \ref{lem-3.8} (replacing $w_1,w_2,\varpi,\mathfrak{c}$ with $w_{3,1},w_{3,2},w_{1,4},\tau_3$) that there exist $w_{3,3}\in\gamma_{12}^{3}$ and $w_{3,4}\in \beta_0[w_{3,2},w_{1,4}]$, such that
the six-tuple $\Omega_3=[w_{3,1},w_{3,2},\gamma_{12}^3,w_{1,4},w_{3,3},w_{3,4}]$ satisfies Property \ref{Property-A} with constant $\tau_3$ $($see Figure \ref{fig-11-12}$)$. This proves the statement \eqref{23-08-29-6}.

Since $w_{1,2}=w_{2,2}\in \gamma_{12}^2$ and $w_{3,2}=w_{2,4}\in \beta_0[w_{2,4},\varpi_2]$, statement \eqref{23-08-29-7}
follows from Claim \ref{23-08-20}-{\rm I}(\ref{23-08-29-5}).

Since both $\Omega_1$ and $\Omega_2$ satisfy Property \ref{Property-A} with constant $\tau_3$, we know from Lemmas \ref{lem23-01}(\ref{23-12-1-Add}) and Lemma \ref{lem-22-12-1} that
\be\label{24-2-26-7}
k_D(w_{2,1},w_{2,3})\leq 22\mu_1,
\ee
and since
both $\Omega_2$ and $\Omega_3$ satisfy Property \ref{Property-A} with constant $\tau_3$, it follows from Lemma \ref{lem-22-12-0} that
\be\label{24-2-26-6}
k_D(w_{2,3},w_{3,3})=k_D(w_{3,1},w_{3,3})\leq \frac{1}{2}(\tau_3-k_D(w_{2,2},w_{2,3}))+28\mu_1.
\ee
Then \eqref{8-27-1}, the choice of $w_{2,1}$ in \eqref{23-08-09-1} and \eqref{24-2-26-7} show that
\[
\begin{aligned}
	k_D(w_{2,2},w_{2,3})&\geq k_D(w_{2,2},w_{1,3})-k_D(w_{1,3},w_{2,1})-k_D(w_{2,1},w_{2,3})
	\\ &\geq -22\mu_1-\varrho_{n}-\tau_2+\tau_3.
\end{aligned}
\]
Substituting the above estimate in  \eqref{24-2-26-6} gives $$k_D(w_{3,1},w_{3,3})=k_D(w_{2,3},w_{3,3})\leq 39\mu_1+\frac{1}{2}(\varrho_{n}+\tau_2).$$
This proves the statement \eqref{24-2-26-5}, and hence, the proof of claim is complete.\medskip

One more six-tuple is constructed in the following claim.\medskip

\noindent {\bf Claim \ref{23-08-20}-{\rm III}}: {\it Let $w_{4,1}=w_{3,3}$, $w_{4,2}=w_{3,4}$, $\gamma_{12}^4\in \Lambda_{w_{4,1}w_{4,2}}^\lambda$ and $\varpi_4=w_{1,4}$. Then   }
\begin{enumerate}
	\item\label{24-2-26-3}
	$k_D(w_{4,1},w_{4,2})=\tau_3$;
	\item\label{24-2-26-4}  {\it
		there are $w_{4,3}\in\gamma_{14}^{2}$
		and $w_{4,4}\in\beta_0[w_{4,2},\varpi_4]$ such that   }
	\begin{enumerate}
		\item\label{23-12-11-1}  {\it
			the six-tuple $\Omega_4=[w_{4,1},w_{4,2},\gamma_{12}^4,\varpi_4,w_{4,3},w_{4,4}]$ satisfies Property \ref{Property-A} with constant $\tau_3$;    }
		\item\label{24-2-16-3}   {\it
			$k_D(w_{1,2},w_{4,2})=k_D(w_{2,2},w_{4,2})\geq \tau_3-\frac{1}{4}$, and  }
		\item\label{24-2-26-8}
		$k_D(w_{4,1},w_{4,3})=k_D(w_{4,1},w_{4,3})\leq 48\mu_1+\frac{1}{4}(\varrho_{n}+\tau_2)$.
\end{enumerate}\end{enumerate}

We have shown that all six-tuples $\Omega_1,$ $\Omega_2$ and $\Omega_3$ satisfy Property \ref{Property-A}. Thus, statement \eqref{24-2-26-3} of Claim \ref{23-08-20}-{\rm III} follows from Property \ref{Property-A}(\ref{3-3.1}).
Moreover, as Claim \ref{23-08-20}-{\rm II}\eqref{23-08-29-7} shows
$$k_D(w_{3,1},w_{3,3})=k_D(w_{2,3},w_{3,3})\leq 39\mu_1+\frac{1}{2}(\varrho_{n}+\tau_2)\stackrel{\eqref{24-2-28-4}}{\leq}-56\mu_1-\frac{\varrho_{n}+\tau_2}{2}+\tau_2.$$
we know from Lemma \ref{lem-22-09-30}$(3)$ that for any $w\in\gamma_{12}^{4}$,
\be\label{24-5-20-1}
k_D(w,w_{1,4})>\tau_3.
\ee
Thus we know from Lemma \ref{lem-3.8} (replacing $w_1,w_2,\varpi,\mathfrak{c}$ with $w_{4,1},w_{4,2},w_{1,4},\tau_3$), together with the statement \eqref{24-2-26-3} of Claim \ref{23-08-20}-{\rm III} and \eqref{24-5-20-1}, that there exist $w_{4,3}\in\gamma_{12}^{4}$ and $w_{4,4}\in \beta_0[w_{4,2},\varpi_4]$, such that
the six-tuple $\Omega_4=[w_{4,1},w_{4,2},\gamma_{12}^4,\varpi_4,w_{4,3},w_{4,4}]$ satisfies Property \ref{Property-A} with constant $\tau_3$ $($see Figure \ref{fig-11-12}$)$. This proves the statement \eqref{23-12-11-1}.

Since $w_{2,2}=w_{2,1}$ and $w_{4,2}=w_{3,4}\in \beta_0[w_{2,4},w_{1,4}]$, we see from Claim \ref{23-08-20}-{\rm I}(\ref{23-08-29-5}) that statement \eqref{24-2-16-3} holds.

Since both $\Omega_3$ and $\Omega_4$ satisfy Property \ref{Property-A}, we obtain from Lemma \ref{lem-22-12-0} (with $\Omega_1=\Omega_3$ and $\Omega_2=\Omega_4$) and the above estimate that
\beqq
\begin{aligned}
	k_D(w_{3,3},w_{4,3})&=k_D(w_{4,1},w_{4,3})\leq \frac{1}{2}\big(\tau_3-k_D(w_{3,2},w_{3,3})\big)+28\mu_1\\
	&\leq \frac{1}{2}k_D(w_{2,3},w_{3,3})+28\mu_1\stackrel{{\rm Claim\ \ref{23-08-20}-}{\rm II}(\ref{24-2-26-8})}{\leq} 48\mu_1+\frac{1}{4}(\varrho_{n}+\tau_2).
\end{aligned}
\eeqq
Since $w_{3,3}=w_{4,1}$, the above inequality shows that statement \eqref{24-2-26-8} holds, and hence, the proof of claim is complete.
\medskip

Let $u_1=w_{4,1}$, $u_2=w_{4,2}$, $\alpha_{12}=\gamma_{12}^4$, $\varpi_1=w_{1,4}$, $u_3=w_{4,3}$ and $u_4=w_{4,4}$. Then stataments (1) and (2) in Lemma \ref{23-08-20} follow from Claim \ref{23-08-20}-{\rm III} \eqref{24-2-26-3}, \eqref{23-12-11-1} and \eqref{24-2-16-3} . Since by \eqref{24-4-20-1}, $45\mu_1+\frac{1}{4}\varrho_{n}+\frac{1}{4}\tau_2\leq \varrho_{n+1}$ when $n\geq 2$,
statement (3) in Lemma \ref{23-08-20} follows from Claim \ref{23-08-20}-{\rm III} \eqref{24-2-26-8}. The proof of lemma is thus complete.
\epf

\subsection{Final proof of Theorem \ref{22-11-16}}\label{subsec:final proof of auxiliary theorem}

\begin{figure}[htbp]
	\begin{center}
		\begin{tikzpicture}[scale=1]
			\draw (-6,-0.5) to [out=65,in=155] (0,0) to [out=-25,in=205] (5.6,0.28);
			\node[below] at(-1.6,0.3) {$\alpha$};
			\filldraw  (-6,-0.5)node[below] {\tiny $\zeta_1=x$} circle (0.04);
			\filldraw  (5.6,0.28)node[below right,xshift=-0.3cm] {\tiny $\zeta_{m}=y_1$} circle (0.04);
			
			\draw (-6,-0.5) to [out=95,in=180] node[left,pos=0.35,xshift=0cm] { $\gamma_1$} (-0.4,9) to [out=0,in=95] (5.6,0.28);
			\filldraw  (-5.985,2)node[below right] {\tiny $\xi_1=w_{1,2}$} circle (0.04);
			\filldraw  (-1,8.97)node[above] {\tiny $w_{1,3}\;\;\;\;$} circle (0.04);
			\filldraw  (0,9)node[above] {\tiny $w_{1,1}$} circle (0.04);
			\filldraw  (1.5,8.61)node[above] {\tiny $\;\;x_0$} circle (0.04);
			\filldraw  (5.325,2.3)node[above left] {\tiny $\xi_{m}$} circle (0.04);
			\node[above,xshift=0.3cm] at(3.4,7.5) {$\gamma_1\in \Lambda_{xy_{1}}^{\lambda}$};
			
			\draw (-5.985,2) to [out=60,in=155] (0,2) to [out=-25,in=210] (5.32,2.3);
			\node at(-1.6,2) {$\beta_0$};
			\filldraw  (-4.7,3.05)node[below] {\tiny $w_{2,2}$} circle (0.04);
			\filldraw  (-3.5,3.14)node[below] {\tiny $w_{3,2}$} circle (0.04);
			\filldraw  (-2.3,2.905)node[below] {\tiny $w_{m_1,2}\;\;$} circle (0.04);
			\filldraw  (-0.5,2.23)node[below] {\tiny $w_{m_1,4}\;\;\;$} circle (0.04);
			\filldraw  (0.7,1.72)node[below] {\tiny $w_{3,4}$} circle (0.04);
			\filldraw  (1.9,1.45)node[below] {\tiny $w_{2,4}$} circle (0.04);
			\filldraw  (3.7,1.6)node[below] {\tiny $w_{1,4}$} circle (0.04);

			\draw (-4.7,3.05) to [out=85,in=180]node[left,pos=0.1,xshift=0.13cm] {\tiny $\gamma_{12}^{2}$} (-0.5,8);
			\filldraw  (-0.5,8)node[below] {\tiny $\;\;\;w_{2,1}$} circle (0.04);
			\filldraw  (-1.2,7.94)node[above] {\tiny $w_{2,3}\;\;\;\;$} circle (0.04);
			
			\draw (1.9,1.45) to [out=90,in=-20]node[left,pos=0.1,xshift=0.15cm] {\tiny $\gamma_{34}^{2}$} (-1.2,7.94);
			
			\draw (-3.5,3.14) to [out=85,in=200]node[left,pos=0.15,xshift=0.13cm] {\tiny $\gamma_{{12}}^{3}$} (-0.8,7);
			\filldraw  (-0.8,7)node[below] {\tiny $\;\;\;\;\;w_{3,1}$} circle (0.04);
			\filldraw  (-1.3,6.8)node[above] {\tiny $w_{3,3}\;\;\;\;\;\;$} circle (0.04);
			
			\draw (0.7,1.72) to [out=85,in=-40]node[left,pos=0.15,xshift=0.15cm] {\tiny $\gamma_{34}^{3}$} (-1.3,6.8);
			
			\node[rotate=84] at(-1.4,6) {$\cdots$};
			
			\draw (-2.3,2.905) to [out=85,in=220] node[pos=0.2,rotate=0,xshift=-0.225cm] {\tiny $\gamma_{12}^{m_1}$}(-1.2,5.5);
			\filldraw  (-1.2,5.5)node[right] {\tiny $w_{m_1,1}$} circle (0.04);
			\filldraw  (-1.5,5.2)node[left] {\tiny $w_{m_1,3}$} circle (0.04);
			
			\draw (-0.5,2.23) to [out=85,in=-40]node[left,pos=0.23,rotate=0,xshift=0.15cm] {\tiny $\gamma_{34}^{m_1}$} (-1.5,5.2);
			
			\draw (3.7,1.6) to [out=100,in=-20]node[left,pos=0.07,xshift=0.15cm] {\tiny $\gamma_{34}^{1}$} (-1,8.97);
		\end{tikzpicture}
	\end{center}
	\caption{Illustration for the proof of Theorem \ref{22-11-16}} \label{fig-8-11}
\end{figure}

In this section, we complete the proof of  Theorem \ref{22-11-16}. Towards a contradiction, set
\be\label{24-2-28-2}
\tau_2=\frac{1}{64}\log \mu_5\;\;\mbox{and}\;\;\tau_3=\frac{1}{16}\log \mu_5.
\ee

For $i\in \{1, \ldots, m\}$, let $\xi_i$ be the points defined in Section \ref{sub-5.3}, and $w_{1,2}=\xi_1$ (see Figure \ref{fig-8-11}).
Since Lemma \ref{cl3.2-1}$(i)$ and the choice of $\xi_i$ in \eqref{2022-09-30-3} ensure that
\be\label{23-08-20-1}
k_D(\xi_{i},x_0)= k_D(w_{1,2},x_0)\geq \log \frac{\mu_5-2}{2+\mu_6}-\frac{1}{2}\stackrel{\eqref{24-2-28-2}}{>}2\tau_3+\frac{1}{2},
\ee
we infer from the above estimate that there exists $w_{1,1}\in\gamma_1[w_{1,2},x_0]$ such that
\be\label{2022-10-12-2}
k_D(w_{1,1},w_{1,2})=\tau_3.
\ee

Let  $\gamma_{12}^1=\gamma_1[w_{1,1},w_{1,2}]$. Then for any $w\in\gamma_{12}^1$, we know from Lemma \ref{lem-2.3} that
\beqq
k_D(w,\xi_{m}) \geq
\ell_k(\gamma_1[w,\xi_{m}])-\frac{1}{2} >\ell_k(\gamma_1[\xi_{m},x_0])-\frac{1}{2}
\geq k_D(\xi_{m},x_0)-\frac{1}{2}\stackrel{\eqref{23-08-20-1}}{>}2\tau_3.
\eeqq
Thus it follows from \eqref{2022-10-12-2} and Lemma \ref{lem-3.8} that there exist $w_{1,3}\in \gamma_{12}^1$ and $w_{1,4}\in \beta_0[w_{1,2},\xi_{m}]$ such that for $\varpi=w_{1,4}$, the six-tuple $[w_{1,1},w_{1,2},\gamma_{12},\varpi,w_{1,3},w_{1,4}]$ satisfies Property \ref{Property-A} with constant $\tau_3$ (see Figure \ref{fig-8-11}).
Then by Lemmas \ref{lem23-01}(\ref{23-12-1-Add}) and Lemma \ref{lem-22-12-1}, it holds
\beqq\label{23-08-08-1}k_D(w_{1,1},w_{1,3})\leq 22\mu_1\stackrel{\eqref{24-4-20-1}}{\leq}\varrho_2,
\eeqq
and thus, we infer from (\ref{2022-10-12-2}) and Lemma \ref{23-08-20} that
there exist $w_{2,1}\in D$, $w_{2,2}\in\beta_0[w_{1,2},w_{1,4}]$, $\gamma_{12}^2\in \Lambda_{w_{2,1},w_{2,2}}^\lambda$, $\varpi_2=w_{1,4}\in \beta_0[w_{2,2},w_{1,4}]$, $w_{2,3}\in \gamma_{12}^{2}$ and $w_{2,4}\in \beta_0[w_{2,2},\varpi_2]$ such that the six-tuple $[w_{2,1},w_{2,2},\gamma_{12}^2,\varpi_2,w_{2,3},w_{2,4}]$ satisfies Property \ref{Property-A} with constant $\tau_3$ (see Figure \ref{fig-8-11}). Moreover, it holds
\beqq\label{9-3-1}
k_D(w_{2,1},w_{2,2})=\tau_3, \;\; k_D(w_{1,2},w_{2,2})\geq \tau_3-\frac{1}{4} \;\;\mbox{and}\;\;  k_D(w_{2,1},w_{2,3})\leq \varrho_3.
\eeqq

Now, applying Lemma \ref{23-08-20} once again, we know that there exist $w_{3,1}\in D$, $w_{3,2}\in\beta_0[w_{2,2},w_{2,4}]$, $\gamma_{12}^3\in \Lambda_{w_{3,1}w_{3,2}}^\lambda$, $\varpi_3=w_{2,4}\in \beta_0[w_{3,2},w_{2,4}]$, $w_{3,3}\in \gamma_{12}^{3}$ and $w_{3,4}\in \beta_0[w_{3,2},\varphi_3]$ such that the six-tuple $[w_{3,1},w_{3,2},\gamma_{12}^3, \varphi_3,w_{3,3},w_{3,4}]$ satisfies Property \ref{Property-A} with constant $\tau_3$ (see Figure \ref{fig-8-11}). Moreover, it holds
\beqq\label{9-3-3}
k_D(w_{3,1},w_{3,2})=\tau_3, \;\; k_D(w_{2,2},w_{3,2})\geq \tau_3-\frac{1}{4} \;\;\mbox{and}\;\;  k_D(w_{3,1},w_{3,3})\leq \varrho_4.
\eeqq

Repeating this procedure for $m_1$ times, where $m_1=\big[\frac{\ell_k(\beta_0)}{\tau_3-\frac{1}{4}}\big]+3$, we get six sequences $\{w_{i,1}\}_{i=1}^{m_1}$, $\{w_{i,2}\}_{i=1}^{m_1}$, $\{\gamma_{12}^i\}_{i=1}^{m_1}$, $\{\varphi_i\}_{i=1}^{m_1}$, $\{w_{i,3}\}_{i=1}^{m_1}$ and $\{w_{i,4}\}_{i=1}^{m_1}$ (see Figure \ref{fig-8-11}) such that
\begin{enumerate}
	\item $w_{1,1}\in\gamma_1[x_0,\xi_1]$, $w_{1,2}=\xi_1$, $w_{1,3}\in \gamma_1[w_{1,1},w_{1,2}]$ and $w_{1,4}\in \beta_0[w_{1,2},\xi_{m}]$;
	\item for each $i\in\{2,\cdots,m_1\}$, $w_{i,1}\in D$, $w_{i,2}\in\beta_0[w_{i-1,2},w_{i-1,4}]$, $\gamma_{12}^i\in \Lambda_{w_{i,1}w_{i,2}}^\lambda$, $\varpi_i=w_{i-1,4}$, $w_{i,3}\in \gamma_{12}^{i-1}$ and $w_{i,4}\in \beta_0[w_{i,2},\varpi_i]$;
	\item  for each $i\in\{1,\cdots,m_1\}$, the six-tuple $[w_{i,1},w_{i,2},\gamma_{12}^i,\varpi_i,w_{i,3},w_{i,4}]$ satisfies Property \ref{Property-A} with constant $\tau_3$, where $w_{0,4}=\xi_m$, and
	\item\label{23-08-23-1} for each $i\in\{1,\cdots, m_1-1\}$,
	$k_D(w_{i,2},w_{i+1,2})\geq \tau_3-\frac{1}{4}$.
\end{enumerate}
Then summing over all $i$ in \eqref{23-08-23-1} gives
$$\sum_{i=1}^{m_1-1}k_D(w_{i,2},w_{i+1,2})\geq \Big(\tau_3-\frac{1}{4}\Big)(m_1-1)\geq\ell_k(\beta_0)+\tau_3-\frac{1}{4},$$
which is impossible as
$$\sum_{i=1}^{m_1-1}k_D(w_{i,2},w_{i+1,2})\leq \ell_k(\beta_0).$$
The proof of  Theorem \ref{22-11-16} is thus complete.
\medskip

\textbf{Conflict of interest}. On behalf of all authors, the corresponding author states that there is no conflict of interest.

\smallskip 
\textbf{Data availability}. There is no research data associated with this article.

\medskip
\textbf{Acknowledgment.}
C.-Y. Guo is supported by the Young Scientist Program of the Ministry of Science and Technology of China (No.~2021YFA1002200), the NSF of China (No.~12311530037), the Taishan Scholar Project, the NSF of Shandong Province (No.~ZR2022YQ01) and the Jiangsu Provincial Scientific Research Center of Applied Mathematics (No.~BK20233002). M. Huang and X. Wang are partly supported by NSF of China (No.12371071 and No.12571081).

We would like to thank Prof.~Pekka Koskela for helpful discussions and suggestions that greatly improved our exposition.

\end{document}